\documentclass[a4paper]{article}

\usepackage[plainpages=false, colorlinks=true, 
            linkcolor=black, urlcolor=black, citecolor=black]{hyperref}
\usepackage{geometry}

\usepackage{tikz}
\usepackage{overpic}
\usepackage{amsmath, amssymb, amsthm,  amsfonts, mathrsfs, cite}
\usepackage{multirow}
\usepackage{color}
\usepackage{booktabs}
\usepackage{tikz}
\usepackage{bm}
\usepackage{hyperref}
\usepackage{scrextend}
\usepackage{arydshln}
\usepackage {url}
\usepackage[font=scriptsize, labelfont=bf]{caption}

\usepackage{graphicx}
\usepackage{epstopdf}

\usepackage{tikz-cd}

\usepackage{algpseudocode}

\usepackage{array,  tabularx}
\usepackage{algorithm, algcompatible}
\usepackage{booktabs} 
\usepackage{paralist} 
\usepackage{verbatim} 
\usepackage{subfig} 
\usepackage{xcolor, marginnote, enumitem}
\usepackage[numbers]{natbib}
\numberwithin{equation}{section}
\theoremstyle{theorem}
\newtheorem{lemma}{Lemma}
\newtheorem{theorem}{Theorem}
\newtheorem{proposition}{Proposition}

\theoremstyle{remark}
\newtheorem{remark}{Remark}

\theoremstyle{definition}

\DeclareMathOperator{\dom}{dom}

\DeclareMathOperator{\diag}{Diag}

\DeclareMathOperator{\test}{test}

\DeclareMathOperator{\TGV}{TGV}

\newcommand{\mE}{\mathbb{E}}

\newcommand{\lookUp}[1]{}

\newcommand{\mF}{\mathcal{F}}
\newcommand{\mG}{\mathcal{G}}
\newcommand{\mH}{\mathcal{H}}

\newcommand{\mL}{\mathcal{L}}

\newcommand{\mI}{\mathcal{I}}
\newcommand{\mT}{\mathcal{T}}

\newcommand{\mK}{\mathcal{K}}

\newcommand{\mwy}{\mathbf{w}_y}
\newcommand{\mrx}{\mathbf{R}_x} 
\newcommand{\mry}{\mathbf{R}_y} 
\newcommand{\mrrx}{\mathbf{r}_x} 
\newcommand{\mrry}{\mathbf{r}_y} 
\newcommand{\mv}{\mathbf{v}}
\newcommand{\mty}{\tilde{\mathbf{y}}}
\newcommand{\my}{\mathbf{y}}
\newcommand{\tbx}{\tilde {\bar{x}}}
\newcommand{\tby}{\tilde {\bar{y}}}

\newcommand{\bitem}{\begin{itemize}}
	\newcommand{\eitem}{\end{itemize}}

\newcommand{\bpm}{\begin{pmatrix}}

\DeclareMathAlphabet{\mathbfit}{OML}{cmm}{b}{it}

\usepackage{xspace}
\usepackage{bold-extra}
\usepackage[most]{tcolorbox}

\colorlet{texcscolor}{blue!50!black}
\colorlet{texemcolor}{red!70!black}
\colorlet{texpreamble}{red!70!black}
\colorlet{codebackground}{black!25!white!25}

\begin{document}
	\title{A stochastic preconditioned Douglas-Rachford splitting method for saddle-point problems}
	\author{Yakun Dong\thanks{	Institute for Mathematical Sciences, Renmin University of China. Current address: Faculty of Mathematics, University of Vienna, Oskar Morgenstern-Platz 1, 1090 Vienna, Austria. \href{mailto:yakun.dong@univie.ac.at}{yakun.dong@univie.ac.at}}
		\and Kristian Bredies\thanks{Institute for Mathematics and Scientific Computing, University of Graz, Heinrichstraße 36, A-8010 Graz,
			Austria. \href{mailto:kristian.bredies@uni-graz.at}{kristian.bredies@uni-graz.at} }
		\and Hongpeng Sun\thanks{Institute for Mathematical Sciences,
			Renmin University of China, 100872 Beijing, China. \href{mailto:hpsun@amss.ac.cn}{hpsun@amss.ac.cn}}
	}

\maketitle

\begin{abstract}
In this article, we propose and study a stochastic and relaxed preconditioned Douglas--Rachford splitting method to solve saddle-point problems that have separable dual
variables. We prove the almost sure convergence of the iteration sequences in Hilbert spaces for a class of convex-concave and nonsmooth saddle-point problems. We also provide the sublinear convergence rate for the ergodic sequence concerning the expectation of the restricted primal-dual gap functions. Numerical experiments show the high efficiency of the proposed stochastic and relaxed preconditioned Douglas--Rachford splitting methods.  
\end{abstract}
\paragraph{Key words.}{Douglas--Rachford splitting, stochastic
optimization, saddle-point problems, over-relaxation, almost sure convergence, preconditioning techniques.}
\paragraph{MSCodes.}{65K10,
65F08,
49K35,
90C25.
}

\section{Introduction}
The Douglas--Rachford (henceforth abbreviated as DR) splitting method is an iterative splitting method for finding a zero point of the sum of maximal set-valued monotone operators, 
which is a powerful optimization tool for convex optimization and has plenty of applications \cite{HBPL2, BSCC, Condat1, he1, EP}. The method originated from splitting linear operators for solving heat conduction equations \cite{DR} and was extended to a more general nonlinear maximal monotone operator splitting framework \cite{LM, BF}. The convergence of the DR splitting method in Hilbert spaces can be found in \cite{EP} through a proximal point framework. We further refer to \cite{arias, he2, he3,guohan} for general DR splitting methods for convex problems and \cite{Artacho, luke,patrinos, tkpong} for nonconvex problems.   The DR splitting method is essentially an implicit and unconditionally stable splitting algorithm, which means the DR splitting method is convergent with arbitrary positive step sizes.
Solving implicit equations including linear equations can be computationally demanding for DR methods. In \cite{BSCC}, within the generalized proximal-point framework, an efficient preconditioned DR framework with a convergence guarantee aiming at a specific class of convex-concave saddle-point problems is developed.  It can incorporate the classical preconditioning techniques for linear systems into the nonlinear DR iterations.
We refer to \cite{BS} for applications of preconditioned DR and \cite{BS1} for its accelerated variants. 

Due to a wide range of applications in various fields including big data science and machine learning, there are a lot of studies on stochastic first-order algorithms \cite{DPB,CPE} and especially for saddle-point problems \cite{AFC,CERS,CPE}. In this paper, we study a stochastic preconditioned DR splitting framework which includes the stochastic DR without preconditioners \cite{CPE} as a special case for solving the following saddle-point problems with separable dual variables \cite{AFC, CPE}. We refer to \cite{hanyan} for the randomized DR methods aiming at solving linear systems.
\begin{equation}\label{eq:saddle-point-sum}
\min_{x\in  \dom{\mF}}
\max_{y\in \dom{\mG}} \mathcal{L}(x,y), \quad 
\mathcal{L}(x,y)   :=\mathcal{F}(x) + \sum_{i=1}^n
(\langle \mathcal{K}_i x,y_i\rangle-\mathcal{G}_i(y_i)).
\end{equation}
where $\mathcal{F}: \mathbb{X}  \rightarrow \mathbb{R}_{\infty}:=\mathbb{R} \cup\{+\infty\}$ and each $\mG_i:  \mathbb{Y}_i  \rightarrow \mathbb{R}_{\infty}$ are proper, convex and lower semi-continuous. Here, $\mathbb{X}$ and $ \mathbb{Y}_i$ are real Hilbert spaces and each $\mK_i : \mathbb{X}  \rightarrow \mathbb{Y}_i$ is linear and continuous.  The space $\mathbb{Y} $ is a product space, i.e.,  $\mathbb{Y}=\prod^n_{i=1} \mathbb{Y}_i$ with $ y=(y_1,y_2,\ldots,y_n)^T$.
The saddle-point problem \eqref{eq:saddle-point-sum} can be reformulated as the following common primal-dual formulation \cite{HBPL2}: 
\begin{equation}\label{eq:saddle-point}
\min_{x\in  \dom{\mF}}
\max_{y\in  \dom{\mG}}\mathcal{L}(x,y), \quad 
\mathcal{L}(x,y):= \mathcal{F}(x)+
\langle \mathcal{K}x,y\rangle-\mathcal{G}(y),
\end{equation}
where the dual function $\mG $ is separable, i.e.,
$\mG(y) = \sum^n_{i=1} \mG_i(y_i)$ with  $ \mathcal{G}: \mathbb{Y}  \rightarrow \mathbb{R}_{\infty}$ and $\mK = [\mK_1, \mK_2, \ldots, \mK_n]^T:  \mathbb{X}  \rightarrow \mathbb{Y}$. When strong duality holds \cite{HBPL2}, the solutions of \eqref{eq:saddle-point-sum} are exactly the solutions of the following primal and dual problems:
\begin{equation}\label{eq:primal-dual-problem}
	\min _{x \in \mathbb{X}}
	\mathcal{P}(x), \   \mathcal{P}(x) := \mathcal{F}(x)+ \sum_{i=1}^n\mathcal{G}_i^{*}(\mathcal{K}_i x), \quad  \max _{y \in \mathbb{Y}}\mathcal{D}(y), \ \mathcal{D}(y):=-\mathcal{F}^{*} (-\mathcal{K}^{*} y )-\sum_{i=1}^n\mathcal{G}_i(y_i),
\end{equation}
with $\mathcal{F}^{*}$ and $\mathcal{G}_i^{*}$ denoting the Fenchel conjugate functionals of $\mathcal{F}$ and $\mathcal{G}_i$, respectively.



Compared to related work, our contributions are as follows.
{First, different from the stochastic first-order primal-dual method \cite{CERS} or stochastic classical DR splitting method without preconditioners \cite{CPE},  preconditioners are introduced. The challenges in the convergence analysis are addressed for the proposed stochastic DR splitting method. Preconditioners already introduced challenges for the deterministic preconditioned DR and extra effort is necessary for the corresponding convergence analysis \cite{BSCC, BS1}. Second, we incorporate a general over-relaxed strategy into the stochastic and preconditioned  DR splitting method with different relaxation parameters on different variables. It is well known that over-relaxation techniques can accelerate the deterministic DR splitting algorithm \cite{EP,Condat1,luke, BrediesELN}. For the over-relaxed, stochastic, and preconditioned DR (hereafter abbreviated as SRPDR)},
the almost sure convergence framework of the stochastic fixed point iteration based on firmly nonexpansive operators developed in \cite{CPE} cannot be applied directly here.  We circumvent these difficulties by analyzing the relation between the stochastic updates and the one-step deterministic updates based on the current stochastic iteration.
We mainly employed the celebrated stochastic quasi-Fejér monotonicity \cite{CPE,CPE1, YMV} and Opial's lemma \cite{Op1}.
{Third}, we show the {SRPDR} splitting method has a $\mathcal{O}(1/K)$ rate for the ergodic sequence respecting restricted primal-dual gap functions where $K$ is the iteration number.
Since taking the expectation of the supremum of the primal-dual gap {concerning transitional variables} is a very subtle issue \cite{AFC}, a detailed analysis is devoted to this part. {Furthermore, we also show the SRPDR splitting method has a $\mathcal{O}(1/K)$ rate for primal error when Lipschitz continuity is satisfied for dual functions.} Finally, we show the high efficiency of the proposed {SRPDR} compared to state-of-art methods including the stochastic primal-dual method \cite{CERS} with some numerical experiments {on synthetic and real datasets}. 

The rest of the paper is organized as follows. In section \ref{sec:saddle-pdr}, we give a brief introduction of the preconditioned deterministic DR and {the proposed stochastic, over-relaxed, and preconditioned DR} splitting method along with some necessary analysis tools. In section \ref{sec:almost:sure}, we discuss the almost sure convergence of the proposed {SRPDR} splitting method to a fixed point using stochastic quasi-Fejér monotonicity \cite{CPE} and Opial's lemma \cite{Op1}. In section \ref{sec:sublinear:gap}, we show the sublinear convergence rate of ergodic sequences respecting the expectation of the restricted primal-dual gap {and the primal error}. In section \ref{sec:sto:pdrq}, we extend these results to the stochastic {relaxed} and preconditioned DR method for quadratic-linear problems (henceforth abbreviated as {SRPDRQ}). In section \ref{sec:numer}, we present detailed numerical experiments for total generalized variation (TGV) regularized problems and {binary classification with real datasets} to show the competitive performance of the proposed stochastic and preconditioned DR methods. Finally, we give some conclusions in section \ref{sec:conclusion}.

\section{Preconditioned DR framework}\label{sec:saddle-pdr}

Now, let us turn to the deterministic preconditioned Douglas--Rachford splitting method (PDR), which aims at solving the saddle-point problem \eqref{eq:saddle-point}. It reads as follows \cite{BSCC, BS, BS1}:
\begin{equation}\label{PDR}\tag {PDR}
\left\{\begin{array}{ll}
x^{k+1}=N_1^{-1}[(N_1-I)x^k+\bar x^k-\sigma\mathcal{K}^\ast   y^{k+1}],   \\
y^{k+1}=N_2^{-1}[(N_2-I)  y^k+ \bar y^k+\tau\mathcal{K} x^{k+1}],   \\
x_{\test}^{k+1}=(I+\sigma\partial\mathcal{F})^{-1}[2x^{k+1}-\bar x^k],\\
y_{\test}^{k+1}=(I+\tau\partial\mathcal{G})^{-1}[2  y^{k+1}- \bar y^k],\\
\bar x^{k+1}=\bar x^k+x_{\test}^{k+1}-x^{k+1},\\
\bar y^{k+1}= \bar y^k+  y_{\test}^{k+1}-  y^{k+1},
\end{array}
\right.
\end{equation}
where $\sigma $ and $\tau$ are arbitrary positive primal and dual step sizes, $N_1 \succeq  I $, $N_2\succeq  I$ with $N_1:  \mathbb{X}  \rightarrow  \mathbb{X} $ and $N_2:  \mathbb{Y}  \rightarrow  \mathbb{Y} $ being bounded, linear and self-adjoint operators, and $x_{\test}^{k+1}$ and $y_{\test}^{k+1}$ are transitional variables.
Denoting by 
$z = (x,y) \in \mathbb{X} \times \mathbb{Y}$ and $\bar z = (\bar x, \bar y) \in \mathbb{X} \times \mathbb{Y}$,
the iteration \ref{PDR} with $N_1=I$, $N_2=I$ and $\sigma=\tau$ becomes the original DR splitting method \cite{EP, LM}. It was based on the following maximal monotone operator splitting
\[
0 \in  Az + Bz, \ \  \text{with}\quad A = \begin{pmatrix}
0 & \mK^* \\ -\mK & 0
\end{pmatrix}, 
\quad
B =
\begin{pmatrix}
\partial \mF & 0 \\ 0 & \partial \mG
\end{pmatrix},
\]
by the first-order optimality conditions of \eqref{eq:saddle-point}.
Note that the first two equations of \ref{PDR} are implicit. 
However, if we choose $N_2 =I$ and an operator $M$ such that
\begin{equation}\label{N1:m}
N_{1} = M - \sigma \tau \mK^{*}\mK
\end{equation}
where $M: X \rightarrow X $ is linear, continuous, self-adjoint and satisfies $M\succeq T := I + \sigma \tau  \mK^*\mK$, we can update $x^{k+1}$ first with the following classical preconditioned iteration 
\begin{equation}\label{preconditioner:from}
x^{k+1} = x^{k} + M^{-1}(b^{k} - Tx^{k}), \quad b^{k} = \bar x^k - \sigma \mK^*\bar y^k,
\end{equation}
and then obtain the linear update $y^{k+1}=\bar y^k + \tau \mK x^{k+1}$ according to the second equation of \ref{PDR} followed by $x^{k+1}$ \cite{BSCC, BS, BS1}. Classical preconditioning techniques can be employed with \eqref{preconditioner:from} which can be seen as an iterative method for solving $Tx^{k+1}=b^k$. This linear equation arises from the linear update step in the original Douglas--Rachford splitting method where $N_1=I$ and $N_2=I$. Preconditioned DR splitting can bring high efficiency with great flexibility and convergence guarantees using finitely many iterations of any feasible preconditioner \cite{BSCC, BS, BS1}. A preconditioner is called \emph{feasible} if $M \succeq  T$ as in \eqref{preconditioner:from} includes a lot of classical efficient preconditioners \cite{BSCC}. Let us denote by $u:=(x, y,\bar x, \bar y)\in \mathbb{U}: =(\mathbb{X}\times \mathbb{Y})^2$ and introduce an associated positive semidefinite and possibly degenerate inner product as
\begin{equation}
\langle u, u^{\prime} \rangle_{\mathbf{M}}:= \langle \mathbf{M} u,u^{\prime} \rangle =\frac{1}{\sigma} [ \langle (N_{1}-I ) x, x^{\prime} \rangle+ \langle\bar{x}, \bar{x}^{\prime} \rangle ]+\frac{1}{\tau} [ \langle (N_{2}-I ) y, y^{\prime} \rangle+ \langle\bar{y}, \bar{y}^{\prime} \rangle ],
\end{equation}
with  $\mathbf{M}:=\diag[\tfrac1\sigma (N_1-I), \tfrac1\tau (N_2-I), \tfrac{1}{\sigma} I, \tfrac{1}{\tau} I]$. 
Denoting by $\mT : \mathbb{U} \rightarrow \mathbb{U}$ as the iteration operator that maps $u^k=(x^k,y^k,\bar x^k, \bar y^k)$ to $u^{k+1}=(x^{k+1},y^{k+1},\bar x^{k+1}, \bar y^{k+1})$ according to  \ref{PDR}, we can rewrite the iteration \ref{PDR} as 
\begin{equation}\label{eq:pdr:T}
u^{k+1}: = \mT u^k.
\end{equation}

Now, we present the following lemma and proposition, which are useful for analyzing the stochastic PDR method in the subsequent sections. The proofs for the case where \(\sigma = \tau\) can be found in \cite{BSCC, BS}, while the proofs for the case where \(\sigma \neq \tau\) can be obtained easily.
\begin{lemma}[{\cite[Lemma 2.1]{BS}}]\label{lem:inner-product-1}
	For ${u^k={(x^k, y^k,\bar x^k,\bar y^k)}}$, we have
	\begin{equation}
	\frac{1}{\sigma}\langle x^{k+1}-\bar x^{k+1},x_{\test}^{k+1}\rangle +\frac{1}{\tau}\langle y^{k+1}-\bar y^{k+1}, y_{\test}^{k+1}\rangle
	=
	\langle 	u^k-u^{k+1},u^{k+1}\rangle_\mathbf{M}.
	\end{equation}
\end{lemma}
Hence, employing Lemma \ref{lem:inner-product-1} yields the following proposition.
\begin{proposition}[{\cite[Lemma 2.2]{BS}}]\label{pro:prop-pdr}
	If ${u={(x, y,\bar x,\bar y)}}$ with $(x, y) \in \dom \mF \times \dom \mG$, $\bar x=x+\sigma\mathcal{K}^{\ast} y$, and $\bar y= y-\tau\mathcal{K}x $, then we have
	\begin{align}
	\mathcal{L}(x_{\test}^{k+1}, y)-\mathcal{L}(x,y_{\test}^{k+1})&\leq\langle u^k-u^{k+1},u^{k+1}-u\rangle_\mathbf{M}\notag\\
	&= \frac{1}{2  }  ( \|u^k-u\|_{\mathbf{M}}^2 - \|u^{k+1}-u\|_{\mathbf{M}}^2 -\|u^{k+1}-u^k\|_{\mathbf{M}}^2 ). \label{eq:polari}
	\end{align}
\end{proposition}
The equality in \eqref{eq:polari} comes from the polarization identity with positive semi-definite weight.

Next, let us introduce the relaxed preconditioned Douglas-Rachford splitting method,
\begin{equation}\label{eq:pdr:relax:deter}
    u^{k+1}=(I-I_{\varrho_k})u^k+I_{\varrho_k}\mathcal Tu^k,
\end{equation}
where $I_{\varrho_k}:=\textrm {Diag}[\rho_k I,\rho_k I,\rho_{k,x}I,\rho_{k,y}I]$. Throughout the paper, we assume the dimension of each identity operator $I$ depends on the corresponding variable. Suppose \(\rho_k\), \(\rho_{k,x}\), and \(\rho_{k,y}\) are non-decreasing numbers with lower and upper limits \(\rho_l\) and \(\rho_u\), respectively, such that \(0 < \rho_l < \rho_u < 2\). With $\mathcal{ T}$ being the iteration operator of \ref{PDR} given by \eqref{eq:pdr:T}, we then define \( u^{k+1}_t = \mathcal{T}u^k := (x^{k+1}_t, y^{k+1}_t, \bar{x}^{k+1}_t, \bar{y}^{k+1}_t) \). In the computation, we also introduce the following new transitional variables:
\begin{equation}
\left\{\begin{array}{ll}
    x_{\test}^{k+1}=(I+\sigma\partial\mathcal{F})^{-1}[2x_t^{k+1}-\bar x^k], \\
     y_{\test}^{k+1}=(I+\tau\partial\mathcal{G})^{-1}[2 y_t^{k+1}- \bar y^k].
    \end{array}\right.
\end{equation}
Writing \eqref{eq:pdr:relax:deter} component-wisely, we get the relaxed and preconditioned Douglas-Rachford iteration as follows
\begin{equation}\label{rPDR}\tag{RPDR}
	\left\{\begin{array}{ll}
		x^{k+1}=(1-\rho_k)x^k+\rho_k x_t^{k+1},   \\
		 y^{k+1}=(1-\rho_k) y^k+\rho_k y_t^{k+1},   \\
		\bar x^{k+1}=\bar x^k+ \rho_{k,x}[ x_{\test}^{k+1}-x_t^{k+1}],\\
		 \bar y^{k+1}= \bar y^k+\rho_{k,y}[y_{\test}^{k+1}- y_t^{k+1}].\\
	\end{array}
	\right.
\end{equation}
Now we define the iteration operator \ref{rPDR} as $\mathcal{T}_R$ such that $u^{k+1} := \mathcal{T}_R u^{k}$. Note that \ref{PDR} is a special case of \ref{rPDR} when $I_{\varrho_k} = I$. Let us define 
$$\mathbf{M}_{\varrho_k} := \mathbf{M} I^{-1}_{\varrho_k} = \diag\left[\frac{1}{\sigma \rho_k}(N_1 - I), \frac{1}{\tau \rho_k}(N_2 - I), \frac{1}{\sigma \rho_{k,x}} I, \frac{1}{\tau \rho_{k,y}} I \right].$$
We will present the following proposition, with the proof provided in the Appendix.
\begin{proposition}\label{pro:prop-rpdr}
	If ${u={(x, y,\bar x, \bar  y)}}$ with $(x, y) \in \dom \mF \times \dom \mG$ and $\bar x=x+\sigma\mathcal{K}^{\ast}  y$, $ \bar y= y-\tau\mathcal{K}x$, then we have
\begin{equation}
	\mathcal{L}(x_{\test}^{k+1}, y)-\mathcal{L}(x,y_{\test}^{k+1})\leq
	\frac{1}{2  }  ( \|u^k-u\|_{\mathbf{M}_{\varrho_k}}^2 - \|u^{k+1}-u\|_{\mathbf{M}_{\varrho_k}}^2 -\|u^{k+1}-u^k\|_{  \mathbf{M}  (2I-I_{\varrho_k})I^{-2}_{\varrho_k}}^2 ). \label{eq:polari-rpdr}
	\end{equation}
\end{proposition}
 An interesting question arises: given a sequence $\{u^{k_i}\}_{k_i}$ generated by \ref{rPDR} with w-$\lim_{i \rightarrow \infty}$ $(\bar{x}^{k_i}, \bar{y}^{k_i}) = (\bar{x}^{*}, \bar{y}^{*})$, where w-$\lim$ denotes the weak limit, under what conditions can we ensure that $\{\mathcal{T}_Ru^{k_i}\}$ converges weakly to a fixed point of $\mathcal{T}_R$. The following lemma provides a condition useful for the almost sure convergence analysis. The proof is also presented in the Appendix.
\begin{lemma}\label{lemma:fix-rpdr}
	The iteration \ref{rPDR} for problem \ref{eq:saddle-point} possesses the following properties:
	\begin{itemize}
		\item [\emph{(i)}] A point $u^{*}= (x^{*}, y^{*}, \bar{x}^{*}, \bar{y}^{*} )$ is a fixed point of \ref{rPDR} if and only if $ (x^{*}, y^{*} )$ is a saddle point for \eqref{eq:saddle-point} and $x^{*}+\sigma \mathcal{K}^{*} y^{*}=\bar{x}^{*}, y^{*}-\tau \mathcal{K} x^{*}=\bar{y}^{*}$.
		\item[\emph{(ii)}] If w-$\displaystyle{\lim_{j  \rightarrow \infty}} (\bar{x}^{k_{j}}, \bar{y}^{k_{j}} )= (\bar{x}^{*}, \bar{y}^{*} )$ and $\displaystyle{\lim _{j  \rightarrow \infty}}( (N_{1}-I)(x^{k_{j}}-x^{k_{j}+1}), (N_{2}-I)(y^{k_{j}}-y^{k_{j}+1}),\bar{x}^{k_{j}}-\bar{x}^{k_{j}+1}, \bar{y}^{k_{j}}-\bar{y}^{k_{j}+1})=(0,0,0,0)$, then $ (x^{k_{j}+1}, y^{k_{j}+1}, \bar{x}^{k_{j}+1}, \bar{y}^{k_{j}+1} )$ converges weakly to a fixed point.
	\end{itemize}
\end{lemma}

\begin{remark}
    If we choose $I_{\varrho_k}:=\textrm {Diag}[ I, I,\rho_{k,x}I,\rho_{k,y}I]$ in relaxed preconditioned Douglas-Rachford method \eqref{eq:pdr:relax:deter}, we get partial relaxed preconditioned Douglas-Rachford method \eqref{p-RPDR}. Further, if we choose $N_1 = I$, $N_2 = I$, and $\rho_{k,x} = \rho_{k,y} $, we derive the original relaxed Douglas-Rachford splitting that was studied in \cite{EP}, of which the weak convergence can be seen in \cite[Theorem 7]{EP} with an ergodic convergence rate. 
\begin{equation}\tag{p-RPDR}\label{p-RPDR}
\left\{
\begin{aligned}
x^{k+1}&=N_1^{-1}[(N_1-I)x^k+\bar x^k-\sigma\mathcal{K}^\ast   y^{k+1}],   \\
y^{k+1}&=N_2^{-1}[(N_2-I)  y^k+ \bar y^k+\tau\mathcal{K} x^{k+1}],   \\
\bar{x}^{k+1} &= \bar{x}^k + \rho_{k,x}\left[(I + \sigma \partial F)^{-1} (2x^{k+1} - \bar{x}^k) - x^{k+1}\right], \\
\bar{y}^{k+1} &= \bar{y}^k + \rho_{k,y}\left[(I + \tau \partial G)^{-1} (2y^{k+1} - \bar{y}^k) - y^{k+1}\right].
\end{aligned}
\right.
\end{equation}
To clarify, we refer to the general \ref{rPDR} as the fully relaxed preconditioned Douglas-Rachford method (f-RPDR) when $I_{\varrho_k} := \textrm{Diag}[\rho_k I, \rho_k I, \rho_{k,x} I, \rho_{k,y} I]$, where $\rho_k \neq 1$, $\rho_{k,x} \neq 1$, and $\rho_{k,y} \neq 1$.
 
\end{remark}

In the special case where $\mF$ is a quadratic-linear form as follows:
\begin{equation}\label{eq:pdrq:F}    
\mathcal F(x)=\langle \frac{1}{2}Qx-f,x\rangle,
\end{equation}
where $Q$ is a self-adjoint, linear, and bounded positive semi-definite operator and $f \in X$,  we can get a more compact preconditioned Douglas--Rachford splitting method for quadratic problems, which is called PDRQ \cite{BSCC, BS, BS1}. Denoting by $ u =(x,  y, \bar y)\in\mathbb{U}_Q:=\mathbb{X}\times \mathbb{Y} \times \mathbb{Y} $, $N_1\succeq 0 $, $N_2\succeq I $, the PDRQ iteration reads as,
\begin{equation}\label{pdrq}\tag{PDRQ}
\left\{\begin{array}{l}
x^{k+1}= (N_{1}+\sigma Q )^{-1} [N_{1} x^{k}-\sigma (\mathcal{K}^{*} y^{k+1}-f ) ], \\
y^{k+1}=N_{2}^{-1} [ (N_{2}-I ) y^{k}+\bar{y}^{k}+\tau \mathcal{K} x^{k+1} ], \\
y_{\test}^{k+1}=(I+\tau \partial G)^{-1} [2 y^{k+1}-\bar{y}^{k} ],\\
\bar{y}^{k+1}=\bar{y}^{k}+y_{\test}^{k+1}-y^{k+1}.
\end{array} \right.
\end{equation}

Without the updates for $\bar x^{k+1}$, the \ref{pdrq} iterations are more compact and efficient compared to \ref{PDR} \cite{BSCC}. Moreover, if we choose $N_2 =I$ and $N_1$ such that
\begin{equation}\label{N1:m1}
N_{1} = M_Q - \sigma Q- \sigma \tau \mK^{*}\mK,
\end{equation}
where $M_Q: X\rightarrow X$ is linear, continuous, self-adjoint and \emph{feasible}  meaning $M_Q \succeq T_Q:=\sigma Q + \sigma \tau  \mK^*\mK$, then we can update $x^{k+1}$ first with the following classical preconditioned iteration 
\begin{equation}\label{preconditioner:from:Q}
x^{k+1} = x^{k} + M_Q^{-1}(b^{k} - Tx^{k}), \quad b^{k} = \sigma f - \sigma \mK^*\bar y^k.
\end{equation}
The linear update of $y^{k+1}=\bar y^k +\tau \mK x^{k+1}$ is then according to the second line in \ref{pdrq} \cite{BSCC, BS, BS1}.
Denoting 
\begin{displaymath}
\mathbf M_Q=\diag[\tfrac1\sigma N_1, \tfrac1\tau (N_2-I), \tfrac{1}{\tau} I],
\end{displaymath}
we have the associated inner product again,
\begin{equation}
\langle u, u^{\prime} \rangle_{\mathbf{M}_{Q}}=\frac{1}{\sigma} \langle N_{1} x, x^{\prime} \rangle+\frac{1}{\tau} [ \langle (N_{2}-I ) y, y^{\prime} \rangle+ \langle\bar{y}, \bar{y}^{\prime} \rangle ].
\end{equation}
Let $\mT_Q$ denote the iteration operator associated with \ref{pdrq}, we can rewrite the iteration \ref{pdrq} as 
\begin{equation}\label{eq:pdrq:T}
u^{k+1}: = \mT_Q u^k.
\end{equation}
Now, we introduce the relaxed PDRQ method,
\begin{equation}\label{eq:pdrq:relax:deter}
    u^{k+1}=(I-I_{q,\varrho_k})u^k+I_{q,\varrho_k}\mT_Q u^k,
\end{equation}
where $I_{q,\varrho_k}:=\textrm {Diag}[\rho_k I,\rho_k I,\rho_{k,y} I]$. Suppose $\rho_k$,  $\rho_{k,y}$ are non-decreasing and have lower and upper limits $ \rho_l, \rho_u$, and $0< \rho_l < \rho_u< 2 $. We 
 also set $ u^{k+1}_t=\mathcal{T}_{{Q}}u^k:=(x^{k+1}_t, y^{k+1}_t,\bar y^{k+1}_t) $ and 
\begin{equation*}
 y_{\test}^{k+1}=(I+\tau\partial\mathcal{G})^{-1}[2 y_t^{k+1}- y^k],
\end{equation*}
Writing \eqref{eq:pdrq:relax:deter} component-wisely, we get the relaxed and preconditioned Douglas-Rachford iteration for the quadratic problem as follows 
\begin{equation}\label{rPDRQ}\tag{RPDRQ}
	\left\{\begin{array}{ll}
		x^{k+1}=(1-\rho_k)x^k+\rho_k x_t^{k+1},   \\
		 y^{k+1}=(1-\rho_k) y^k+\rho_k y_t^{k+1},   \\
		 \bar y^{k+1}= \bar y^k+\rho_{k,y}[y_{\test}^{k+1}- y_t^{k+1}].\\
	\end{array}
	\right.
\end{equation}
We now define the iteration operator \ref{rPDRQ} as $\mathcal{T}_{RQ}$ such that $u^{k+1} := \mathcal{T}_{RQ} u^{k}$.  Let $\mathbf{M}_{Q\varrho_k}:=\mathbf{M}_{Q}I^{-1}_{q,\varrho_k}=\diag[\frac{1}{\sigma\rho_k} N_1, \frac{1}{\tau\rho_k} (N_2-I),  \tfrac{1}{\tau\rho_{k,y}} I]$.
Similarly, we have the following proposition and lemma for the iteration sequence $\{u^k\}_k$ produced by \ref{rPDRQ}. 
\begin{proposition}\label{pro:prop-pdrq}
	If ${u ={(x,  y, \bar y)}}$, with $ \bar y=  y-\tau\mathcal{K}x $, then we have  the following estimate:
	\begin{align}
	\mathcal{L}(x_{\test}^{k+1}, y)-\mathcal{L}(x,y_{\test}^{k+1}) \leq\frac{1}{2  }  (  \|u^k-u\|_{\mathbf{M}_{Q\varrho_k}}^2 - \|u^{k+1}-u\|_{\mathbf{M}_{Q\varrho_k}}^2 -\|u^{k+1}-u^k\|_{\mathbf{M}_Q (2I-I_{q,\varrho_k})I^{-2}_{q,\varrho_k}}^2 ).
	\end{align}
\end{proposition}

\begin{lemma}\label{lemma:fixpdrq}
	The iteration \ref{rPDRQ} for \ref{eq:saddle-point} with $\mF$ as in \eqref{eq:pdrq:F} possesses the following properties:
	\begin{itemize}
		\item [\emph{(i)}] A point $u^{*}= (x^{*}, y^{*}, \bar{y}^{*} )$ is a fixed point of \eqref{eq:pdrq:T} if and only if $ (x^{*}, y^{*} )$ is a saddle point for \eqref{eq:saddle-point} and $Qx^*+\mathcal{K}^*y^*=f $, $ y^{*}-\tau \mathcal{K} x^{*}=\bar{y}^{*}$.
		\item[\emph{(ii)}] If w-$\displaystyle{\lim_{j  \rightarrow \infty}} \bar{y}^{k_{j}}=\bar{y}^{*}$ and $\displaystyle{\lim _{j  \rightarrow \infty}}(  N_{1}(x^{k_{j}}-x^{k_{j}+1}), (N_{2}-I)(y^{k_{j}}-y^{k_{j}+1}), \bar{y}^{k_{j}}-\bar{y}^{k_{j}+1})=(0,0,0)$, then $ (x^{k_{j}+1}, y^{k_{j}+1},  \bar{y}^{k_{j}+1} )$ converges weakly to a fixed point.
	\end{itemize}
\end{lemma}

Now, we introduce the restricted primal-dual gap function and the Bregman distance for the gap estimates. Suppose that the bounded set $\mathbb B:= \mathbb B_1\times\mathbb B_2\subset \dom\mathcal F \times \dom\mathcal G $  contains a saddle point. 
We give the following restricted primal-dual gap $ \mathfrak{G}_{\mathbb B_1\times\mathbb B_2}$ associated with  $\mathbb B$:
\begin{equation}\label{eq:restricted-primal-dual-gap}
\mathfrak{G}_{\mathbb B_1\times\mathbb B_2}(x,y):=\sup_{\tilde y\in \mathbb B_2}\mathcal{L}(x,\tilde y)-\inf_{\tilde x\in  \mathbb B_1}\mathcal{L}(\tilde x,y):=\sup_{\tilde z \in\mathbb{B}} \mH (z, \tilde z), \quad \mH (z, \tilde z): = \mathcal{L}(x,\tilde y)-\mathcal{L}(\tilde x,y),
\end{equation}
where $z=(x,y)$ as before and $\tilde z:=(\tilde x, \tilde y)$. For any proper, convex, and lower semi-continuous function $ \zeta: \mathbb{X}  \rightarrow \mathbb{R}_{\infty}$, the Bregman distance  at $s$ with $ q \in \partial \zeta(s)$ is defined as $$
   D_{\zeta}^{q}(t,s)=\zeta(t)-\zeta(s)-\langle q,t-s\rangle.  
$$ If $\tilde z$ is a saddle point, we can interpret $\mH (z, \tilde z)$ as the following Bregman distance for $h(z):=\mF(x)+\mG(y)$ with respect to the point $q:= (-\mK^{*} \tilde y, \mK \tilde x )$ which is included in the subgradient $\partial h (\tilde  z  )$ by optimality, i.e.,
\begin{equation*}
\mH (z, \tilde z):=D_{h}^{q} (z, \tilde z )=D_{\mF}^{-\mK^{*} \tilde y} (x, \tilde x  )+D_{\mG}^{\mK \tilde x} (y,\tilde y ).
\end{equation*}
The following theorem tells us that over-relaxation can bring certain acceleration for the deterministic PDR method and motivates us to focus on the stochastic and over-relaxed PDR. The proof is given in the Appendix.
\begin{theorem}
\label{thm:convergence-rpdr}
	If $u=(x,y,\bar x,\bar y),$ with $ \bar x=x+\sigma\mathcal K^{\ast}y,$ $\bar y=y-\tau\mathcal Kx$, and $(x,y)$ being a saddle-point of \eqref{eq:saddle-point}. Let $u^k$ be the iteration sequence generated by \ref{rPDR}. Denoting $K_1 = K+1$, let $z^k_{\test}=(x^k_{\test},y^k_{\test}) $ be the transitional variables of \ref{rPDR}, then the ergodic sequence $ z_{\test,K}=\frac{1}{K+1}\sum_{k=0}^{K} z^k_{\test} $ converges with rate $\mathcal{O}({1}/{K})$ in a
	restricted primal-dual gap sense, i.e., for any $ z=(x,y) \in  \mathbb B_1 \times\mathbb B_2 =\mathbb B \subset \dom\mathcal{F}\times\dom\mathcal{G}$, a bounded domain and containing a saddle-point, it holds that
\begin{equation}\label{eq:rpdr-gap-convergence}
	\mathfrak{G}_{z\in\mathbb B}(z_{\test,K})=\sup_{z\in\mathbb B}\mH(z_{\test,K}, z) \leq 
 \frac{1}{2K_1 \rho_{l} } \sup_{(x,y) \in \mathbb{B}_1 \times \mathbb{B}_1} \|u^0-u\|_{\mathbf{M}}^2 = \mathcal{O}(\frac{1}{K}).
	\end{equation}
\end{theorem}

With these preparations for the deterministic relaxed PDR method, i.e., \ref{rPDR} (or \ref{rPDRQ}), let us now shift our focus to the stochastic and relaxed PDR method (or stochastic and relaxed PDRQ method).

\subsection{The stochastic relaxed preconditioned DR (SRPDR) method}
We begin by introducing the stochastic setting. Let $(\Omega, \mathfrak{F}, \mathcal{P})$ denote the underlying probability space.
Denote $\{i_0, i_1, i_2, \ldots\}$ as index selection variables, i.e., independent random variables on the latter space with values in $\{1,\ldots,n\}$ which are identically distributed. We denote the
probability of selecting an index $i\in\{1,\ldots,n\} $ by $p_i$ where we assume $0<p_{i}<1$ for each $i$ and recall that $\sum_{i=1}^{n} p_{i}=1$. Introduce $P=\diag[p_1,p_2,\ldots,p_n]$ and observe that $P^{-1}=\diag[\frac{1}{p_1},\frac{1}{p_2},\ldots,\frac{1}{p_n}]$.  
A Borel random vector $y$ in a separable Hilbert space $\mathbb{H}$ is a measurable map from $\Omega$ into $\mathbb{H}$ where $\Omega$ is equipped with $\mathfrak{F}$ and $\mathbb{H}$ is equipped with the Borel $\sigma$-algebra, i.e., the $\sigma$-algebra generated by the open sets of $\mathbb{H}$ (see \cite{CPE} or  \cite[Chapter 2]{LT}). Further, 
$(\mathfrak{F}_{k})_{k \in \mathbb{N}}$ denotes the filtration generated by the indices $ \{i_{0}, i_{1}, \ldots, i_{k} \}$, 
i.e., the $\sigma$-algebras generated by $i_j^{-1}\{i\}$ for $j=0,\ldots,k$ and $i = 1,\ldots,n$ \cite[Chapter 7]{SAN}. Then, $\mathfrak{F}_k \subset \mathfrak{F}_{k+1}$ for each $k \in \mathbb{N}$. Henceforth, we denote $\mathbb{E}^{k}:=\mathbb{E} [\cdot \mid \mathfrak{F}_{k} ]$ as the conditional expectation with respect to $\mathfrak{F}_{k}$, and by $\mathbb{E}$ the full expectation. Recall that if a Borel random vector $y$ is $\mathfrak{F}_{k}$-measurable, i.e., the preimages of Borel sets are in $\mathfrak{F}_k$, then $\mathbb{E}[y \mid \mathfrak{F}_{k} ] = y$. 

Let us first extend the relaxed preconditioned Douglas--Rachford splitting method, i.e., \ref{rPDR} to a stochastic setting. For arbitrary positive step sizes $\sigma $, $\tau$ and $(x^0,y^0,\bar{x}^0,\bar{y}^0)\in(\mathbb{B}_1 \times \mathbb{B}_2)^2$, $y^0_{\test}\in\mathbb{B}_2 $, the stochastic relaxed PDR (\ref{sto-rpdr}) reads as follows
\begin{equation}\label{sto-rpdr} \tag{SRPDR}
\left\{\begin{array}{ll}
\tilde{x}_t^{k+1}=  \left(N_{1}+\sigma \tau \mathcal{K}^{*} N_{2}^{-1} \mathcal{K}\right)^{-1}\left(\left(N_{1}-I\right) x^{k}+\bar{x}^{k}-\sigma \mathcal{K}^{*} N_{2}^{-1}\left(\left(N_{2}-I\right) y^{k}+\bar{y}^{k}\right)\right), \\
x^{k+1} = (1-\rho_k)x^k+\rho_k \tilde{x}_t^{k+1},  \\
y^{k+1}= \left\{
\begin{array}{ll}
(1-\rho_k)y_i^k+\rho_k \tilde{y}^{k+1}_{t,i}, \quad  \tilde{y}^{k+1}_{t,i} := N_{2,i}^{-1}[(N_{2,i}-I) y_i^k+\bar y_i^k+\tau\mathcal{K}_i \tilde{x}_t^{k+1}]& i= i_k,  \\
y^k_i &\forall i \neq i_k,\\
\end{array}  \right.\\
x_{\test}^{k+1}=(I+\sigma\partial\mathcal{F})^{-1}[2\tilde{x}_t^{k+1}-\bar x^k],\\
y_{\test}^{k+1}= \left\{
\begin{array}{ll}
(I+\tau\partial\mathcal G_i)^{-1}[2\tilde{y}^{k+1}_{t,i}-\bar y^k_i]  & i= i_k,  \\
y_{\test,i}^{k} & \forall i \neq i_k,\\
\end{array}  \right.\\
\bar x^{k+1}=\bar x^k+ \rho_{k,x} (x_{\test}^{k+1}-\tilde{x}_t^{k+1} ),\\
\bar y^{k+1}= \left\{
\begin{array}{ll}
\bar y_i^k+\rho_{k,y}(y_{\test,i}^{k+1}-\tilde{y}_{t,i}^{k+1}) & i= i_k,  \\
\bar y^k_i & \forall i \neq i_k, 
\end{array}  \right.
\end{array}
\right.
\end{equation}
where $ \tilde{u}^{k+1}_t=\mathcal{T}u^k:=(\tilde{x}^{k+1}_t, \tilde{y}^{k+1}_t,\tbx^{k+1}_t,\tby^{k+1}_t) $ and $\mathcal{T}$ is the iteration operator of \ref{PDR}.
 $N_1 \succeq I $ and $N_2=\diag[N_{2,1}, \ldots,  N_{2,n}]$ with each $N_{2,i}: \mathbb{Y}_i \rightarrow \mathbb{Y}_i$ being self-adjoint, bounded and linear such that $N_{2,i} \succeq I$. Again, let \( I_{\varrho_k} := \text{Diag}[\rho_k I, \rho_k I, \rho_{k,x} I, \rho_{k,y} I] \), where \(\rho_k\), \(\rho_{k,x}\), and \(\rho_{k,y}\) are non-decreasing numbers with lower and upper limits \(\rho_l\) and \(\rho_u\), respectively, such that \(0 < \rho_l < \rho_u < 2\). Although
the updates of $x^{k+1}$ in the above \ref{sto-rpdr} method seem complicated, it is convenient for analysis. Specifically, the updates of $\tilde{x}_t^{k+1} $ is based on preconditioned iteration given in \eqref{preconditioner:from} when $N_2=I$, and $N_1$ as in \eqref{N1:m}.

To clarify, we refer to the general \ref{sto-rpdr} as the stochastic fully relaxed preconditioned Douglas-Rachford method (f-SRPDR) when $I_{\varrho_k} = \textrm{Diag}[\rho_k I, \rho_k I, \rho_{k,x} I, \rho_{k,y} I]$, where $\rho_k \neq 1$, $\rho_{k,x} \neq 1$, and $\rho_{k,y} \neq 1$. While when $I_{\varrho_k} = \textrm{Diag}[I, I, \rho_{k,x} I, \rho_{k,y} I]$, where  $\rho_{k,x} \neq 1$, and $\rho_{k,y} \neq 1$, we call it stochastic partial relaxed preconditioned Douglas-Rachford method (p-SRPDR). In addition, by letting $I_{\varrho_k} = I$, the stochastic PDR, i.e., \ref{sto-pdr} reads as follows
\begin{equation}\label{sto-pdr} \tag{SPDR}
\left\{\begin{array}{ll}
x^{k+1}=  \left(N_{1}+\sigma \tau \mathcal{K}^{*} N_{2}^{-1} \mathcal{K}\right)^{-1}\left(\left(N_{1}-I\right) x^{k}+\bar{x}^{k}-\sigma \mathcal{K}^{*} N_{2}^{-1}\left(\left(N_{2}-I\right) y^{k}+\bar{y}^{k}\right)\right),   \\
y^{k+1}= \left\{
\begin{array}{ll}
N_{2,i}^{-1}[(N_{2,i}-I) y_i^k+\bar y_i^k+\tau\mathcal{K}_i x^{k+1}]& i= i_k,\\
y^k_i &\forall i \neq i_k,\\
\end{array}  \right.\\
x_{\test}^{k+1}=(I+\sigma\partial\mathcal{F})^{-1}[2x^{k+1}-\bar x^k],\\
y_{\test}^{k+1}= \left\{
\begin{array}{ll}
(I+\tau\partial\mathcal G_i)^{-1}[2y^{k+1}_{i}-\bar y^k_i]  & i= i_k,  \\
y_{\test,i}^{k} & \forall i \neq i_k,\\
\end{array}  \right.\\
\bar x^{k+1}=\bar x^k+(I+\sigma\partial\mathcal{F})^{-1}[2x^{k+1}-\bar x^k]-x^{k+1},\\
\bar y^{k+1}= \left\{
\begin{array}{ll}
\bar y_i^k+y_{\test,i}^{k+1}-y_i^{k+1} & i= i_k,  \\
\bar y^k_i & \forall i \neq i_k. 
\end{array}  \right.
\end{array}
\right.
\end{equation}
Considering a simpler case when $N_1 = I$ and $N_2 = I$, we will derive the stochastic DR method introduced in \cite{CPE}. 

With these preparations, now let us turn to our algorithms' almost sure convergence and restricted gap estimates. Since we can treat \ref{sto-pdr} as a special case of \ref{sto-rpdr}, our analysis will focus solely on the convergence properties and gap estimates for \ref{sto-rpdr}.

\section{Almost sure convergence of stochastic \ref{rPDR}}\label{sec:almost:sure}
This section presents the almost sure convergence of iteration sequences generated by \ref{sto-rpdr} to a solution of \eqref{eq:saddle-point}. To avoid any ambiguity, henceforth,  the notation  $u^k=(x^k, y^k,\bar x^k,\bar y^k)$ along with $ (x^k_{\test},y^k_{\test})$ always represents the stochastic sequences produced by \ref{sto-rpdr}. To write the iteration in a more compact form, we introduce the notation $\tilde u^{k+1}$ for one step \ref{rPDR} update based on $u^k$, i.e., 
\begin{equation} \label{eq:tildeu:after:u}
\tilde u^{k+1}:= \mathcal{T}_Ru^k, \quad \text{with} \  \ \tilde u^{k+1}= (\tilde x^{k+1},\tilde y^{k+1},\tbx^{k+1}, \tby^{k+1}), 
\end{equation}
where $\mathcal{T}_R$ is the iteration operator of \ref{rPDR}, and let $(\tilde x^{k+1}_{\test},\tilde y^{k+1}_{\test})$ be the transitional variables of this one step \ref{rPDR} update.  Further, we introduce the notation $\tilde u_t^{k+1}$ for one step deterministic \ref{PDR} update based on $u^k$, i.e., 
\begin{equation} \label{eq:tildeu:after:u}
\tilde u_t^{k+1}:= \mathcal{T}u^k, \quad \text{with} \  \ \tilde u_t^{k+1}= (\tilde x_t^{k+1},\tilde y_t^{k+1},\tbx_t^{k+1}, \tby_t^{k+1}),
\end{equation}
where $\mathcal{T}$ is the same as in  \eqref{eq:pdr:T},  and let $(\tilde x^{k+1}_{t,\test},\tilde y^{k+1}_{t,\test})$ be the transitional variables of this one step \ref{PDR} update.
Hence the stochastic update $u^{k+1}$ in \ref{sto-rpdr} can then be written as
\begin{equation}\label{eq:same}
x^{k+1} = \tilde x^{k+1},\quad  x_{\test}^{k+1} = \tilde x_{\test}^{k+1} = \tilde x_{t,\test}^{k+1},\quad \bar x^{k+1} = \tbx^{k+1},
\end{equation}
and
\begin{equation} \label{eq:sto-different}
    y^{k+1}_i = \begin{cases}
        \tilde y^{k+1}_i & \text{if} \ i = i_k, \\
        y^k_i & \text{if} \ i \neq i_k,
    \end{cases} \quad
    y^{k+1}_{\test,i} = \begin{cases}
        \tilde y^{k+1}_{\test,i} ={\tilde y^{k+1}_{t,\test,i}} & \text{if} \ i = i_k, \\
        y^k_{\test,i} & \text{if} \ i \neq i_k,
    \end{cases}
    \quad
    \bar y^{k+1}_{i} = \begin{cases}
        \tby^{k+1}_{i} & \text{if} \ i = i_k, \\
        \bar y^k_{i} & \text{if} \ i \neq i_k.
        \end{cases}
\end{equation}
Note that $x^k, x^k_{\test}, \bar x^k, y^k, y^k_{\test}, \bar y^k$ according to \eqref{eq:same} and \eqref{eq:sto-different} are $\mathfrak{F}_k$-measurable Borel random vectors by continuity of $\mathcal{T}_R$ and construction.

The stochastic iteration sequence $\{ u^k\}_k$ and the corresponding one-step deterministic update $\{\tilde u^{k+1}\}_{k+1}$ can be shown in the following diagram
\begin{equation}
\begin{aligned}
& u^0 \rightarrow u^1 \rightarrow u^2 \rightarrow u^3 \rightarrow\cdots  \rightarrow u^{k} \xrightarrow{\ref{sto-rpdr}} u^{k+1} \cdots\\
&\downarrow \ \; \; \quad  \downarrow \; \;\;\;  \quad  \downarrow\; \;\;\;  \quad \downarrow \; \;\; \qquad \qquad  \downarrow{\mathcal{T}_R}\\
& \tilde u^1 \, \ \quad \tilde u^2 \, \ \quad \tilde u^3 \, \ \quad \tilde u^4 \, \ \quad\cdots  \, \ \quad \tilde u^{k+1}= \mathcal{T}_Ru^k \cdots.\\
\end{aligned}
\end{equation}
Let us begin the convergence analysis with the following estimate.
\begin{lemma}\label{lem:almost}For any $u=(x, y,\bar x, \bar y)$ with $(x, y) \in \dom \mF \times \dom \mG$ and $\bar x=x+\sigma\mathcal{K}^{\ast} y$ and $\bar y= y-\tau\mathcal{K}x $, it holds that
	\begin{align}\label{eq:H-tildey}
	\mH \bigl((x_{\test}^{k+1},\tilde y^{k+1}_{\test}), (x,y) \bigr) &\leq\frac{1}{2}\left( \|  u^k- u\|^2_{\mathbf{M}_{p,\varrho_k}}-\mathbb{E}^k\| u^{k+1}-u\|^2_{\mathbf{M}_{p,\varrho_k}}-\mathbb{E}^k\|u^{k+1}-u^k\|^2_{(2I-I_{\varrho_k})I^{-2}_{\varrho_k}\mathbf{M}_p}\right),
	\end{align}
 where $\mathbf{M}_{p}:=\diag[\frac{1}{\sigma} (N_1-I), \frac{1}{\tau}(N_2-I)\mathbf{P}^{-1},\frac{1}{\sigma} I, \frac{1}{\tau}\mathbf{P}^{-1}]$, and $\mathbf{M}_{p,\varrho_k}:=I^{-1}_{\varrho_k}\mathbf{M}_p
 $. Here $\mathbf{P} = \diag[p_1 I,\ldots,p_n I]$,
 consequently, $\mathbf{P}^{-1}= \diag[p_1^{-1} I,\ldots,p_n^{-1} I]$.
\end{lemma}	
\begin{proof}
		First, noting that for $\mathfrak{F}_k$-measurable Borel random vectors $ s^k =(s^k_{i})_{i=1,\ldots, n} $, $\tilde s^{k+1} =(\tilde s^{k+1}_{i})_{i=1,\ldots, n}$ on $\mathbb{Y}$, the update rule
        $$ 
		s_{i}^{k+1}= \left\{\begin{array}{ll}
		\tilde{s}_{i}^ {k+1} & i =i_k, \\
		s_{i}^{k} & \forall i \neq i_k,
		\end{array} \right.
		$$
        and any continuous $\varphi: Y \to \mathbb{R}$, it holds that
        \begin{equation}\label{eq:Ek}
		\mathbb{E}^k\varphi(s_i^{k+1})=p_i\varphi(\tilde s_i^{k+1})+(1-p_i)\varphi(s_i^{k}),
		\end{equation}
		which can be rewritten to 
		\begin{equation}
		\varphi (\tilde{s}_{i}^{k+1} )=\frac{1}{p_{i}}\mathbb{E}^{k}\varphi (s_{i}^{k+1} )- (\frac{1}{p_{i}}-1 ) \varphi (s_i^{k} ). \label{eq:trival2}
		\end{equation}
		In particular, for any block diagonal operator $\mathbf{B}=\diag [\mathbf{B}_{1}, \ldots, \mathbf{B}_{n} ]$ of symmetric and positive semidefinite operators $\mathbf{B}_1, \ldots, \mathbf{B}_n$, this implies that 
		\begin{equation}\label{eq:trival}
		\|\tilde{s}^{k+1}-\cdot\|_{\mathbf{B}}^{2} - \|s^{k}-\cdot\|_{\mathbf{B}}^{2} =\mathbb{E}^{k}\|s^{k+1}-\cdot\|_{\mathbf{P}^{-1}\mathbf{ B}}^{2}-\|s^{k}-\cdot\|_{\mathbf{P}^{-1} \mathbf{B}}^{2}.
		\end{equation}
Then with Proposition \ref {pro:prop-rpdr}, $u = (x,y,\bar x, \bar y)$ and $\tilde u^{k+1}$ as in \eqref{eq:tildeu:after:u},  \eqref{eq:same} and \eqref{eq:sto-different}, we arrive at 
		\begin{align}
		\mH \bigl((x_{\test}^{k+1}, \tilde y^{k+1}_{\test}), (x,y) \bigr) &\leq\frac{1}{2} \big(   \|u^k-u\|_{\mathbf{M}_{\varrho_k}}^2 - \|\tilde u^{k+1}-u\|_{\mathbf{M}_{\varrho_k}}^2 -\|\tilde u^{k+1}-u^k\|_{   (2I-I_{\varrho_k})I^{-2}_{\varrho_k}\mathbf{M}   }^2 \big)
		\end{align}
		Next, reformulating $\| \tilde y^{k+1}-\cdot \|_{(N_2-I)\rho^{-1}_k\tau^{-1}}^2 $ and	$\| \tby^{k+1}-\cdot \|_{\rho^{-1}_{k,y}\tau^{-1}I}^2 $ in the above inequality by using  \eqref{eq:trival}, then together with the fact that $x^{k+1}$ and $\bar x^{k+1}$ are $\mathfrak{F}_{k}$ measurable, we get 
        \[
        \begin{aligned}
             \|\tilde u^{k+1} - u\|_{\mathbf{M}_{\varrho_k}}^2 - \|u^{k} - u\|_{\mathbf{M}_{\varrho_k}}^2 &= \mathbb{E}^k\|u^{k+1} - u\|_{\mathbf{M}_{p,\varrho_k}}^2 - \|u^{k} - u\|_{\mathbf{M}_{p,\varrho_k}}^2, \\
             \|\tilde u^{k+1} - u^k\|_{(2I-I_{\varrho_k})I^{-2}_{\varrho_k}\mathbf{M}}^2 &= \mathbb{E}^k\|u^{k+1} - u^k\|_{(2I-I_{\varrho_k})I^{-2}_{\varrho_k}\mathbf{M}_p}^2,
        \end{aligned}
        \]
        and consequently, \eqref{eq:H-tildey}.  
	\end{proof}
Moreover, we see, for instance,
\[ (N_2-I)\mathbf{P}^{-1}=(p_1^{-1}(N_{2,1}-I),\ldots,p_n^{-1}(N_{2,n}-I))=\mathbf{P}^{-1}(N_2-I). \]
For this reason, we will tacitly use that $\mathbf{P}$ (and $\mathbf{P}^{-1}$) always commutes with the linear operators considered in this paper unless otherwise specified.

Let us give some notations first. Denote $\mathcal{Z}^* $ as the set for the saddle points of \eqref{eq:saddle-point-sum} and $\mathbb{F} \subset \mathbb{U}$ as the set for the fixed points for $\mathcal{T}_{R}$. 
Now, let us turn to the almost sure weak convergence of $\{ u^k \}_k$ to a fixed point of $\mathcal{T}_{R}$ as $k \to \infty$. Inspired by \cite{CPE}, we mainly employ the stochastic quasi-Fejér monotonicity \cite{CPE}, Opial's lemma \cite[Lemma 21]{AC} or \cite{Op1}, and Lemma \ref{lemma:fix-rpdr} to show the almost sure convergence of sequence $\{\tilde u^{k+1}, k\in \mathbb{N}\}$ first. We then ``glue" $\{u^{k+1}, k \in \mathbb{N}\}$ to $\{\tilde u^{k+1}, k\in \mathbb{N}\}$ and finish the proof. 

\begin{theorem}\label{thm:almostsure}
	The iteration sequence $\{u^k\}_k = \{(x^k,y^k,\bar x^k, \bar y^k)\}_k$ produced by \ref{sto-rpdr} weakly converges to an $\mathbb{F}$-valued random variable $u^{**} =(x^{**},y^{**},\bar x^{**}, \bar y^{**})$ almost surely, and $(x^{**}, y^{**})$ is a $\mathcal{Z}^*$-valued random variable with 
$x^{**}+\sigma \mathcal{K}^{*} y^{**}=\bar{x}^{**}$ and ${y}^{**}-\tau \mathcal{K} {x}^{**}=\bar{y}^{**}$. Moreover, the transitional sequence $(x^{k}_{\test}, y^{k}_{\test})$ 
    converges weakly to a $\mathcal{Z}^* $-valued random variable $(x^{**},y^{**})$ almost surely.
\end{theorem}
	\begin{proof}
		Inserting $u^*=(x^*,y^*,\bar{x}^*,\bar{y}^*)$ a fixed point of $\mT_R$ in \eqref{eq:H-tildey}. By Lemma \ref{lemma:fix-rpdr}, $(x^*,y^*)$ is a saddle point of  \eqref{eq:saddle-point-sum}. Then it follows for every $(x^*,y^*,\bar{x}^*,\bar{y}^*)\in \mathbb{F}$, $\mathcal{H}((x_{\test}^{k+1},\tilde y_{\test}^{k+1}),(x^*,y^*)) \geq 0$. Employing Lemma \ref{lem:almost} as well as 
  $\mathbf{M}_{p,\varrho_{k+1}} \preceq \mathbf{M}_{p,\varrho_k} $ we have 
\begin{equation*}
		\begin{aligned}\label{eq:e-delta}
		\mathbb E^{k}\|u^{k+1}-u^* \|^2_{\mathbf{M}_{p,\varrho_{k+1}}}  &\leq \mathbb E^{k}\|u^{k+1}-u^* \|^2_{\mathbf{M}_{p,\varrho_k}}\\
		&\leq\|u^{k}-u^* \|^2_{\mathbf{M}_{p,\varrho_k}}-\mathbb E^{k}\|u^{k+1}-u^k \|^2_{(2I-I_{\varrho_k})I^{-2}_{\varrho_k}\mathbf{M}_p},
		\end{aligned}
\end{equation*}
denoting $\Delta^{k}=\|u^{k}-u^* \|^2_{\mathbf{M}_{p,\varrho_k}}$, and due to $\mathbb{E}^k\|u^{k+1} - u^k\|^2_{(2I-I_{\varrho_k})I^{-2}_{\varrho_k}\mathbf{M}_p} = \|\mathcal{T}_Ru^k - u^k\|^2_{(2I-I_{\varrho_k})I^{-2}_{\varrho_k}\mathbf{M}}$, we get,  
		\begin{align}\label{eq:e-delta}
		\mathbb E(\Delta^{k+1}|\mathfrak{F}_k)=\mathbb E^{k}\Delta^{k+1}  
        &\leq\Delta^{k}-\|\mathcal{T}_Ru^{k}-u^k \|^2_{(2I-I_{\varrho_k})I^{-2}_{\varrho_k}\mathbf{M}}.
		\end{align}
		Using the Robbins--Siegmund theorem \cite{RM2} (see also \cite{CPE,AFC}) on \eqref{eq:e-delta}, we conclude that $\sum_{k=0}^{\infty}\| \mathcal{T}_Ru^{k} - u^k \|^2_{(2I-I_{\varrho_k})I^{-2}_{\varrho_k}\mathbf{M}}< \infty$. Due to $\rho_k, \rho_{k,x}, \rho_{k,y}$ are non-decreasing, thus $\frac{2-\rho_k}{\rho_k^2}$, $\frac{2-\rho_{k,x}}{\rho_{k,x}^2}$ and $\frac{2-\rho_{k,y}}{\rho_{k,y}^2}$ are decreasing and the lower bound is $\frac{2-\rho_{u}}{\rho_{u}^2}$. Then it follows that 
		\begin{equation}\label{eq:M:summable}
		\sum_{k=0}^{\infty}\frac{2-\rho_{u}}{\rho_{u}^2}\| \mathcal{T}_Ru^{k} - u^k \|^2_{\mathbf{M}}< \infty, \quad \| \mathcal{T}_Ru^{k} -u^k \|^2_{\mathbf{M}}  \rightarrow 0,
		\end{equation}
        almost surely in $\Omega$. Recalling the definition of $\mathbf{M}_p$, we immediately see that also $\|\mathcal{T}_Ru^k - u^k\|_{\mathbf{M}_p} \to 0$ as $k \to \infty$ almost surely.
  Still with the Robbins--Siegmund theorem or Proposition 2.3 in \cite{CPE}, we conclude that  $\Delta^{k}=\|u^{k}-u^* \|^2_{\mathbf{M}_{p,\varrho_k}}$ converges to a finite random variable almost surely. 
        Moreover, $\|u^k - u^*\|_{\mathbf{M}_p}$ also converges almost surely by the non-decreasing and convergence of  $\rho_k, \rho_{k,x}, \rho_{k,y}$. Using further that $\|u^k - u^*\|_{\mathbf{M}_p} = \|\mathbf{M}_p^{1/2}(u^k - u^*)\|$ for each $k$, 
        we conclude that $\{\mathbf{M}_p^{1/2}u^k\}_k$ is also bounded almost surely and that the limit
		\begin{align}\label{eq:unormlimit}
		l(u^*): =\lim_{k\rightarrow +\infty}\|u^k-u^* \|_{\mathbf{M}_p} 
		\end{align}
		exists almost surely. 
       The existence and boudnedness of $\mathbf{M}_p^{1/2}$ can be guaranteed by the positive semi-definiteness of $N_1-I$ and $N_2-I$ and the positive definiteness of $\mathbf{P}$.  In the following, we will establish the almost sure weak convergence of $\{u^k\}_k$.
		
		Our first goal is to show that $\{\tilde u^{k}\}_k$ converges weakly to  an $\mathbb{F}$-valued random variable, i.e.,  a fixed point of $\mathcal{T}_R$. 
		We will begin with proving the almost sure weak convergence of $\{\mathbf{M}_p^{1/2}\tilde u^{k} \}_k$ with Opial's lemma \cite[Lemma 21]{AC}. This can be done by showing that almost surely, $\|\mathbf{M}_p^{1/2}(\tilde u^{k+1} - u^*) \|$ converges and every weak sequential limit point of $\{\mathbf{M}_p^{1/2} \tilde u^{k+1}\}_k$ is a $(\mathbf{M}_p^{1/2}\mathbb{F})$-valued random variable. Here $\mathbf{M}_p^{1/2}\mathbb{F}$ is defined as follows
\begin{equation*}
\mathbf{M}_p^{1/2}\mathbb{F} :=  \{((N_1-I)^{\frac{1}{2}} x^*/\sigma^{\frac{1}{2}},(N_2-I)^{\frac{1}{2}} \mathbf{P}^{-\frac{1}{2}}y^*/\tau^{\frac{1}{2}} ,\bar{x}^*/\sigma^{\frac{1}{2}}, \mathbf{P}^{-\frac{1}{2}}\bar{y}^*/\tau^{\frac{1}{2}})\ | \  u^*=(x^*,y^*,\bar{x}^*,\bar{y}^*) \in\mathbb{F}\}, 
\end{equation*}
  with  $\mathbf{P}^{-\frac{1}{2}}= \diag[p_1^{-\frac{1}{2}} I,\ldots,p_n^{-\frac{1}{2}} I]$.
  According to the triangle inequality, \eqref{eq:M:summable}, we have 
         \[ 
         \big| \|\mathbf{M}_p^{1/2}(\tilde u^{k+1}- u^*)\| - \|\mathbf{M}_p^{1/2}(u^k - u^*)\| \big| \leq \|\mathbf{M}_p^{1/2}(\tilde u^{k+1} - u^k)\| = \|\mathcal{T}_Ru^k - u^k\|_{\mathbf{M}_p} \rightarrow 0.
         \] 
         Combining the fact of \eqref{eq:unormlimit}, we derive, 
		\begin{equation}\label{eq:tildeunorm}
		\| \mathbf{M}_p^{1/2}(\tilde u^{k+1} -u^*)\|\rightarrow l(u^*)  
		\end{equation}
        almost surely as $k \to \infty$. In particular, $\{\mathbf{M}_p^{1/2} \tilde u^{k} \}_k$ is  bounded almost surely. 
        Next, we show that for each event such that \eqref{eq:M:summable}, \eqref{eq:unormlimit} and \eqref{eq:tildeunorm} hold, each weak accumulation point of 
        $\{\tilde u^k\}_k$ 
        is a fixed point.
        For that purpose, note that \eqref{eq:tildeunorm}, we have for the sequence $\{\mathbf{M}_p^{1/2} \tilde u^k\}_k$ and each $\mathbf{M}_p^{1/2} u^*$  being  a $(\mathbf{M}_p^{1/2} \mathbb{F})$-valued random variable that $\|\mathbf{M}_p^{1/2} \tilde u^k - \mathbf{M}_p^{1/2} u^*\|$ converges. Now, let $\{\mathbf{M}_p^{1/2} \tilde u^{k_j} \}_j$ be a weakly converging subsequence. Then, as $\|\mathbf{M}_p^{1/2}(\mathcal{T}_Ru^k - u^k)\| \to 0$ as $k \to \infty$, $\{\mathbf{M}_p^{1/2} u^{k_j-1}\}_j$ converges weakly to the same limit 
        and looking at the definition of $\mathbf{M}_p$, we conclude that $\bigl\{\bigl((N_1 - I)x^{k_j-1}, (N_2 - I)y^{k_j-1}, \bar x^{k_j-1}, \bar y^{k_j-1} \bigr)\bigr\}_j$ 
        converges weakly. 
        We further have that $\bigl((N_1 - I)(\tilde x^{k_j} - x^{k_j-1}), (N_2 - I)(\tilde y^{k_j} - y^{k_j-1}), \tbx^{k_j} - \bar x^{k_j-1}, \tby^{k_j} - \bar y^{k_j-1} \bigr) \to 0$ as $j \to \infty$ such that with Lemma \ref{lemma:fix-rpdr} (ii), we obtain that $\{(\tilde x^{k_j}, \tilde y^{k_j}, \tbx^{k_j}, \tby^{k_j})\}$ converges weakly to a fixed point $(x^{**}, y^{**}, \bar x^{**}, \bar y^{**})$ of $\mathcal{T}_R$. In particular, by Lemma \ref{lemma:fix-rpdr} (i), $(x^{**}, y^{**})$ is a saddle point for \eqref{eq:saddle-point-sum} and $x^{**} + \sigma \mathcal{K}^* y^{**} = \bar x^{**}, y^{**} - \tau \mathcal{K} x^{**} = \bar y^{**}$. It is then immediate that the weak limit of $\{\mathbf{M}_p^{1/2} \tilde u^{k_j}\}_j$ is a $(\mathbf{M}_p^{1/2}\mathbb{F})$-valued random variable. By Opial's lemma, we have that the whole sequence $\{\mathbf{M}_p^{1/2} \tilde u^{k} \}_k$ converges weakly to a $(\mathbf{M}_p^{1/2}\mathbb{F})$-valued random variable, implying that $ (\tbx^{k},\tby^{k})$ converges weakly to $(\bar x^{**},\bar y^{**})$. Then by \eqref{eq:M:summable} again, we have $\bigl((N_1 - I)(\tilde x^{k} - x^{k-1}), (N_2 - I)(\tilde y^{k} - y^{k-1}), \tbx^{k} - \bar x^{k-1}, \tby^{k} - \bar y^{k-1} \bigr) \to 0$ as $k \to \infty$. Applying Lemma \ref{lemma:fix-rpdr} (ii) again, we conclude that $\{\tilde u^k\}_k$ converges weakly to a $\mathbb{F}$-valued random variable.

        Finally, let us turn to the almost sure weak convergence of $\{u^k\}_k$. First, note that since for fixed $i=1,\ldots,n$, the probability $p_i$ of $i_k = i$ is positive and all $i_k$ are independent random variables, the probability of $i_k = i$ only occurring finitely many times has to be zero. For that reason, almost surely, for each $i=1,\ldots,n$, the sequences $\{k_{i,j}\}_j$ such that $k_{i,j} = \max \ \{k : k \leq j, \ i_{k-1} = i\}$, i.e., the sequences that yield the index of the last activation of $i$, are increasing and satisfy $\lim_{j \to \infty} k_{i,j} = \infty$. Fix an event for which $\{\tilde u^k\}_k$ converges weakly to a fixed point $u^{**} = (x^{**}, y^{**}, \bar x^{**}, \bar y^{**})$ of $\mathcal{T}_R$ and for which $k_{i,j}$ tends to infinity as $j \to \infty$ for all $i=1,\ldots,n$. We will show that $\{u^k\}_k$ converges weakly, as $k \to \infty$, to $u^*$. For that purpose, observe that by \eqref{eq:same}, we have $x^k = \tilde x^k \to x^{**}$, $x^k_{\test} = \tilde x^{k} + \tbx^{k} - \bar x^{k-1} \to x^{**}$ and $\bar x^k = \tbx^k \to \bar x^{**}$ in the weak sense as $k \to \infty$. By \eqref{eq:sto-different} and the definition of $k_{i,j}$, we see that $y^j_i = \tilde y^{k_{i,j}}_i$, $y^j_{\test,i} = \tilde y^{k_{i,j}}_i + \tby^{k_{i,j}}_i - \bar y^{(k_{i,j}) - 1}_i$ and $\bar y^j_i = \tby^{k_{i,j}}_i$ for each $j$. Hence, letting $j \to \infty$, we see that $y^j_i \to y^{**}_i$, $y^j_{\test,i} \to y^{**}_i$ and $\bar y^j_i \to \bar y^{**}_i$ in the weak sense. Since these events occur almost surely, the proof is complete.
        \end{proof}

\begin{table}
\caption{\label{tb:pre} Existing research and our contributions}
\centering
\begin{tabular}{ccc}
   \toprule
   Preconditioner setting & Almost sure convergence & Restricted primal-dual gap \\
   \midrule
    $N_1=I; N_2=I$  & Already have in \cite{CPE} & Presented in this paper \\
   $N_1\succ I; N_2\succ I$ & Easy to prove with \cite{AFC}  & Presented in this paper\\
   $N_1 \succeq I,N_1\neq I; N_2=I$ & Presented in this paper & Presented in this paper\\
   \bottomrule
\end{tabular}
\end{table}
As shown in Table \ref{tb:pre}, our framework covers some results in existing research. The commonly used preconditioners for stochastic (and relaxed) PDR are the cases of $N_1 \succeq I, N_1\neq I$ and $N_2=I$, which brings out the degeneracy of $N_1-I$ and leads to the difficulties in this proof. For instance, the symmetric Gauss-Seidel preconditioner is widely used for an elliptic operator $T$ as in \eqref{N1:m} and $M-T=N_1-I$ is only positive semidefinite \cite[Proposition 2.13]{BSCC}. Next, let us turn to the convergence rate of the restricted primal-dual gap with respect to the transitional variables $(x^k_{\test},y^k_{\test})$. It is also a challenge due to the preconditioners and particularly because of the transitional variables, which are the main differences compared to the gap estimate reported in \cite{AFC}. Since the relations derived from \eqref{eq:e-delta} are concerning $u^k$, i.e., $(x^k,y^k,\bar x^k,\bar y^k)$, without any relations of transitional variables $(x^k_{\test},y^k_{\test})$, we introduce Lemma \ref{lemma:M_norm}, Lemma \ref{lemma:SK} and Lemma \ref{lemma:exp-sum} to address this issue.
\section{Sublinear convergence of stochastic RPDR for restricted primal-dual gap}\label{sec:sublinear:gap}
Taking full expectation 
in \eqref{eq:e-delta}, induction with  $k $ yields
\begin{equation}\label{eq:Delta}
\mathbb{E} [\Delta^{k} ] \leq \Delta^{0}=\|u^{0}-u^* \|^2_{\mathbf{M}_{p, \varrho_0}},\ \
\sum_{i=0}^{k}\frac{2-\rho_u}{\rho_u^2}\mathbb{E} [\|u^{i+1}-u^i \|^2_{\mathbf{M}_p}   ] \leq\sum_{i=0}^{k}\mathbb{E} [\|u^{i+1}-u^i \|^2_{(2I-I_{\varrho_k})I^{-2}_{\varrho_k}\mathbf{M}_p}   ] \leq \Delta^{0}.
\end{equation} 
Since \ref{sto-rpdr} involves $u^{k+1}_t = \mathcal{T}u^{k}$ in the iteration, let us begin with the following lemma which gives the control of $u^{k+1}_t$ by the weighted $\mathbf{M}_p$ norm of $u^k$, where $u^k$ is the iteration sequence of \ref{sto-rpdr}.  Applying \cite[Lemma 2.5]{BSCC} with $\mathbf{M}_p $ norm, we have the following estimate.
\begin{lemma}\label{lemma:M_norm}
	There exists a constant $C > 0$ such that for all $u^{1}, u^{2} \in \mathbb{U}$ and $\tilde u^1=\mathcal{T}u^1 $, $\tilde u^2=\mathcal{T}u^2 $ satisfies
	\begin{equation}\label{eq:mp}
	\|\tilde u^{1}-\tilde u^{2}\|^2 \leq C\|u^{1}-u^{2}\|_{\mathbf{M}_p}^2,
	\end{equation}
	where $\mT$ is the iteration operator of \ref{PDR}, the constant $C $ is only related to $N_1 $, $N_2 $, $\sigma $, $\tau $ and $\mathcal{K} $.
\end{lemma}
The following two lemmas can address the main challenges in $y^{k+1}_{\test} $ for the gap estimate.
\begin{lemma}\label{lemma:SK}
	Supposing $\delta \in(0,1)$ is a constant, given a non-decreasing sequence $ S_k $  and a bounded sequence $h_k$ satisfying 
	$0<S_k \leq h_k+(1-\delta)S_{k-1}$ for  $k\geq 1$,  we  
	have
	\begin{align}
	S_k\leq \sum_{i=0}^{k-1} (1-\delta)^i h_{k-i}+(1-\delta)^{k}S_0 .\label{eq:Sk}
	\end{align}
\end{lemma}	
\begin{proof}
		An induction of $S_k $ gives \eqref{eq:Sk}. Furthermore, if $h_k $ is bounded by $M $, we have $S_k\leq [1-(1-\delta)^{k}]M/{\delta} +(1-\delta)^{k}S_0 \leq  {M}/{\delta}    +(1-\delta)^{k}S_0 $.
\end{proof}

\begin{lemma}\label{lemma:exp-sum}
	The estimates of $\mathbb{E}\|y^{k+1}_{\test}-y^* \|^2 $ and $\mathbb{E}\sum_{k=0}^{K}\| y^{k+1}_{\test}-{y}^{k}_{\test} \|$ are as follows. First, with notations $ s^k $ and $\varphi $ defined in \eqref{eq:Ek}, we deduce
	
	\begin{itemize}
		\item [\emph{(i)}] 
		
		$\mathbb E \varphi(s_i^{k+1})=\mathbb{E}[\sum_{j=0}^{k}p_i(1-p_i)^j\varphi(\tilde s_i^{k+1-j})]+(1-p_i)^{k+1}\varphi(s_i^0),$
		\item [\emph{(ii)}]
		$\mathbb{E}\|s^{k+1}-s^k \|^2_{P^{-1}-I}=\mathbb{E}\|\tilde s^{k+1}-s^k \|^2_{I-P} =\mathbb{E}\|\tilde s^{k+1}-s^{k+1} \|^2.$
	\end{itemize} 
	Then, it follows 
	\begin{itemize}
		\item [{\emph{(iii)}}]
		$
		\mathbb E \| y^{k+1}_{\test}-y^* \|^2 \leq  \frac{(3C+3\tau)\rho_u}{\rho_l} \|u^0-u^*\|_{\mathbf{M}_p}^2+\| y^{0}_{\test}-y^* \|^2, 
		$ where $C  $ is as in Lemma \ref{lemma:M_norm}.

		\item [{\emph{(iv)}}] 
		$\mathbb{E}\sum_{k=1}^{K}\|\tilde y^{k+1}_{\test}-\tilde{y}^{k}_{\test} \|^2 \leq \frac{D\rho^2_u}{(2-\rho_u)\rho_l}\|u^0-u^*\|_{\mathbf{M}_p}^2$, where $ D>0$ is a positive constant.

		\item [{\emph{(v)}}] 
		$\mathbb{E}\sum_{k=0}^{K}\| y^{k+1}_{\test}-{y}^{k}_{\test} \|^2 \leq M_{K}$  with $M_K \leq M_{\infty}$ and $M_{\infty}$ being a positive constant.
		
	\end{itemize}
  \end{lemma}	
	\begin{proof} For (i), since $ \mathbb{E}\mathbb{E}^k[\cdot]=\mathbb{E}[\cdot]  $, a direct induction gives,
		\begin{align} 
		\mathbb E \varphi(s_i^{k+1})
		&=\mathbb E\mathbb E^k  \varphi(s_i^{k+1})
		=\mathbb E[  p_i\varphi(\tilde s_i^{k+1})+(1-p_i)\varphi( s_i^{k})]\notag\\
		&=\mathbb E  p_i\varphi(\tilde s_i^{k+1})+\mathbb E\mathbb E^{k-1}(1-p_i)\varphi( s_i^{k})\notag\\
		&=\mathbb E [ p_i\varphi(\tilde s_i^{k+1})+  p_i(1-p_i)\varphi(\tilde s_i^{k}) +(1-p_i)^2\varphi(s_i^{k-1}) ]=\cdots\notag\\
		&=\mathbb{E}[\sum_{j=0}^{k}p_i(1-p_i)^j\varphi(\tilde s_i^{k+1-j})]+(1-p_i)^{k+1}\varphi(s_i^0),
		\end{align}
		which leads to (i). 
		
		For (ii), the trick is to find the relations among those three items. That is, for every $k $, $\|\tilde s^{k+1}-s^k \|^2=\|s^{k+1}-s^k \|^2+\|\tilde s^{k+1}-s^{k+1} \|^2$ since $\|s_i^{k+1}-s_i^k \|^2$ will be zero except $i=i_k$ and $\|\tilde s_i^{k+1}-s_i^{k+1} \|^2$ will be zero except $i\neq i_k$.  Using $\mathbb{E}\|s^{k+1}-s^k \|^2=\mathbb{E}\|\tilde s^{k+1}-s^k \|_P^2 $, the statement (ii) follows. Now, let us prove $\mathbb{E}\| \tilde y^{k+1}_{\test}-y^* \|^2 $ is uniformly bounded for all $k$, based on \eqref{eq:sto-different} and the iteration formula of \ref{PDR}, we have
  \begin{align}
		\mathbb{E}\| \tilde y^{k+1}_{\test} -y^* \|^2
	& = \mathbb{E}\| \tilde y^{k+1}_{t,\test} -y^* \|^2 = \mathbb{E}\| \tby_t^{k+1} -\bar{y}^k   +\tilde y_t^{k+1}-y^* \|^2\notag\\
&\leq 3\mathbb{E} \| \tby_t^{k+1} -\bar y^*\|^2 
  +3\mathbb{E}\| \tilde y_t^{k+1}-y^* \|^2+3\mathbb{E}\| \bar{y}^k-\bar y^*    \|^2 \notag\\
&\leq 3  \mathbb{E} \| \mathcal{T}u^k-u^*\|^2+3\tau\rho_{k,y}\mathbb{E}\|\bar{y}^k-\bar y^*  \|_{\tau^{-1}\rho_{k,y}^{-1}\mathbf{P}^{-1}}^2\notag\\  
&\leq 3C  \mathbb{E} \|u^k-u^*\|_{\mathbf{M}_{p}}^2   +3\tau\rho_{u}\mathbb{E} \|u^{k} -u^*\|_{\mathbf{M}_{p,\varrho_k}}^2  \notag\\
&\leq 3C\rho_u  \mathbb{E} \|u^k-u^*\|_{\mathbf{M}_{p,\varrho_k}}^2   +3\tau\rho_{u}\mathbb{E} \|u^{k} -u^*\|_{\mathbf{M}_{p,\varrho_k}}^2  \notag\\
  &= (3C+3\tau)\rho_u\mathbb{E}  \Delta_k \notag  \leq (3C+3\tau) \rho_u \Delta^
0\leq \frac{(3C+3\tau)\rho_u}{\rho_l} \|u^0-u^*\|_{\mathbf{M}_p}^2
  \end{align}where the third inequality is deduced from Lemma \ref{lemma:M_norm} while the last one is from \eqref{eq:Delta}. 
		Applying (i) with $\varphi(\cdot)=\|\cdot-y^*\|^2 $ on $y^{k+1}_{\test} $ and noting $ 0<\sum_{j=0}^{k}p_i(1-p_i)^j< 1$, $0<(1-p_i)^{k+1}< 1 $, then $ 0 \prec \sum_{j=0}^{k}\mathbf{P}(I-\mathbf{P})^j\prec I$, $0\prec(I-\mathbf{P})^{k+1} \prec I $, we arrive at
		\begin{align}
		\mathbb E \| y^{k+1}_{\test}-y^* \|^2&=\mathbb{E}\sum_{j=0}^{k}\| \tilde y^{k+1-j}_{\test}-y^* \|_{\mathbf{P}(I-\mathbf{P})^j}^2+\| y^{0}_{\test}-y^* \|_{(I-\mathbf{P})^{k+1}}^2,\notag\\
		&\leq \frac{(3C+3\tau)\rho_u}{\rho_l} \|u^0-u^*\|_{\mathbf{M}_p}^2+\| y^{0}_{\test}-y^* \|^2\label{eq:exp-ytest}.
		\end{align}

		Now, let us turn to (iv). Remembering the relation between $\tilde y^{k+1}_{\test}$ and $\tilde y^{k+1}_{t,\test}$ as in \eqref{eq:sto-different}, the updates of \ref{PDR} and the non-expansiveness of the resolvent operator $(I+\tau\partial\mathcal G)^{-1}$, it can be readily checked that
		\begin{align}
		&\mathbb{E}\sum_{k=1}^K\|\tilde y^{k+1}_{\test}-\tilde{y}^{k}_{\test} \|^2 = \mathbb{E}\sum_{k=1}^K\|\tilde y^{k+1}_{t,\test}-\tilde{y}^{k}_{t,\test} \|^2\leq\mathbb{E}\sum_{k=1}^K\|(2\tilde y_t^{k+1}-\bar y^{k})- (2\tilde y_t^{k}-\bar y^{k-1})\|^2 \notag \\
		&= \mathbb{E}\sum_{k=1}^K\|A_1(N_1-I)(x^{k}-x^{k-1})+ A_1(\bar x^{k}-\bar x^{k-1}) +A_2(N_2-I)(y^{k}-y^{k-1}) \notag \\
		& \quad  +A_3(\bar y^{k}-\bar y^{k-1})  
		\|^2 \leq D \mathbb{E}\sum_{k=1}^K \|u^{k} -u^{k-1}\|_{\mathbf{M}_{p}}^2, \label{eq:ytestmkleq}
		\end{align}
		where $A_1=\tau A_2\mK N_1^{-1}$, $A_2=2(N_2+\tau\sigma\mK N_1^{-1}\mK^*)^{-1}$,  $A_3=A_2-I $, and
		\begin{align}
  D=4\max({\sigma(\|A_1\|^2\|N_1-I\|),\ {\sigma}\|A_1\|^2,\ \tau\|A_2\|^2\|N_2-I\|, \ {\tau}\|A_3\|^2}).
		\end{align}
Combining \eqref{eq:ytestmkleq} and \eqref{eq:Delta}, we obtain (iv). 	
				
		For (v), since $\mathbb{E} \| y^{k+1}_{\test}-{y}^{k}_{\test} \|^2= \mathbb{E}\|\tilde y^{k+1}_{\test}-{y}^{k}_{\test} \|_P^2$, we henceforth focus on $\mathbb{E}\|\tilde y^{k+1}_{\test}-{y}^{k}_{\test} \|^2$.  By the inequality  $\|a+b\|^2\leq(1+\epsilon)\|a\|^2+(1+{1}/{\epsilon})\|b\|^2 $, for any $\epsilon>0 $ and $K\geq1$, we arrive at
		\begin{align}
		& \sum_{k=1}^{K}\mathbb{E}\|\tilde y^{k+1}_{\test}-{y}^{k}_{\test} \|^2
		\leq\sum_{k=1}^{K}[\mathbb{E}(1+\epsilon)\|\tilde y^{k+1}_{\test}-\tilde{y}^{k}_{\test}\|^2+\mathbb{E}(1+\frac{1}{\epsilon})\| \tilde y^{k}_{\test}-y^{k}_{\test}\|^2  ]\notag \\
		 &=\sum_{k=1}^{K}\mathbb{E}(1+\epsilon)\|\tilde y^{k+1}_{\test}-\tilde{y}^{k}_{\test}\|^2+\sum_{k=1}^{K}\mathbb{E}(1+\frac{1}{\epsilon})\| \tilde y^{k}_{\test}-y^{k-1}_{\test}\|_{I-\mathbf{P}}^2  \label{eq:ytestepsilon} \\
		&=\sum_{k=1}^{K}\mathbb{E}(1+\epsilon)\|\tilde y^{k+1}_{\test}-\tilde{y}^{k}_{\test}\|^2  +(\sum_{k=1}^{K-1}\mathbb{E}(1+\frac{1}{\epsilon})\| \tilde y^{k+1}_{\test}-y^{k}_{\test}\|_{I-\mathbf{P}}^2+\mathbb{E}(1+{1}/{\epsilon})\| \tilde y^{1}_{\test}-y^{0}_{\test}\|_{I-\mathbf{P}}^2), \notag 
		\end{align}
		where the first equality is obtained by (ii). For the convenience of proofs, letting $\epsilon ={1}/{\underline{p}}  $ with $\underline{p} : = \min_{i}\{p_i, \ i=1, \ldots, n \}$ henceforth, then $(1+{1}/{\epsilon})(I-\mathbf{P})\leq (1+{1}/{\epsilon})(1-\underline{p})I=(1-\underline{p}^2)I$.
		
		For employing Lemma \ref{lemma:SK}, choosing $\delta=\underline{p}^2$ and $S_0=0$, we denote for $K\geq1 $,
		\[
		S_K:=\sum_{k=1}^{K}\mathbb{E}\|\tilde y^{k+1}_{\test}-{y}^{k}_{\test} \|^2, \ \   h_K:=\sum_{k=1}^{K}\mathbb{E}(1+\epsilon)\|\tilde y^{k+1}_{\test}-\tilde{y}^{k}_{\test}\|^2 +\mathbb{E}(1+{1}/{\epsilon})\| \tilde y^{1}_{\test}-y^{0}_{\test}\|_{I-\mathbf{P}}^2.
		\]
		The equation \eqref{eq:ytestepsilon} thus turns to
		\begin{align}
		S_K-(1-\delta)S_{K-1} \leq  h_K. \label{eq:SK-1}
		\end{align} 
		For $\mathbb{E}(1+{1}/{\epsilon})\| \tilde y^{1}_{\test}-y^{0}_{\test}\|_{I-\mathbf{P}}^2 $, we have
		\begin{align}
		&\mathbb{E}(1+{1}/{\epsilon})\| \tilde y^{1}_{\test}-y^{0}_{\test}\|_{I-\mathbf{P}}^2 \leq \mathbb{E}\|\tilde y^1_{\test}-y^0_{\test} \|_{(1-\delta)I}^2 =\|\tilde y^1_{\test}-y^0_{\test} \|^2_{(1-\underline{p}^2)I}.\label{eq:ytest1exact}
		\end{align}
		Then by (iv), $h_K$ is bounded, i.e.,
		\begin{align}
		h_K \leq  (1+\underline{p}^{-1})\frac{D\rho^2_u}{(2-\rho_u)\rho_l}\|u^0-u^*\|_{\mathbf{M}_p}^2+ \| \tilde y_{\test}^1- y^0_{\test}\|_{(1-\underline{p}^2)I}^2. \label{eq:hkbounded}
		\end{align}
		 Thus by \eqref{eq:SK-1}, \eqref{eq:hkbounded}, Lemma \ref{lemma:SK} and the fact $S_0=0$, we have 
		\begin{align}
		S_K&\leq \frac{1-(1-\underline{p}^2)^{K-1}}{\underline{p}^2} \left((1+\underline{p}^{-1})\frac{D\rho^2_u}{(2-\rho_u)\rho_l}\|u^0-u^*\|_{\mathbf{M}_p}^2+ \| \tilde y_{\test}^1- y^0_{\test}\|_{(1-\underline{p}^2)I}^2\right) \\
     &\leq \frac{1}{\underline{p}^2}\left((1+\underline{p}^{-1})\frac{D\rho^2_u}{(2-\rho_u)\rho_l}\|u^0-u^*\|_{\mathbf{M}_p}^2+ \| \tilde y_{\test}^1- y^0_{\test}\|_{(1-\underline{p}^2)I}^2\right).\notag
		\end{align}
		Hence
		\begin{align}
		&\mathbb{E}\sum_{k=0}^{K}\|\tilde y^{k+1}_{\test}-{y}^{k}_{\test} \|^2=\mathbb{E}\sum_{k=1}^{K}\|\tilde y^{k+1}_{\test}-{y}^{k}_{\test} \|^2+\mathbb{E}\|\tilde y^{1}_{\test}-{y}^{0}_{\test} \|^2\notag\\
		&\leq \frac{1-(1-\underline{p}^{2})^{K-1}}{\underline{p}^{2}}\left((1+\underline{p}^{-1})\frac{D\rho^2_u}{(2-\rho_u)\rho_l}\|u^0-u^*\|_{\mathbf{M}_p}^2 + \| \tilde y_{\test}^1- y^0_{\test}\|_{(1-\underline{p}^{2})I}^2\right)+\|\tilde y^{1}_{\test}-{y}^{0}_{\test} \|^2.
		\notag
		\end{align}
		Finally, with (ii), we derive
		\begin{align}
		\mathbb{E}&\sum_{k=0}^{K}\| y^{k+1}_{\test}-{y}^{k}_{\test} \|^2= \mathbb{E}\sum_{k=0}^{K}\|\tilde y^{k+1}_{\test}-{y}^{k}_{\test} \|_\mathbf{P}^2 \leq \mathbb{E}\sum_{k=0}^{K}\|\tilde y^{k+1}_{\test}-{y}^{k}_{\test} \|^2\notag\\
		&\leq \frac{1-(1-\underline{p}^2)^{K-1}}{\underline{p}^2}\left((1+1/\underline{p})\frac{D\rho^2_u}{(2-\rho_u)\rho_l}\|u^0-u^*\|_{\mathbf{M}_p}^2+ \| \tilde y_{\test}^1- y^0_{\test}\|_{(1-\underline{p}^2)I}^2\right)+\|\tilde y^{1}_{\test}-{y}^{0}_{\test} \|^2:=M_K,\label{eq:Eytestsum} \\
		&\leq {1}/{\underline{p}^2}\left((1+1/\underline{p})\frac{D\rho^2_u}{(2-\rho_u)\rho_l}\|u^0-u^*\|_{\mathbf{M}_p}^2+ \| \tilde y_{\test}^1- y^0_{\test}\|_{(1-\underline{p}^2)I}^2\right)+\|\tilde y^{1}_{\test}-{y}^{0}_{\test} \|^2:=M_{\infty}. \notag
		\end{align}
\end{proof}
Let us extend a lemma in \cite{AFC} originated from \cite{NJL} to positive semidefinite weights.
\begin{lemma}\label{lemma:V}
	With notations $ s^k $ defined in Lemma \ref{lem:almost}, picking up a point $ s^0_{\dagger}\in \mathbb Y$ at will, let us define an iteration sequence
	$
	v^{k+1}=s^{k}-\tilde{s}^{k+1}-\mathbf{P}^{-1} (s^{k}-s^{k+1} ) $
	with
	$
	{s}_{\dagger}^{k+1}={s}_{\dagger}^{k}-\mathbf{P} v^{k+1} $.  Introduce a  positive semidefinite, self-adjoint, and linear weight $\mathfrak D = \emph{Diag}[\mathfrak d_1, \ldots, \mathfrak d_n  ]$ with $\mathfrak d_i: \mathbb{Y}_i \rightarrow \mathbb{Y}_i$, and a non-increasing positive weight $ \Theta^k = \emph{Diag}[\vartheta_1^k I, \ldots,  \vartheta_n^k I  ]$ with $\vartheta_i^k I: \mathbb{Y}_i \rightarrow \mathbb{Y}_i$ and  $0< \vartheta_l \leq \vartheta_i^{k+1}  \leq \vartheta_i^{k} \leq  \vartheta_u $ for $i=1, \ldots, n$ and $k=0, \ldots,K-1$. Here $\vartheta_l$ and $\vartheta_u$ are positive constants.  We also assume that the weights $\Theta^k$, $\mathfrak D$ and the probability matrix   $\mathbf{P}$ can communicate with each other, i.e.,  $\mathbf{P} \Theta^k = \Theta^k \mathbf{P}$,  $\mathfrak D \Theta^k = \Theta^k\mathfrak D $ together with  $\mathfrak D \mathbf{P} = \mathbf{P}\mathfrak D$ as before. Then, it follows that for any $s \in\mathbb Y$
	\begin{align}
	&\sum_{k=0}^{K} \langle{s}_{\dagger}^{k}-s, v^{k+1} \rangle_{\mathfrak{D}\Theta^k } \leq \frac{1}{2}\|{s}_{\dagger}^{0}-s\|_{\mathfrak{D} \Theta^0 \mathbf{P}^{-1}}^{2}-\frac{1}{2}\|{s}_{\dagger}^{K+1}-s\|_{\mathfrak{D}\Theta^K \mathbf{P}^{-1}}^{2}+\sum_{k=0}^{K} \frac{1}{2}\|v^{k+1}\|_{\mathfrak{D} \Theta^k \mathbf{P}}^{2},\label{eq:daggery}
	\end{align}
	where the above inequality in \eqref{eq:daggery} becomes equality while $\Theta^k \equiv \Theta^0$. We further have  
	\begin{align}\mathbb{E} [\sum_{k=0}^{K} \frac{1}{2}\|v^{k+1}\|_{\mathfrak{D} \Theta^k \mathbf{P}}^{2} ] \leq \mathbb{E}[\sum_{k=0}^{K} \frac{1}{2}\|s^k-s^{k+1}\|_{\mathfrak{D}  \Theta^k  \mathbf{P}^{-1}}^{2} ].\label{eq:v-barv}
	\end{align}
	Moreover, $v^{k}$ and ${s}_{\dagger}^{k}$ are $\mathfrak{F}_{k}$ measurable and $\mathbb{E}^{k} [v^{k+1} ]=0$. It follows that 
	\begin{align}
	\mathbb{E}^{k}\langle {s}_{\dagger}^{k},v^{k+1}\rangle=0.\label{eq:inner=0}
	\end{align}
\end{lemma}	
	\begin{proof}
		By the definition of $s^{k+1}_{\dagger} $, we have the estimate,
		\begin{align}\label{eq:proof-dagger}
		\frac{1}{2}\|{s}_{\dagger}^{k+1}-s\|_{ \mathbf{P}^{-1}\mathfrak{D} \Theta^k}^{2} &=\frac{1}{2}\|{s}_{\dagger}^{k}-s\|_{\mathbf{P}^{-1}\mathfrak{D}\Theta^k}^{2}- \langle \mathbf{P} v^{k+1}, {s}_{\dagger}^{k}-s \rangle_{\mathbf{P}^{-1}\mathfrak{D}\Theta^k}+\frac{1}{2}\|\mathbf{P} v^{k+1}\|_{\mathbf{P}^{-1}\mathfrak{D}\Theta^k}^{2} \notag\\
		&=\frac{1}{2}\|{s}_{\dagger}^{k}-s\|_{\mathbf{P}^{-1}\mathfrak{D}\Theta^k}^{2}- \langle v^{k+1}, {s}_{\dagger}^{k}-s \rangle_{\mathfrak{D}\Theta^k}+\frac{1}{2}\|v^{k+1}\|_{\mathbf{P}\mathfrak{D}\Theta^k}^{2}.
		\end{align}
		Summing the equality \eqref{eq:proof-dagger}, we arrive at 
  \[
  \sum_{k=0}^{K} \langle{s}_{\dagger}^{k}-s, v^{k+1} \rangle_{\mathfrak{D}\Theta^k }  = \frac{1}{2} \sum_{k=0}^{K} \left[\|{s}_{\dagger}^{k}-s\|_{\mathbf{P}^{-1}\mathfrak{D}\Theta^k}^{2} - \|{s}_{\dagger}^{k+1}-s\|_{\mathbf{P}^{-1}\mathfrak{D}\Theta^k}^{2}\right] + \frac{1}{2} \sum_{k=0}^{K}\|v^{k+1}\|_{\mathbf{P}\mathfrak{D}\Theta^k}^{2}.
  \]
   Noting that $\Theta^{k} \succeq \Theta^{k+1}$ together with $\mathbf{P}\mathfrak{D}\Theta^{k} \succeq \mathbf{P}\mathfrak{D}\Theta^{k+1}$ for $k=0,\ldots, K-1$, we obtain \eqref{eq:daggery} from the above equation. For the second result, similar to  $\mathbb{E}^k\|v^{k+1}\|^2 = \mathbb{E}^k\|\mathbb{E}^k\xi-\xi\|^2 \leq \mathbb{E}^k\|\xi\|^2$, with $\xi^k = \mathbf{P}^{-1} (s^{k}-s^{k+1} )$ and $\mathbb{E}^k\xi = s^{k}-\tilde{s}^{k+1}$, we obtain
		\begin{align}
		&\mathbb{E} [\sum_{k=0}^{K} \|v^{k+1}\|_{\mathfrak{D} \mathbf{P}\Theta^k}^{2} ] =\sum_{k=0}^{K} \mathbb{E} [\mathbb{E}^{k} [\|v^{k+1}\|_{\mathfrak{D} \mathbf{P} \Theta^k}^{2} ] ]\notag \\
		&\leq \sum_{k=0}^{K}  \mathbb{E} [\mathbb{E}^{k} [\|\mathbf{P}^{-1} (s^{k+1}-s^{k} )\|_{\mathfrak{D} \mathbf{P}\Theta^k}^{2} ] ]=\sum_{k=0}^{K}  \mathbb{E} [\|s^{k+1}-s^{k}\|_{\mathfrak{D} \mathbf{P}^{-1}\Theta^k}^{2} ].
		\end{align}
		By direct calculation, $\mathbb{E}^{k} [v^{k+1} ]=0$. Since ${s}_{\dagger}^{k}$ is $\mathfrak{F}_{k}$ measurable, we conclude $\mathbb{E}^{k}\langle {s}_{\dagger}^{k},v^{k+1}\rangle=\langle {s}_{\dagger}^{k}, \mathbb{E}^{k}v^{k+1}\rangle=0$.
	\end{proof}
 
With the preceding lemmas, we now prove the sublinear convergence of \ref{sto-rpdr}. 
\begin{theorem}\label{thm:convergence-sdr}
	If $u=(x,y,\bar x,\bar y),$ with $ \bar x=x+\sigma\mathcal K^{\ast}y,$ $\bar y=y-\tau\mathcal Kx$. Let $u^k$ be the iteration sequence generated by \ref{sto-rpdr}, and let $z^k_{\test}=(x^k_{\test},y^k_{\test}) $ be the corresponding transitional variables. Denote $K_1: = K+1$, then the ergodic sequence $ z_{\test,K}=\frac{1}{K+1}\sum_{k=0}^K z^k_{\test} =\frac{1}{K_1}\sum_{k=0}^K z^k_{\test}$ converges with rate $\mathcal{O}({1}/{K})$ in an expected
	restricted primal-dual gap sense, i.e., for any $ z=(x,y) \in  \mathbb B_1 \times\mathbb B_2 =\mathbb B \subset \dom\mathcal{F}\times\dom\mathcal{G}$, a bounded domain that contains a saddle-point, it holds that 
	\begin{equation}\label{eq:stodr-gap-convergence}
	\mathbb{E}\mathfrak{G}_{z\in\mathbb B}(z_{\test,K})=\mathbb{E} [\sup_{z\in\mathbb B}\mH(z_{\test,K}, z)]\leq\frac{C_{\mathbb B}}{K_1}=\mathcal{O}(\frac{1}{K}),
	\end{equation}
	where the constant is given by 
	\begin{align}
	&C_{\mathbb B}=	 \sup_{z\in\mathbb{B}}     
	\frac{1}{2\rho_l}\Big[
	\|x^{0}-x\|_{\frac{N_1-I}{\sigma}}^2+\|y^{0}-y\|_{\frac{(N_2-I)\mathbf{P}^{-1}}{\tau}}^2+\|\bar x^{0}-x-\sigma\mK^*y\|^2_{\frac{I}{\sigma}}+\|\bar y^{0}-y+\tau\mK x\|^2_{\frac{\mathbf{P}^{-1}}{\tau}}\notag\\
	&+\|y^0_{\dagger} -y \|^2_{\frac{(N_2-I)\mathbf{P}^{-1}}{\tau}}+\|\bar y^0_{\dagger} -y+\tau\mK x \|^2_{\frac{\mathbf{P}^{-1}}{\tau}} \Big]+\sup_{z\in\mathbb{B}} \frac{1}{2}(\|{y}_{\dagger,\test}^{0}-\mathcal{K}x\|_{\mathbf{P}^{-1}}^{2}+\|\mathcal{K}x\|^2) \notag\\
&+\frac{\rho_u^2}{2(2-\rho_u)\rho^2_l}\| u^0-u^*\|^2_{\mathbf{M}_p}
 +\frac{(1-\underline{p})^2}{ \underline{p}^2}(\frac{(3C+3\tau)\rho_u}{\rho_l}\|u^0-u^*\|^2_{\mathbf{M}_p}+2\| y^{0}_{\test}-y^* \|^2)+\frac{1}{2\underline{p}}M_K\notag\\
	&+\mG_{\mathbf{P}^{{-1}}-I}(y^0_{\test}) -	\sum_{i=1}^{n}(\frac{1}{p_{i}}-1 )\mG_{i} (y_{i}^{*}) +	\sum_{i=1}^{n} (\frac{1}{p_{i}}-1 ) \|\mathcal K_{i} x^{*}\| \sqrt{\frac{(3C+3\tau)\rho_u}{\rho_l} \|u^0-u^*\|_{\mathbf{M}_p}^2+\| y^{0}_{\test,i}-y_i^* \|^2} . \notag 
	\end{align}
Here $ C$ is produced in Lemma \ref{lemma:M_norm}, $ M_K$ is from \eqref{eq:Eytestsum}, $(x^*,y^*) $ is a saddle point, and $y^0_{\dagger} $, $\bar y^0_{\dagger} $, $y^0_{\dagger,\test} $ are picked from $\mathbb{B}_2$ randomly. The same rate also holds for the Bregman distance $\mH(z_{\test,K}, (x^*,y^*) ) $. 
	\end{theorem}
	\begin{proof}
		 Using the following two identities,
		\begin{align*}
		&\mathcal{H} ( (x_{\test}^{k+1}, \tilde{y}_{\test}^{k+1}), z )=\mF (x_{\test}^{k+1} )+ \langle \mK x_{\test}^{k+1}, y \rangle-\mG(y)-\mF(x)- \langle \mK x, \tilde{y}_{\test}^{k+1} \rangle+\mG (\tilde{y}_{\test}^{k+1} ), \\
		&\mathcal{H} ( (x_{\test}^{k+1},{y}_{\test}^{k+1}), z )=\mF (x_{\test}^{k+1} )+ \langle \mK x_{\test}^{k+1}, y \rangle-\mG(y)-\mF(x)- \langle \mK x, {y}_{\test}^{k+1} \rangle+\mG ({y}_{\test}^{k+1} ), 
		\end{align*}
	we obtain,
		\begin{align}\label{eq:h-y} 
		\mH((x_{\test}^{k+1}, {y}_{\test}^{k+1}), z ) =\mathcal{H} ((x_{\test}^{k+1}, \tilde y_{\test}^{k+1}), z )+ \langle \mK x, \tilde y_{\test}^{k+1}-{y}_{\test}^{k+1} \rangle-\mG (\tilde{y}_{\test}^{k+1} )+\mG (y_{\test}^{k+1} ).
		\end{align}
		We now separate $ \mH((x_{\test}^{k+1}, {y}_{\test}^{k+1}), z )$ into three parts including $\mathcal{H} ((x_{\test}^{k+1}, \tilde y_{\test}^{k+1}), z ) $, $\langle \mK x, \tilde y_{\test}^{k+1}-{y}_{\test}^{k+1} \rangle $ and $\mG (y_{\test}^{k+1} ) -\mG (\tilde{y}_{\test}^{k+1} ) $. By the convexity of $\mathcal{H}$ and the property of supremum, we have 
		\begin{align}
		&\mathbb{E} [\sup_{z\in\mathbb B}\mH(z_{\test,K}, z)]\leq \mathbb{E} [\sup_{z\in\mathbb B}  \frac{1}{K_1} \sum_{k=0}^{K} \mH(z^K_{\test}, z)] 
		\leq \mathbb{E} [\sup_{z\in\mathbb B}  \frac{1}{K_1}\sum_{k=0}^{K} \mathcal{H} ((x_{\test}^{k+1}, \tilde y_{\test}^{k+1}), z )  ] \notag \\
		&+\mathbb{E} [\sup_{z\in\mathbb B} \frac{1}{K_1}\sum_{k=0}^{K} \langle \mK x, \tilde y_{\test}^{k+1}-{y}_{\test}^{k+1} \rangle   ]
		+\mathbb{E} [\sup_{z\in\mathbb B} \frac{1}{K_1}\sum_{k=0}^{K}(  \mG (y_{\test}^{k+1} ) -\mG (\tilde{y}_{\test}^{k+1} ) ) ]. \label{eq:lastver:convexity}
		\end{align}
		We will discuss these three parts one by one. First, denote by
		\begin{align}
		\mathbf{x}^k:=(x^k,\bar{x}^k),\quad  \mathbf{y}^k:=(y^k,\bar{y}^k),\quad\tilde{ \mathbf{y}}^k:=(\tilde y^k,\tby^k),\quad  \mathbf{x}:=(x,\bar{x}),\quad\mathbf{y}:=(y,\bar{y}),\label{eq:newy}
		\end{align}
		and the following diagonal weights $\mathbf w_x: =\diag[\tfrac1\sigma (N_1-I), \tfrac{1}{\sigma} I]$, $\mathbf w_y :=\diag[\tfrac1\tau (N_2-I), \tfrac{1}{\tau} I]$, $\mrx^k :=\diag[\tfrac1{\rho_k} I, \tfrac{1}{\rho_{k,x}} I]$, $\mry^k :=\diag[\tfrac1{\rho_k} I, \tfrac{1}{\rho_{k,y}} I]$, $\mrrx^k :=\diag[\frac{2-\rho_k}{\rho_k^2}I, \frac{2-\rho_{k,x}}{\rho_{k,x}^2} I]$, and $\mrry^k :=\diag[\frac{2-\rho_k}{\rho_k^2}I, \frac{2-\rho_{k,y}}{\rho_{k,y}^2} I]$. Therefore, $\|u^k-u\|^2_{\mathbf{M_p}}=\|\mathbf{x}^k-\mathbf{x}\|^2_{\mathbf w_x}+\|\mathbf{y}^k-\mathbf{y}\|^2_{\mathbf w_y\mathbf{P} }$, in which $\mathbf w_y\mathbf{P}:=\diag[\tfrac1\tau(N_2-I)\mathbf{P}, \tfrac1\tau \mathbf{P}] $. 
		For $\mathcal{H} ((x_{\test}^{k+1}, \tilde y_{\test}^{k+1}), z) $, we focus on the iteration $u^k \rightarrow\tilde{u}^{k+1} $. Combining Proposition  \ref{pro:prop-pdr} and \eqref{eq:same}, we obtain
		\begin{equation}\label{eq:h-tilde}
		\mathcal{H} ((x_{\test}^{k+1}, \tilde{y}_{\test}^{k+1}), z) \leq \Delta_{\mathbf x}^k+\Delta_{\mathbf y}^k,
		\end{equation}
		where $\Delta_{\mathbf x}$ and $\Delta_{\mathbf y}$ are defined as
		\begin{equation}\label{eq:delta}
		\begin{array}{c}		\Delta_{\mathbf x}^k:=
		\frac{1}{2} \Big(  \|\mathbf{x}^{k}-\mathbf{x}\|_{\mathbf w_x \mrx^k}^{2}
		-\|\mathbf{x}^{k+1}-\mathbf{x}\|_{\mathbf w_x \mrx^k}^{2}-\|\mathbf{x}^{k+1}-\mathbf{x}^{k}\|_{\mathbf w_x \mrrx^k}^{2} \Big), \\
		\Delta_{\mathbf y}^k:=\frac{1}{2} \Big( \|\mathbf{y}^{k}-\mathbf{y}\|_{\mathbf w_y \mry^k}^{2}
		-\|\tilde{\mathbf{y}}^{k+1}-\mathbf{y}\|_{\mathbf w_y \mry^k}^{2}-\|\tilde{ \mathbf{y}}^{k+1}-\mathbf{y}^{k}\|_{\mathbf w_y \mrry^k}^{2} \Big).\\
		\end{array}
		\end{equation} 
		Let us first reformulate $ \Delta_{\mathbf y}^k$ as follows
		\begin{align}\label{eq:deltay}
		\Delta_{\my}^k &=\frac{1}{2} \Big[-\| \tilde{\my}^{k+1}-\my^{k}\|_{\mwy \mrry^k}^{2}-\|\tilde{\my}^{k+1}-\my \|_{\mwy \mry^k}^{2}+\|\my^{k}-\ y\|_{\mwy \mry^k}^{2} \notag\\
		& + (\|\my^{k}-\my\|_{\mwy \mathbf{P}^{-1}\mry^k}^{2}-\|\my^{k+1}-\my\|_{\mwy \mathbf{P}^{-1}\mry^k}^{2} )- (\|\my^{k}-\my\|_{\mwy \mathbf{P}^{-1}\mry^k}^{2}-\|\my^{k+1}-\my\|_{\mwy \mathbf{P}^{-1}\mry^k}^{2} ) \Big] \notag\\
		&=-\frac{1}{2}\|\my^{k+1}-\my\|_{\mwy \mathbf{P}^{-1}\mry^k}^{2}+\frac{1}{2}\|\my^{k}-\my\|_{\mwy \mathbf{P}^{-1}\mry^k}^{2}-\frac{1}{2}\|\tilde{\my}^{k+1}-\my^{k}\|_{\mwy\mrry^k}^{2}+\epsilon_{\my}^{k},
		\end{align}
		where $\epsilon_{\my}^{k} $ is defined as follows with $\mv^{k+1}:=\my^{k}-\mty^{k+1}-\mathbf{P}^{-1} (\my^{k}-\my^{k+1} ) $ for integer $k \geq 0$,
		\begin{align}
		\epsilon_{\my}^{k} &=\frac{1}{2} \Big[\|\my^{k}-\my\|_{\mwy \mry^k}^{2}-\|\tilde{\my}^{k+1}-\my\|_{\mwy \mry^k}^{2}- (\|\my^{k}-\my\|_{\mwy \mathbf{P}^{-1}\mry^k}^{2}-\|\my^{k+1}-\my\|_{\mwy \mathbf{P}^{-1}\mry^k}^{2} ) \Big] \notag\\
		&=\frac{1}{2} \Big[\|\my^{k}\|_{\mwy  \mry^k}^{2}-\|\tilde{\my}^{k+1}\|_{\mwy \mry^k}^{2}- (\|\my^{k}\|_{\mwy \mathbf{P}^{-1}\mry^k}^{2}-\|\my^{k+1}\|_{\mwy \mathbf{P}^{-1}\mry^k}^{2} ) -2 \langle \my, \mv^{k+1} \rangle_{\mwy\mry^k} \Big].\label{eq:epsy}
		\end{align}
		Now, with \eqref{eq:h-tilde}, \eqref{eq:delta} and \eqref{eq:deltay}, we arrive at 
		\begin{align}
		\mathcal{H} ((x_{\test}^{k+1}, \tilde y_{\test}^{k+1}), z )  
		&\leq \frac{1}{2}\Big( \|\mathbf{x}^{k}-\mathbf{x}\|_{\mathbf w_x \mrx^k}^{2}
		-\|\mathbf{x}^{k+1}-\mathbf{x}\|_{\mathbf w_x \mrx^k}^{2}-\|\mathbf{x}^{k+1}-\mathbf{x}^{k}\|_{\mathbf w_x \mrrx^k}^{2}\notag\\
		&+\|\my^{k}-\my\|_{\mwy \mathbf{P}^{-1}\mry^k}^{2}-\|\my^{k+1}-\my\|_{\mwy \mathbf{P}^{-1}\mry^k}^{2}-\|\tilde{\my}^{k+1}-\my^{k}\|_{\mwy \mrry^k}^{2}\Big)+\epsilon_{\my}^{k}\notag\\
		&= \frac{1}{2}\Big(\|u^{k}-u\|^2_{ \mathbf{M}_{p,\varrho_k}}-\|u^{k+1}-u\|^2_{  \mathbf{M}_{p,\varrho_k}}-\|u^{k+1}-u^k\|^2_{(2I-I_{\varrho_k})I^{-2}_{\varrho_k}\mathbf{M}_p}\Big)\notag\\
		&- \langle \my, \mv^{k+1} \rangle_{\mwy \mry^k}+\frac{1}{2}\Big(\| \my^{k+1}-\my^k \|^2_{\mwy \mathbf{P}^{-1}\mrry^k}-\|\mty^{k+1}-\my^k \|^2_{\mwy \mrry^k}\notag\\
		&+\|\my^{k}\|_{\mwy \mry^k}^{2}-\|\mty^{k+1}\|_{\mwy \mry^k}^{2}-\|\my^{k}\|_{\mwy \mathbf{P}^{-1}\mry^k}^{2}+\|\my^{k+1}\|_{\mwy \mathbf{P}^{-1}\mry^k}^{2} \Big),\label{eq:h-y-2}
		\end{align}
		where the equality is by adding and subtracting $\frac{1}{2}\|\my^{k+1} -\my^k\|^2_{\mwy \mathbf{P}^{-1}\mrry^k} $ and the definition of $\epsilon_{\my}^{k}$ in \eqref{eq:epsy}. In \eqref{eq:h-y-2}, dropping  $-\|u^{k+1}-u^k\|^2_{\mathbf{M}_p (2I-I_{\varrho_k})I^{-2}_{\varrho_k}}$ and summing  over $k=0,\ldots,K $ yields,
		\begin{align}\label{eq:h-y-sum}
		\sum^{K}_{k=0}\mathcal{H} ((x_{\test}^{k+1}, \tilde y_{\test}^{k+1}), z )
		&\leq -\frac{1}{2}\|u^{K+1}-u\|^2_{\mathbf{M}_{p,\varrho_K} }+\frac{1}{2}\|u^{0}-u\|^2_{\mathbf{M}_{p,\varrho_0}}+\mathcal E^K +\sum_{k=0}^{K} [  -\langle \my,\mv^{k+1} \rangle_{\mwy \mry^k}],
		\end{align}
		where we use $\mathbf{M}_{p,\varrho_{k+1}} \leq \mathbf{M}_{p,\varrho_k} $ for all $k=0,\dots,K$. $	\mathcal E^K$ is defined as follows
		\begin{align*}
		\frac{1}{2}\sum_{k=0}^{K}  \Big[  &(
		\| \my^{k+1}-\my^k \|^2_{\mwy \mathbf{P}^{-1} \mrry^k}-\|\mty^{k+1}-\my^k \|^2_{\mwy \mrry^k} )  \\
		&+  (\|\my^{k}\|_{\mwy  \mry^k}^{2}-\|\mty^{k+1}\|_{\mwy  \mry^k}^{2}-\|\my^{k}\|_{\mwy \mathbf{P}^{-1} \mry^k}^{2}+\|\my^{k+1}\|_{\mwy \mathbf{P}^{-1} \mry^k}^{2} )\Big].
		\end{align*}
		Noting there are no items involving $(x,y)$ in $\mathcal E^K$, we can take expectation directly without concerning taking the supremum over $(x,y)$. With \eqref{eq:trival} and the fact $\mathbb{E} [\mathbb{E}^{k}[\cdot] ]=\mathbb{E}[\cdot]$, we get
		\begin{align}\label{eq:esum-e}
		\mE	\mathcal E^K/K_1	&=	\frac{1}{2K_1} \mE	 \sum_{k=0}^{K}  \Big[ \| \my^{k+1}-\my^k \|^2_{\mwy \mathbf{P}^{-1}\mrry^k} -\mE^k\|\my^{k+1}-\my^k \|^2_{\mwy \mathbf{P}^{-1}\mrry^k} 
		\notag\\
		&+ \mE^k(\|\my^{k}\|_{\mwy \mathbf{P}^{-1}\mry^k}^{2}-\|\my^{k+1}\|_{\mwy \mathbf{P}^{-1}\mry^k}^{2})-(\|\my^{k}\|_{\mwy \mathbf{P}^{-1}\mry^k}^{2}-\|\my^{k+1}\|_{\mwy \mathbf{P}^{-1}\mry^k}^{2} )\Big]=0. 
		\end{align}
		For the sum of inner products in \eqref{eq:h-y-sum}, introducing $\my^{k+1}_{\dagger}=\my^{k}_{\dagger}-\mathbf{P}\mv^{k+1}$ and selecting a $\my^0_{\dagger} = (y^0_{\dagger},\bar y^0_{\dagger})\in\mathbb{B}_2\times\mathbb{B}_2$
		randomly, we have 
		\[
		- \langle \my, \mv^{k+1} \rangle_{\mwy \mry^k}= \langle \my_{\dagger}^{k}-\my, \mv^{k+1} \rangle_{\mwy \mry^k}- \langle \my_{\dagger}^{k}, \mv^{k+1} \rangle_{\mwy \mry^k}. 
		\]
		In \eqref{eq:h-y-sum}, dropping $-\frac{1}{2}\|u^{K+1}-u\|^2_{\mathbf{M}_p \mathbf{R}^k}$, dividing both sides by $K_1$, taking supremum with respect to $z $, and employing \eqref{eq:daggery} on the above inner product, we arrive at
		\begin{equation}\label{eq:h-y-sup}
		\begin{aligned}
		\mE [ 	&\sup_{z\in\mathbb{B}}\frac{1}{K_1}\sum_{k=0}^{K}\mathcal{H} ((x_{\test}^{k+1}, \tilde y_{\test}^{k+1}), z ) ] \leq   \sup_{z\in\mathbb{B}}     
		\frac{1}{2K_1} [  \|u^{0}-u\|^2_{\mathbf{M}_{p,\varrho_0}} +\|\my^0_{\dagger} -\my \|^2_{\mwy \mathbf{P}^{-1} \mry^0} ]\\
		&-\mE\|\my^{K+1}_{\dagger} -\my \|^2_{\mwy \mathbf{P}^{-1} \mry^k}+ \mE\sum_{k=0}^{K}\frac{1}{K_1}  [  \frac{1}{2}\|\mv^{k+1} \|^2_{\mwy \mathbf{P} \mry^k} - \langle \my_{\dagger}^{k},  \mv^{k+1} \rangle_{\mwy \mry^k}  ].\\
		\end{aligned}
		\end{equation}
		 By \eqref{eq:v-barv}, non-decreasing of $\rho_{k}$, $\rho_{k,y}$ and \eqref{eq:Delta} we have, \begin{equation}\label{eq:v-delta}
		 \begin{aligned}
		 &\mE\sum_{k=0}^{K}  \|\mv^{k+1} \|^2_{\mwy \mathbf{P}\mry^k} \leq \mathbb{E}[\sum_{k=0}^{K} \|\my^k-\my^{k+1}\|_{\mwy \mathbf{P}^{-1}\mry^k}^{2}]\leq \mathbb{E}[\sum_{k=0}^{K} \frac{1}{\rho_l}\|\my^k-\my^{k+1}\|_{\mwy \mathbf{P}^{-1}}^{2}]\\&\leq\mathbb{E}[\sum_{k=0}^{K} \frac{1}{\rho_l}\|u^k-u^{k+1}\|_{\mathbf{M}_p}^{2}]\leq \frac{\rho_u^2}{(2-\rho_u)\rho_l}\Delta^0 = \frac{\rho_u^2}{(2-\rho_u)\rho_l^2}\| u^0-u^*\|^2_{\mathbf{M}_p}.
   \end{aligned}    
		 \end{equation}
Thus in \eqref{eq:h-y-sup}, dropping $-\mE\|\my^{K+1}_{\dagger} -\my \|^2_{\mwy \mathbf{P}^{-1}\mry^k}$, together with \eqref{eq:inner=0} and \eqref{eq:v-delta}, we can obtain 
  
		\begin{equation}\label{eq:h-y-sup2}
  \begin{aligned}
		&\mE [ 	\sup_{z\in\mathbb{B}}\frac{1}{K_1}\sum_{k=0}^{K}\mathcal{H} ((x_{\test}^{k+1}, \tilde y_{\test}^{k+1}), z ) ] \leq        
		\frac{1}{2K_1\rho_l} \big[ \sup_{z\in\mathbb{B}}(\|u^{0}-u\|^2_{\mathbf{M}_{p}}+\|y^0_{\dagger} -y \|^2_{\frac{(N_2-I)\mathbf{P}^{-1}}{\tau}}\\
  &+\|\bar y^0_{\dagger} -y+\tau\mK x \|^2_{\frac{\mathbf{P}^{-1}}{\tau}} )+\frac{\rho_u^2}{(2-\rho_u)\rho_l}\| u^0-u^*\|^2_{\mathbf{M}_p} \big]. 
  & 
  \end{aligned}
		\end{equation}
		Second, let us give the upper bounds for $	\mE [ 	\sup_{z\in\mathbb{B}}\frac{1}{K_1}\sum_{k=0}^{K}\langle \mK x, \tilde y_{\test}^{k+1}-{y}_{\test}^{k+1} \rangle] $. Denote by $v_{\test}^{k+1}=y_{\test}^{k}-\tilde y_{\test}^{k+1}-\mathbf{P}^{-1} (y_{\test}^{k}-y_{\test}^{k+1} ) $, $y^{k+1}_{\dagger,\test}=y^{k}_{\dagger,\test}-\mathbf{P}v_{\test}^{k+1}$ and select a $y^0_{\dagger,\test}\in\mathbb{B}_2$ randomly. By the definition of $y^{k}_{\dagger,\test} $, $v^{k}_{\test} $ and \eqref{eq:daggery}, we obtain
		\begin{align}
		&\sum_{k=0}^{K}\langle \mK x, \tilde y_{\test}^{k+1}-{y}_{\test}^{k+1} \rangle\notag\\
		=&\sum_{k=0}^{K}-\langle \mK x, y_{\test}^{k}-\tilde{y}_{\test}^{k+1} -\mathbf{P}^{-1}(y_{\test}^{k}-{y}_{\test}^{k+1} )\rangle +\sum_{k=0}^{K}\langle \mK x,  (\mathbf{P}^{-1}-I)(y_{\test}^{k+1}-{y}_{\test}^{k} )\rangle\notag\\
		=&\sum_{k=0}^{K}\langle y^{k}_{\dagger,\test}-\mK x, v^{k+1}_{\test}\rangle-\sum_{k=0}^{K}\langle y^{k}_{\dagger,\test},v^{k+1}_{\test} \rangle +\langle \mK x,  (\mathbf{P}^{-1}-I)(y_{\test}^{K+1}-{y}_{\test}^{0} )\rangle\notag\\
		\leq &\frac{1}{2}\|{y}_{\dagger,\test}^{0}-\mathcal{K}x\|_{\mathbf{P}^{-1}}^{2}-\frac{1}{2}\|{y}_{\dagger,\test}^{K+1}-\mathcal{K}x\|_{\mathbf{P}^{-1}}^{2}+\sum_{k=0}^{K} \frac{1}{2}\|v_{\test}^{k+1}\|_{\mathbf{P}}^{2}-\sum_{k=0}^{K}\langle y^{k}_{\dagger,\test},v^{k+1}_{\test} \rangle\notag\\
		&+\langle \mK x,  (\mathbf{P}^{-1}-I)(y_{\test}^{K+1}-{y}_{\test}^{0} )\rangle\label{eq:kxyt}.
		\end{align}
		For the last term in \eqref{eq:kxyt}, we have \begin{align}
		\langle \mK x, (\mathbf{P}^{-1}-I ) (y_{\test}^{K+1}-y_{\test}^{0} ) \rangle \leq \frac{1}{2}\|\mK x\|^{2}+\frac{(1-\underline{p})^2}{2 \underline{p}^2}\|y_{\test}^{K+1}-y_{\test}^{0}\|^{2}\label{eq:ytestk+1}.
		\end{align}
		In \eqref{eq:kxyt}, dropping $-\frac{1}{2}\|{y}_{\dagger,\test}^{K+1}-\mathcal{K}x\|^{2}$ and by \eqref{eq:ytestk+1},  \eqref{eq:inner=0}, we get 
		\begin{align}
		\mE [ 	\sup_{z\in\mathbb{B}}\frac{1}{K_1}
		\sum_{k=0}^{K}\langle \mK x, \tilde y_{\test}^{k+1}-{y}_{\test}^{k+1} \rangle]&\leq\sup_{z\in\mathbb{B}}[\frac{1}{2K_1}\|{y}_{\dagger,\test}^{0}-\mathcal{K}x\|_{\mathbf{P}^{-1}}^{2}+\frac{1}{2K_1}\|\mathcal{K}x\|^2] \notag\\
		&+\frac{1}{2K_1}\mathbb{E}\sum_{k=0}^{K}\|v_{\test}^{k+1}\|_{\mathbf{P}}^{2}+\mathbb{E}\frac{(1-\underline{p})^2}{2 \underline{p}^2K_1}\|y_{\test}^{K+1}-y_{\test}^{0}\|^{2}\label{eq:inner}.
		\end{align}
		For the last item in \eqref{eq:inner}, inserting $y^*$ and employing Lemma \ref{lemma:exp-sum}, we get, 
		\begin{align}
		\mathbb{E}\frac{(1-\underline{p})^2}{2 \underline{p}^2K_1}\|y_{\test}^{K+1}-y_{\test}^{0}\|^{2}&\leq
		\mathbb{E}\frac{(1-\underline{p})^2}{ \underline{p}^2K_1}\|y_{\test}^{K+1}-y^*\|^{2}+\frac{(1-\underline{p})^2}{ \underline{p}^2K_1}\|y^*-y_{\test}^{0}\|^{2}\notag\\
		&\leq \frac{(1-\underline{p})^2}{ \underline{p}^2K_1}(\frac{(3C+3\tau)\rho_u}{\rho_l} \|u^0-u^*\|_{\mathbf{M}_p}^2+2\| y^{0}_{\test}-y^* \|^2).
		\end{align}
		Noting $\mathbb{E}\sum_{k=0}^{K}[\|v_{\test}^{k+1}\|_{\mathbf{P}}^{2}] \leq \mathbb{E}\sum_{k=0}^{K} \|y_{\test}^k-y_{\test}^{k+1}\|_{\mathbf{P}^{-1}}^2$ by \eqref{eq:v-barv} together with \eqref{eq:Eytestsum}, finally, we arrive at
		\begin{align}
		&\mE [ 	\sup_{z\in\mathbb{B}}\frac{1}{K_1}
		\sum_{k=0}^{K}\langle \mK x, \tilde y_{\test}^{k+1}-{y}_{\test}^{k+1} \rangle]\leq\sup_{z\in\mathbb{B}}[\frac{1}{2K_1}\|{y}_{\dagger,\test}^{0}-\mathcal{K}x\|_{ \mathbf{P}^{-1}}^{2}+\frac{1}{2K_1}\|\mathcal{K}x\|^2] \notag\\
		&+ \frac{(1-\underline{p})^2}{ \underline{p}^2K_1}(\frac{(3C+3\tau)\rho_u}{\rho_l} \|u^0-u^*\|_{\mathbf{M}_p}^2+2\| y^{0}_{\test}-y^* \|^2)+  \frac{1}{2\underline{p}K_1}M_K \label{eq:inner2}.
		\end{align}
		Third,  we show $	\mE [ 	\sup_{z\in\mathbb{B}}\frac{1}{K_1} \sum_{k=0}^K(\mG (y_{\test}^{k+1} ) -\mG (\tilde{y}_{\test}^{k+1} ))] $ is also bounded. In fact,
		\begin{align}
		\mE [ 	\sup_{z\in\mathbb{B}}\frac{1}{K_1} \sum_{k=0}^K(\mG (y_{\test}^{k+1} ) -\mG (\tilde{y}_{\test}^{k+1} ))]=	\mE \frac{1}{K_1} \sum_{k=0}^K[ \mG (y_{\test}^{k+1} ) -\mG (\tilde{y}_{\test}^{k+1} )].\label{eq:G}
		\end{align}
		Employing \eqref{eq:trival2} on $\mG(\tilde y^{k+1}_{\test}) $, we get $\mG(\tilde y^{k+1}_{\test})=\mathbb{E}^k\mG_{\mathbf{P}^{-1}}(y^{k+1}_{\test})-\mG_{\mathbf{P}^{-1}-I}(y^{k}_{\test})$. Here $ \mG_{\mathbf{P}^{-1}}(y^{k+1}_{\test}) := \sum^n_{i=1}1/p_i \mG_i(y^{k+1}_{\test,i})$ and $ \mG_{\mathbf{P}^{-1}-I}(y^{k+1}_{\test}) := \sum^n_{i=1}(1/p_i-1)\mG_i(y^{k+1}_{\test,i})$. Thus, using the fact $\mathbb{E}\mathbb{E}[\cdot]=\mathbb{E}[\cdot]$ again, we arrive at
		\begin{align}
		&\mE \frac{1}{K_1} \sum_{k=0}^K[ \mG (y_{\test}^{k+1} ) -\mG (\tilde{y}_{\test}^{k+1} )]=  \frac{1}{K_1} \sum_{k=0}^K\mE[ \mG (y_{\test}^{k+1} ) -\mathbb{E}^k\mG_{\mathbf{P}^{-1}}(y^{k+1}_{\test})+\mG_{\mathbf{P}^{-1}-I}(y^{k}_{\test})]\notag\\
		&=\frac{1}{K_1} \sum_{k=0}^K[-\mathbb{E} \mG_{\mathbf{P}^{-1}-I}(y^{k+1}_{\test})+\mathbb{E}\mG_{\mathbf{P}^{-1}-I}(y^{k}_{\test})]=\frac{1}{K_1}( \mG_{\mathbf{P}^{-1}-I}(y^{0}_{\test}) -\mathbb{E}\mG_{\mathbf{P}^{-1}-I}(y^{{ K}+1}_{\test})).\label{eq:g2}
		\end{align}
		Since $\mK_{i} x^{*} \in \partial \mG_{i} (y^*_{i} )$, one has
		\[
		\mG_i (y_{\test,i}^{K+1} ) \geq \mG_i (y_{i}^{*} )+ \langle \mK_{i} x^{*}, y_{\test,i}^{K+1}-y_{i}^{*} \rangle \geq \mG_i (y_{i}^{*} )-\|\mathcal{K}_{i} x^{*}\|\|y_{\test,i}^{K+1}-y_{i}^{*}\|.
		\]
		We thus obtain
		\begin{align}\label{eq:gytest}
		\mathbb{E}[\mG_{\mathbf{P}^{-1}-I} (y_{\test}^{K+1} )& ] =\sum_{i=1}^{n} (\frac{1}{p_{i}}-1 ) \mathbb{E} [\mG_{i} (y_{\test,i}^{K+1} ) ]  \geq 	\sum_{i=1}^{n} (\frac{1}{p_{i}}-1 ) (\mG_{i} (y_{i}^{*} )-\|\mathcal K_{i} x^{*}\| \mathbb{E} [\|y_{\test,i}^{K+1}-y_{i}^{*}\|  ] ) \notag\\
		& \geq 	\sum_{i=1}^{n} (\frac{1}{p_{i}}-1 ) (\mG_{i} (y_{i}^{*} )-\|\mathcal K_{i} x^{*}\| \sqrt{\mathbb{E} [\|y_{\test,i}^{K+1}-y_i^{*}\|^2  ]} ) \notag\\
		& \geq 	\sum_{i=1}^{n} (\frac{1}{p_{i}}-1 ) (\mG_{i} (y_{i}^{*} )-\|\mathcal K_{i} x^{*}\| \sqrt{\frac{(3C+3\tau)\rho_u}{\rho_l} \|u^0-u^*\|_{\mathbf{M}_p}^2+\| y^{0}_{\test,i}-y^*_i \|^2} ),	
		\end{align}
		where the second inequality is by Jensen's inequality and the third inequality is by Lemma \ref{lemma:exp-sum}.
		We thus have 
		\begin{align}
		\mE \frac{1}{K_1} \sum_{k=0}^K[ \mG (y_{\test}^{k+1} ) -\mG (\tilde{y}_{\test}^{k+1} )]&\leq\frac{1}{K_1} \Big(\mG_{\mathbf{P}^{-1}-I}(y^{0}_{\test})-	\sum_{i=1}^{n}(\frac{1}{p_{i}}-1 )\mG_{i} (y_{i}^{*}) \notag\\&+	\sum_{i=1}^{n} (\frac{1}{p_{i}}-1 ) \|\mathcal K_{i} x^{*}\| \sqrt{\frac{(3C+3\tau)\rho_u}{\rho_l} \|u^0-u^*\|_{\mathbf{M}_p}^2+\| y^{0}_{\test,i}-y^*_i \|^2} \Big).
		\label{eq:g22}
		\end{align}
		Finally, with \eqref{eq:lastver:convexity} and combining \eqref{eq:h-y-sup2}, \eqref{eq:inner2} and \eqref{eq:g22}, we get \eqref{eq:stodr-gap-convergence} with the constant $C_{\mathbb B}$.
		For the Bregman distance, replacing $z=(x,y)$ by $(x^*,y^*)$, one can obtain the same sublinear convergence rate for the ergodic sequences.
	\end{proof}	
In the next theorem, we investigate the convergence rate of expected primal error when Lipschitz smoothness is satisfied for $\mG$.
\begin{theorem}\label{thm:loss-srpdr}
     For the ergodic sequence $ z_{\test,K}=(x_{\test,K},y_{\test,K})$ defined in Theorem \ref{thm:convergence-sdr}, $K_1 = K+1$ and any $ z=(x,y) \in  \dom\mathcal{F}\times\dom\mathcal{G}$, if \(\mG^*\) is \(\beta\)-Lipschitz continuous with respect to the norm \(\|\cdot\|_{\frac{I}{\tau}}\), and \(y^0\), \((\bar{y}^0 + \tau \mK x^*)\), \(y^0_{\dagger}\), \((\bar{y}^0_{\dagger} + \tau \mK x^*) \in \dom \mG\), and if \(y^*_{\dagger} \in \dom \mathcal{G}\), such that \(\sigma \mK^* y^*_{\dagger} = \bar{x}^0 - x^*\), then the expected primal error also exhibits a convergence rate of \(\mathcal{O}(1/K)\), i.e,
\begin{equation} \label{eq:primalerror}
   \mathbb{E}[ \mathcal{F}(x_{\test,K})+
\mG^*(\mK x_{\test,K}) - \mathcal{F}(x^*)- \mG^*(\mK x^*)]\leq\frac{C_{p,2}}{K_1} = \mathcal{O}(1/K),
\end{equation} where
 \begin{align}
       &C_{p,2} = C_{p}  + \frac{1}{2\rho_l}( \|x^{0}-x^*\|_{\frac{N_1-I}{\sigma}}^2+4\sigma\tau\beta^2\|\mK^*\|^2 + 8\|N_2-I\|\beta^2/\underline{p} +8\beta^2/\underline{p}  )\notag \\
& \quad \qquad  +\frac{1}{2}(\|{y}_{\dagger,\test}^{0}-\mathcal{K}x^*\|_{\mathbf{P}^{-1}}^{2}+\|\mathcal{K}x^*\|^2), \notag  \ \text{with} \\
&C_{p}:=\frac{\rho_u^2}{2(2-\rho_u)\rho^2_l}\| u^0-u^*\|^2_{\mathbf{M}_p}
 +\frac{(1-\underline{p})^2}{ \underline{p}^2}(\frac{(3C+3\tau)\rho_u}{\rho_l}\|u^0-u^*\|^2_{\mathbf{M}_p}+2\| y^{0}_{\test}-y^* \|^2)+\frac{1}{2\underline{p}}M_K\notag\\
&+\mG_{\mathbf{P}^{{-1}}-I}(y^0_{\test}) -	\sum_{i=1}^{n}(\frac{1}{p_{i}}-1 )\mG_{i} (y_{i}^{*}) +	\sum_{i=1}^{n} (\frac{1}{p_{i}}-1 ) \|\mathcal K_{i} x^{*}\| \sqrt{\frac{(3C+3\tau)\rho_u}{\rho_l} \|u^0-u^*\|_{\mathbf{M}_p}^2+\| y^{0}_{\test,i}-y_i^* \|^2} . \notag 
	\end{align}
 
\begin{proof}
     Let us first define a new distance function $\bar\mH(z_{\test,K}, z)$,
    \begin{align}\notag
        &\bar\mH(z_{\test,K}, z) := \mH(z_{\test,K}, z) - \frac{1}{2\rho_l K_1}(
	\|x^{0}-x\|_{\frac{N_1-I}{\sigma}}^2+\|y^{0}-y\|_{\frac{(N_2-I)\mathbf{P}^{-1}}{\tau}}^2+\|\bar x^{0}-x-\sigma\mK^*y\|^2_{\frac{I}{\sigma}}\notag\\
	&+\|\bar y^{0}-y+\tau\mK x\|^2_{\frac{\mathbf{P}^{-1}}{\tau}}+\|y^0_{\dagger} -y \|^2_{\frac{(N_2-I)\mathbf{P}^{-1}}{\tau}}+\|\bar y^0_{\dagger} -y+\tau\mK x \|^2_{\frac{\mathbf{P}^{-1}}{\tau}} )- \frac{1}{2K_1}(\|{y}_{\dagger,\test}^{0}-\mathcal{K}x\|_{\mathbf{P}^{-1}}^{2}+\|\mathcal{K}x\|^2),\notag
    \end{align}
where $ y^0_{\dagger} $, $ \bar y^0_{\dagger} $ and ${y}_{\dagger,\test}^{0}$ are given as in Theorem \ref{thm:convergence-sdr}. Following the same proof of Theorem \ref{thm:convergence-sdr},  we can drive, $ \mathbb{E} \sup_{z\in\mathbb B} \bar \mH(z_{\test,K}, z) \leq\frac{C_{p}}{K_1}$. Then substituting $z = (x^*,y)$ for any $y\in\mathbb{B}_2$ in $\bar \mH(z_{\test,K}, z)$, we derive $ \mathbb{E} \sup_{y\in\mathbb B_2} \bar \mH(z_{\test,K}, (x^*,y)) \leq\frac{C_{p}}{K_1}$. By the definition of $y^*_{\dagger} $, we arrive at $\|\bar x^{0}-x^*-\sigma\mK^*y\|^2_{\frac{I}{\sigma}} \leq \sigma\|\mK^*\|^2\|y^*_{\dagger} -  y\|^2  $. With \cite[Corollary 17.19]{HBPL2} and noting that $\mG^*_i$ is $\beta$-Lipschitz continuous on the norm $\|\cdot\|_{\frac{I}{\tau}}$, we have $\|y^*_{\dagger} -  y\|^2 \leq 4\tau\beta^2 $, 
 $\|y^{0}-y\|_{\frac{(N_2-I)\mathbf{P}^{-1}}{\tau}}^2 \leq 4\|N_2-I\|\beta^2/\underline{p}$, $ \|\bar y^{0}-y+\tau\mK x^*\|^2_{\frac{\mathbf{P}^{-1}}{\tau}}\leq 4\beta^2/\underline{p}$,    $\|y^0_{\dagger} -y \|^2_{\frac{(N_2-I)\mathbf{P}^{-1}}{\tau}} 
 \leq 4\|N_2-I\|\beta^2/\underline{p}$ and $\|\bar y^0_{\dagger} -y+\tau\mK x^* \|^2_{\frac{\mathbf{P}^{-1}}{\tau}}\leq 4\beta^2/\underline{p}$. We thus obtain
 \begin{align}
&\mathbb{E} \sup_{y\in\mathbb B_2}  \mH(z_{\test,K}, (x^*,y))\notag \\
&= \mathbb{E} \sup_{y\in\mathbb B_2} [\mathcal{F}(x_{\test,K})+
\langle \mathcal{K}x_{\test,K},y\rangle-\mathcal{G}(y) - \mathcal{F}(x^*)-
\langle \mathcal{K}x^*,y_{\test,K}\rangle+\mathcal{G}(y_{\test,K})]\leq\frac{C_{p,2}}{K_1}. \label{eq:primalerror-1} 
 \end{align} 
By Lipschitzness of $\mG^*$, we pick $y \in \partial \mG^*(\mK x_{\test,K})\neq \mathbf{0}$ s.t., $\langle \mathcal{K}x_{\test,K},y\rangle-\mathcal{G}(y) = \mG^*(\mK x_{\test,K}) $, and by Fenchel-Young inequality, $\mathcal{G}(y_{\test,K}) - \langle \mathcal{K}x^*,y_{\test,K}\rangle \geq - \mG^*(\mK x^*)$. Employing these inequalities in \eqref{eq:primalerror-1}, we derive \eqref{eq:primalerror}.
\end{proof}   
\end{theorem}

\section{Stochastic PDR for quadratic problems (SRPDRQ)}\label{sec:sto:pdrq}
In this section, we will focus on the following relaxed and stochastic PDRQ method directly for compactness
\begin{equation} \label{sto-pdrq} \tag {SRPDRQ}
\left\{\begin{array}{ll}
\tilde x^{k+1}=\left(N_{1}+\sigma Q +\mathcal{K}^{*} N_{2}^{-1} \mathcal{K}\right)^{-1}\left(N_{1} x^{k}+\sigma f-\sigma \mathcal{K}^{*} N_{2}^{-1}\left(\left(N_{2}-I\right) y^{k}+\bar{y}^{k}\right)\right)\\
x^{k+1} = (1-\rho_k)x^k+\rho_k \tilde{x}_t^{k+1}, \\
y^{k+1}= \left\{
\begin{array}{ll}
(1-\rho_k)y_i^k+\rho_k \tilde{y}^{k+1}_{t,i}, \quad  \tilde{y}^{k+1}_{t,i} :=N_{2,i}^{-1}[ (N_{2,i}-I)y_i^k+\bar y_i^k+\tau\mathcal{K}_i \tilde x^{x+1} ] & i=i_k\\
y^k_i &\forall  i \neq i_k\\
\end{array} \right. \\
y_{\test}^{k+1}= \left\{
\begin{array}{ll}
(I+\tau\partial\mathcal G_i)^{-1}[2\tilde{y}^{k+1}_{t,i}-\bar y^k_i]  & i= i_k,  \\
y_{\test,i}^{k} & \forall i \neq i_k,\\
\end{array}  \right.\\
\bar y^{k+1}= \left\{
\begin{array}{ll}
\bar y_i^k+\rho_{k,y}(y_{\test,i}^{k+1}-\tilde{y}_{t,i}^{k+1}) & i= i_k,  \\
\bar y^k_i & \forall i \neq i_k, 
\end{array}  \right.\\
\end{array}
\right.
\end{equation}
where $N_1\succeq 0 $, $N_2\succeq I $ and $N_2 $ is separable, i.e., $N_2y=(N_{2,1}y_1, \ldots,N_{2,n}y_n)^T$. Here $\rho_k$,  $\rho_{k,y}$ are non-decreasing and have lower and upper limit $ \rho_l, \rho_u$, and $0< \rho_l< \rho_u< 2 $ which is similar as the \ref{sto-rpdr}. It can be observed that the updates for \(\bar{x}^{k+1}\) are not required in the \ref{sto-pdrq} iterations, unlike in \ref{sto-rpdr}. As a result, \ref{sto-pdrq} can be more efficient in terms of iteration time.

To clarify, we refer to the general \ref{sto-pdrq} as the stochastic fully relaxed preconditioned Douglas-Rachford method for quadratic problems (f-SRPDRQ) when $\rho_k \neq 1$ and $\rho_{k,y} \neq 1$. While when $\rho_{k} = 1$, and $\rho_{k,y} = 1$, we call it the stochastic preconditioned Douglas-Rachford method for quadratic problems (SPDRQ).

If we set $N_2=I$ and $N_1$ as in \eqref{N1:m1}, we derive the preconditioned iteration \eqref{preconditioner:from:Q}. Define $x_{\test}^{k+1}=x^{k+1}$ and $\mT_Q $ as in \eqref{eq:pdrq:T}, we do one-step deterministic update based on the $u^k$ and get $\tilde u^{k+1} $, i.e., 
\begin{equation} \label{eq:tildeu:after:updrq}
\tilde u^{k+1}:= \mathcal{T}_Qu^k, \quad \text{with} \  \ \tilde u^{k+1}= (\tilde x^{k+1},\tilde y^{k+1}, \tby^{k+1}),
\end{equation}
where $x^{k+1}=\tilde x^{k+1} $ and $\bar x^{k+1}=\tbx^{k+1} $. Noting $I_{q,\varrho_k}:=\textrm {Diag}[\rho_k,\rho_k,\rho_{k,y}]$  as before, with Proposition \ref {pro:prop-pdrq} and Lemma \ref{lemma:fixpdrq}, we can derive the following lemma and theorem similarly.
\begin{lemma}\label{lem:almostpdrq} For any $(x,y)\in \mathbb{B}_1\times \mathbb{B}_2 $, a bounded domain, $\bar x=x+\sigma\mathcal{K}^{\ast} y$ and $\bar y= y-\tau\mathcal{K}x $, denoting $\mathbf{M}_{Qp}:=\diag[\tfrac1\sigma N_1, \tfrac1\tau (N_2-I)\mathbf{P}^{-1}, \tfrac1\tau \mathbf{P}^{-1}]$, it holds that
	\begin{align}\label{eq:H-tildeypdrq}
	\mH ((x_{\test}^{k+1},\tilde y^{k+1}_{\test}), (x,y)) &\leq\frac{1}{2}\Big( \|  u^k- u\|^2_{\mathbf{M}_{Qp}I_{q,\varrho_k}^{-1}}-\mathbb{E}^k\| u^{k+1}-u\|^2_{\mathbf{M}_{Qp} I_{q,\varrho_k}^{-1}} \\
  &-\mathbb{E}^k\|u^{k+1}-u^k\|^2_{\mathbf{M}_{Qp}(2I-I_{q,\varrho_k})I_{q,\varrho_k}^{-2}}\Big). \notag
	\end{align}
\end{lemma}
Denote $\mathcal{Z}^* $ as the set for the saddle points of \eqref{eq:saddle-point-sum} and $\mathbb{F} \subset \mathbb{U}$ as the set of the fixed points for $\mathcal{T}_{Q}$.
\begin{theorem}\label{thm:almostsurepdrq}
	The iteration sequence $u^k=(x^k,y^k, \bar y^k)$ produced by \ref{sto-pdrq} converges weakly  to an $\mathbb{F}$-valued random variable $u^{**} =(x^{**},y^{**}, \bar y^{**}) $ almost surely, and $(x^{**}, y^{**})$ is a $\mathcal{Z}^*$-valued random variable with $Qx^{**}+\mathcal{K}^*y^{**}=f $ and ${y}^{**}-\tau \mathcal{K} {x}^{**}=\bar{y}^{**}$. Moreover the sequence $(x^{k+1}_{\test} $, $y^{k+1}_{\test}) $ also converges weakly to the $\mathcal{Z}^*$-valued random variable $(x^{**},y^{**})$ almost surely.
\end{theorem}
For the gap estimate, we give the following two lemmas which are similar to the SRPDR case. 

\begin{lemma}\label{lemma:M_normpdrq}
	There exists a constant $C_Q > 0$ such that for all $u^{1}, u^{2} \in \mathbb{U}_Q$ and $\tilde u^1=\mathcal{T}_Qu^1 $, $\tilde u^2=\mathcal{T}_Qu^2$ and  it leads to
	\begin{equation}\label{eq:mppdrq}
	\|\tilde u^{1}-\tilde u^{2}\|^2 \leq C_Q\|u^{1}-u^{2}\|_{\mathbf{M}_{Qp}}^2
	\end{equation}
	where $C_Q $ is only related to $N_1 $, $N_2 $, $f $, $\sigma $, $\tau $ and $\mathcal{K} $.
\end{lemma}
\begin{lemma}\label{lemma:ytestpdrq}
	Similar to Lemma \ref{lemma:exp-sum}, we give the following estimates of $\mathbb{E}\|y^{k+1}_{\test}-y^* \|^2 $ and  $\mathbb{E}\sum_{k=0}^{K}\| y^{k+1}_{\test}-{y}^{k}_{\test} \|$.
	\begin{itemize}
		\item [\emph{(i)}]
		$
		\mathbb E \| y^{k+1}_{\test}-y^* \|^2
		<
		\frac{(3C_Q+3\tau)\rho_u}{\rho_l}\| u^0-u^*\|^2_{\mathbf{M}_{Q_p}}+\| y^{0}_{\test}-y^* \|^2, 
		$ where $C_Q  $ is given in Lemma \ref{lemma:M_normpdrq}.
		\item [\emph{(ii)}] 
		$\mathbb{E}\sum_{k=1}^{K}\|\tilde y^{k+1}_{\test}-\tilde{y}^{k}_{\test} \|^2 <\frac{D_Q\rho_u^2}{(2-\rho_u)\rho_l}\| u^0-u^*\|^2_{\mathbf{M}_{Q_p}}$, where
		$$D_Q=3\max(\sigma\|N_1\|\|A_1\|^2,\  \tau \|N_2-I\|\|A_2\|^2,\ \tau\|A_3\|^2)$$
		and  $A_2=2(N_2+\tau\sigma\mK (N_1+\sigma Q)^{-1}\mK^*)^{-1}$, $A_1=A_2\tau\mK(N_1+\sigma Q)^{-1}$, $A_3=A_2-I $.
		\item [\emph{(iii)}] 
		$\mathbb{E}\sum_{k=0}^{K}\| y^{k+1}_{\test}-{y}^{k}_{\test} \|^2 <M_{Q,K}$, where
		\[
		M_{Q,K}=\frac{1-(1-\underline{p}^2)^{K-1}}{\underline{p}^2}((1+1/\underline{p})\frac{D_Q\rho_u^2}{(2-\rho_u)\rho_l}\| u^0-u^*\|^2_{\mathbf{M}_{Q_p}}+ \| \tilde y_{\test}^1- y^0_{\test}\|_{1-\underline{p}^2}^2)+\|\tilde y^{1}_{\test}-{y}^{0}_{\test} \|^2. 
		\]
	\end{itemize}
\end{lemma}

Concerning the restricted gap estimate for \ref{sto-pdrq}, we have the following similar theorem. 
\begin{theorem}\label{thm:convergence-spdrq}
	If $u=(x,y,\bar y),$ with  $\bar y=y-\tau\mathcal Kx$, let $u^k$ be the iteration sequence generated by \ref{sto-pdrq}, and let $z^k_{\test}=(x^k,y^k_{\test}) $ be the corresponding transitional variables. Denote $K_1: = K+1$, then the ergodic sequence $ z_{\test,K}=\frac{1}{K+1}\sum_{k=0}^K z^k_{\test} = \frac{1}{K_1}\sum_{k=0}^K z^k_{\test} $ converges with rate $\mathcal{O}({1}/{K})$ in an expected
	restricted primal-dual gap sense, i.e., for any $ z=(x,y) \in  \mathbb B_1 \times\mathbb B_2 =\mathbb B $, a bounded domain, it holds that
	\begin{equation}\label{eq:stodr-gap-convergence2}
	\mathbb{E}\mathfrak{G}_{z\in\mathbb B}(z_{\test,K})=\mathbb{E} [\sup_{z\in\mathbb B}\mH(z_{\test,K}, z)]\leq\frac{C_{\mathbb B}}{K_1} =\mathcal{O}({1}/{K}),
	\end{equation}
	where the constant $C_{\mathbb B}$ is a positive constant depending on $u^0$, $\mK$, $\mathbf{P}$, $\rho_l$, $\rho_u$,  $\mathbf{M}_{Q_p}$, $u^*$, $z$, and $C_Q$, $ M_{Q,K}$.
Here $ C_Q$ is produced from Lemma \ref{lemma:M_normpdrq}, and $ M_{Q,K}$ is from Lemma \ref{lemma:ytestpdrq}. Therefore the same rate holds for the Bregman distance $\mH(z_{\test,K}, (x^*,y^*)) $. 
\end{theorem}





\begin{remark}\label{remark:sto-block}
Theorem \ref{thm:almostsure}, \ref{thm:convergence-sdr}, \ref{thm:almostsurepdrq}, and \ref{thm:convergence-spdrq} can also cover the situation of random sets of indices to be selected and updated. Sample the dual variable $y = (y_1,\dots,y_n)$ properly \cite{sample1,sample2,CERS} into multiple blocks $\mathbb{S}_1$, $\mathbb{S}_2, \dots,\mathbb{S}_b$, such that $\mathbb{S}_1 \cup\mathbb{S}_2 \cup\dots \cup \mathbb{S}_b = \{y_1, y_2,\dots,y_n\}$ and $\mathbb{S}_i \cap \mathbb{S}_j  = \emptyset$ for any $i,j\in\{1,\dots,b \}$, $i\neq j$. We denote the
probability of selecting a block $i\in\{1,\ldots,b\} $ by $p_i$ where we assume $0<p_{i}<1$ for each $i$ and $\sum_{i=1}^{b} p_{i}=1$. Following the same routine, we can derive the same results on almost sure convergence and sublinear convergence of restricted primal-dual gap. 
\end{remark}

\section{Numerical Experiments}\label{sec:numer} In this section, we will illustrate two experiments, one is a deblurring problem, and the other is a Binary classification problem with real data.
\subsection{TGV-Deblur Problem}\label{subsec:tgv}
In the first experiments, we focus on the following image deblurring problem \eqref{eq:TGV-KL} based on the Kullback-Leibler (KL)  divergence  and the anisotropic second-order Total Generalized Variation  (TGV) \cite{BPK} regularization
\begin{equation}\label{eq:TGV-KL}
\min _{u \in U} \mF(u)+\mG_{1}^{*} (\mK_1 u)+\TGV_{\alpha}(u)=	\min _{u \in U, w \in V} \mF(u)+\mG_{1}^{*} (\mK_1 u)+\alpha_{1}\|\nabla u-w\|_{1}+\alpha_{0}\|\mathcal{E} w\|_1.
\end{equation}
A similar denoising problem in multi-modal electron tomography with KL divergence and TGV regularization is studied in \cite{RGMGK}. Here $U$ and $V$ are finite-dimensional Hilbert spaces and $w=(w^1,w^2)$. The two-dimensional problem \eqref{eq:TGV-KL} is standard and all definitions along with discretizations will be made clear afterwards. Let us first reformulate the problem \eqref{eq:TGV-KL} as the saddle-point problem.  $\mG_{1}^{*}$ is the appropriate data term measure with blurred data given by the following KL divergence \cite[Equation (9)]{RGMGK}, which is a variant of \cite[Equation (13)]{kl2}. 
\begin{equation}\label{eq:bredies_kl}
	\mathcal{G}^{*}_{1}(s)=\sum_{i=1}^{N} \begin{cases}
		s_{i}-b_{i} \log  s_i
		& \text { if }  s_{i} \geq 0,\\
		\infty & \text { else.}
	\end{cases}
\end{equation}
where $b_i$ is the blurred image and we use the convention $0\log0=0$ and $-a\log(0)=\infty$ for $a>0$.  It can be readily checked that its Fenchel conjugate is
\begin{equation}
	\mathcal{G}_{1}(z)=\sum_{i=1}^{N} \begin{cases}
		\frac{-b_jz_{j}}{1-z_j}-\frac{b_j}{1-z_j}
		+b_j \log  ( \frac{b_j}{1-z_j}  ),
		& \text { if } z_{j} \leq 1 \text { and } (b_{j}=0 \text { or } z_{j}<1 ), \\
		\infty & \text { else, }
	\end{cases}
\end{equation}
and its resolvent is $
	(I+\tau \partial \mathcal{G}_1)^{-1}(z)=\frac{1}{2} (z+1-\sqrt{ (z-1 )^{2}+ 4 \tau b} )$.
Here $\mF(u)= \mI_{\{0 \leq u\leq 1\}}(u)$ is a box constraint. It forces
$ u $ to be non-negative by the requirement of KL divergence and to have an upper bound $ 1 $. The equivalent saddle-point problem of \eqref{eq:TGV-KL} is as follows
\begin{equation}\label{eq:saddle:pd:tgv}
\begin{aligned}
\min_{u \in U, w \in V } \max _{s \in U, p \in V, q \in W}\langle \mK_1u, s\rangle_U + \langle\nabla u-w, p\rangle_{V}+\langle\mathcal{E} w, q\rangle_{W}  +\mF(u)\\
-\mG_{1}(s) -\sum_{i=1}^{2}\mI_{ \{\|p^i\|_{\infty} \leq \alpha_1 \}}(p^i)
-\sum_{j=1}^{4}\mI_{ \{\|q^j\|_{\infty} \leq \alpha_0 \}}(q^j).
\end{aligned}
\end{equation}
Let us reformulate \eqref{eq:saddle:pd:tgv} into \eqref{eq:saddle-point-sum}. We have the data $\mF(x)=\mF(u)= \mI_{\{0 \leq u\leq 1\}}(u)$ with $x=(u,w)$  and  $\mathcal{G}^{*}_{i}( y_i) = \alpha_1 \| y_i \|_{1},$ $i = 2,3 $, $ \mathcal{G}^{*}_{i}( y_i) = \alpha_0 \| y_i \|_{1}, $ $i = 4,5,6,7 $.  $\mG_i $ is the Fenchel conjugate of $\mG_i^*$. The dual variables $y=(y_1,\ldots,y_7)=(s,p,q)$ with $p=(p^1,p^2) \in V:=U^2 $ and $q=(q^1,q^2,q^3,q^4) \in W:=U^4$ is thus obtained. The linear operator $\mK : X \rightarrow Y $ as in \eqref{eq:saddle-point-sum} is
\begin{equation}
\mK= \left(\begin{array}{cc}
\mK_1 & 0 \\
\nabla & -I \\
0 & \mathcal{E}
\end{array} \right), \quad X = U \times V, \quad Y=U \times V \times W,
\end{equation} 
where $\mK_1 = \kappa  \ast$ is a convolution operator with kernel $\kappa$ for deblurring problems.  

We then go to the details of the discrete setting including the operators and spaces. In \eqref{eq:TGV-KL}, $\mathcal{E} w =(\nabla w+\nabla w^{T} ) / 2$ is a discrete symmetrized derivative. $\nabla$ is a discrete gradient operator with finite forward differences, i.e.,
$
(\nabla u)=(\partial_{x}^{+} u,\partial_{y}^{+} u)^T
$ \cite{CP}. For discretization, set image domain $\Omega \subset \mathbf{Z}^{2}$ as the discretized grid
$\Omega= \{(i, j) \mid i, j \in \mathbb{N}, 0 \leq i \leq d_1-1,0 \leq j \leq d_2-1 \}$,
where $d_1, d_{2}$ are the image dimensions \cite{BSCC}. The image space is defined as $U=\{u: \Omega  \rightarrow \mathbf{R}\}$ with the standard $L^{2}$ scalar product. 
For $ t \in [1, \infty)$, let us introduce the following norms: $  \|u\|_{\infty}=\max _{(i, j) \in \Omega} |u_{i, j} |$,
$$
\begin{array}{l}
\|u\|_{t}= (\sum_{(i, j) \in \Omega} |u_{i, j} |^{t} )^{1 / t},  \ \ 
\|p\|_{1}= \sum_{l=1}^2\sum_{(i, j) \in \Omega} |p_{i, j}^{l}|,  \ \ \|q\|_{1}= \sum_{l=1}^4\sum_{(i, j) \in \Omega} |q_{i, j}^{l}|.
\end{array}
$$
The adjoint operator of  $\nabla:  U  \rightarrow V$ is  $\nabla^{*}=-$ div (e.g., $\operatorname{div}p = \partial_{x}^{-}p^1 + \partial_{y}^{-}p^2$) with homogeneous Neumann and Dirichlet
boundary conditions, respectively and we refer to \cite{CP} for details. 
Next we describe the discretization of the anisotropic symmetric derivative $\mathcal{E} $ in the space $W= \{u: \Omega  \rightarrow S^{2 \times 2} \}
$, which leads to the anisotropic second-order TGV which makes it possible to separate dual variables. We identify elements $q \in W$ and employ the scalar product
$
\langle q, q^{\prime} \rangle_{W}= \langle q^{1}, (q^{\prime} )^{1} \rangle_{U}+ \langle q^{2}, (q^{\prime} )^{2} \rangle_{U}+ \langle q^{3}, (q^{\prime} )^{3} \rangle_{U}+ \langle q^{4}, (q^{\prime} )^{4} \rangle_{U}$.
The symmetrized derivative can be defined as follows
\begin{equation}\label{eq:anisoTGV}
\begin{aligned} \mathcal{E} w = \left[\begin{array}{cc}\partial_{x}^{+} w^{1} & \frac{1}{2} (\partial_{y}^{+} w^{1}+\partial_{x}^{+} w^{2} ) \\ \frac{1}{2} (\partial_{y}^{+} w^{1}+\partial_{x}^{+} w^{2} ) & \partial_{y}^{+} w^{2}\end{array}\right]  = \left[\begin{array}{c}\partial_{x}^{+} w^{1} \\ \partial_{y}^{+} w^{2} \\ \frac{1}{2} (\partial_{y}^{+} w^{1}+\partial_{x}^{+} w^{2} )\\
\frac{1}{2} (\partial_{y}^{+} w^{1}+\partial_{x}^{+} w^{2} )
\end{array} \right] \end{aligned}
\end{equation}
where the second equation has to be understood in terms of the identification $W=U^{4}$. Consequently, the negative adjoint realizes a discrete negative divergence operator according to $\langle\mathcal{E} w, q\rangle_{W}=-\langle w, \operatorname{div} q\rangle_{V}$ for all $w \in V, q \in$
$W$.  With $\Delta w= (\Delta w^{1}, \Delta w^{2} ),$ we can express the operator $\mathcal{E}^{*} \mathcal{E} w=-\frac{1}{2} \Delta w-\frac{1}{2} \delta w$. Here $\operatorname{div} q$ and $\delta w $ are gives as  
\[		
\operatorname{div} q:= \left[\begin{array}{l}\partial_{x}^{-} q^{1}+\frac{1}{2}\partial_{y}^{-} q^{3} +\frac{1}{2}\partial_{y}^{-} q^{4}\\ \frac{1}	{2}\partial_{x}^{-} q^{3}+\frac{1}	{2}\partial_{x}^{-} q^{4}+\partial_{y}^{-} q^{2}\end{array} \right],  \quad \delta w :=   \left[\begin{array}{l}\partial_{x}^{-} \partial_{x}^{+}w^{1}+\partial_{y}^{-}\partial_{x}^{+}w^{2} \\ \partial_{y}^{-} \partial_{y}^{+}w^{2}+\partial_{x}^{-}\partial_{y}^{+}w^{1}\end{array} \right].
\]
Moreover, we confirm that $q^3=q^4$ and only $6$ dual variables will be randomly updated.

\subsubsection{The preconditioner}
Setting $N_2=I $, our goal is to find efficient preconditioners to solve $Tx^{k+1}= b^k $, with $x^{k+1}=(u^{k+1},w^{k+1})^{T} $ and  $b^k=\bar{x}^k-\sigma\mK^*\bar{y}^k=(b_u^k,b_{w}^k)^{T}$ and
\begin{equation}\label{eq:T}
T=\left[\begin{array}{cc} \sigma\tau \mK_1^*\mK_1+I-\sigma\tau\Delta & \sigma\tau \mathrm{div} \\ -\sigma\tau \nabla & (1+\sigma\tau) I+\sigma\tau\mathcal{E}^{*} \mathcal{E}\end{array}\right].
\end{equation}
$T$ is challenging for preconditioning due to $\mathcal{E}^{*} \mathcal{E}$. Fortunately, it can be readily checked that
\begin{equation}\label{eq:T'}
T^{\prime}=\left[\begin{array}{cc}
(1+\sigma\tau) I-\sigma\tau \Delta & \sigma\tau \operatorname{div} \\
-\sigma\tau \nabla & (1+\sigma\tau) I-\sigma\tau \Delta
\end{array}\right]
\end{equation}
satisfies $T^{\prime} -T \succeq 0 $ \cite{BS}. It is clear that feasible preconditioners $M$ for $T^{\prime} $ are also feasible for $T$ (Since $M \succeq T'$, we have $M \succeq T$). Consequently, we obtain
\begin{equation}\label{eq:pre-update-x}
x^{k+1}=x^{k}+M^{-1}\left(b^{k}+\left(T^{\prime}-T\right) x^{k}-T^{\prime} x^{k}\right),
\end{equation}
which is a preconditioned iteration for $T'$. Now we adjust the right-hand side as follows.
\begin{equation}\label{eq:pre-update-b}
\begin{aligned}
\left(b^{\prime}\right)^{k} &:=b^{k}+\left(T^{\prime}-T\right) x^{k}
&=\left[\begin{array}{c}
b_{u}^{k}+(\sigma\tau I -\sigma\tau \mK_1^*\mK_1)u^k \\
b_{w}^{k}+\frac{\sigma\tau}{2}\left(\delta-\Delta\right) w^{k}
\end{array}\right].
\end{aligned}
\end{equation}
Now solving $Tx^{k+1} = b^k$ by preconditioned iteration: $x^{k+1} = x^{k} + M^{-1}(b^k-Tx^k)$ is transformed into a new preconditioned iteration: $x^{k+1} = x^{k} + M^{-1}((b^{\prime})^k-T^{\prime}x^k)$, which is much easier to compute and implement.
We refer to \cite[Section 5.2]{BS} for the block symmetric Gauss--Seidel preconditioners.

\subsubsection{Parameters}
Here let us set $U = \mathbb{R}^{d_{1}\times d_{2} }$ with $d_{1}=512$, $d_{2}=357$ and  $\mathcal{K}_{1} = k \ast$ is a convolution kernel with  $40\times40$ motion blur.
It can be checked that the Bregman distance for the ergodic sequence is bounded for a bounded domain containing a saddle point.  
\par
 Set the sampling to be uniform, i.e., $p_i=\frac{1}{n},\, n=6$. The regularization parameters are  $\alpha_0=10^{-4}$ and $\alpha_1=5\times10^{-5}$. The random iteration sequence is determined by random seed 5. The initial value of all algorithms is set to be zero vector, except the initial primal update term of $u$ in \eqref{eq:saddle:pd:tgv} is set to be the blurred image. We mainly compared with the SPDHG as in \cite{CERS} and we refer to \cite{FercoqBianchi} for the case stochastic updates both on the primal and dual variables. The step sizes for compared algorithms are as follows:\par
\setlength{\hangindent}{2em}
\ \ $\bullet$ 
PDHG, SPDHG: Primal-Dual Hybrid Gradient algorithm (PDHG) \cite{CP} and stochastic PDHG \cite{CERS}.
We use constant step size $\tau=\sigma=1/\|\mathcal{K} \|$, $\|\mathcal{K} \|=\sqrt {14} $ for PDHG and $\tau=0.25$, $ \sigma=0.12$ for SPDHG, which is much faster than the step size in \cite{CERS}. \par

\setlength{\hangindent}{2em}
\ \ $\bullet$ \ref{PDR}, \ref{sto-pdr}: One-step block symmetric Red-Black Gauss--Seidel is the preconditioner for the elliptic system equation. The step sizes are chosen as $\tau=0.1,\, \sigma=5$.\par

\ \ $\bullet$ p-RPDR (f-RPDR), p-SRPDR (f-SRPDR): partial (fully) relaxed PDR, partial (fully) stochastic relaxed PDR. Again, the preconditioner and the step sizes are the same as PDR. The relaxation parameters are all set to 1.9.\par

Now, let us turn to the numerical ergodic Bregman distance.
It can be checked that $x^{k+1}_{\test}$, $s^{k+1}_{\test}$ and $y^{k+1}_{\test,i},$ $i=\{2,\ldots,7\}$ produced by \ref{sto-pdr}  satisfy the constraints of $\mF$, $\mG_1$ and $\mG_{i},$ $i=\{2,\ldots,7\}$ correspondingly. Still, by the convexity of these functions, the ergodic sequence   $(x^{k+1}_{\test,K},y^{k+1}_{\test,K}) $ also satisfies the corresponding constraints. We thus conclude that the Bregman distance of ergodic sequence is finite and can be written as $ \langle \mK x_{\test,K},y^{*} \rangle-\mathcal{G}_1(s^{*}) - \langle \mK x^*,y_{\test,K} \rangle+\mathcal{G}_1(s_{\test,K})$.
Here $(x^*,y^*) $ is the approximate saddle point which is calculated for more than $10^5$ epochs.

\begin{figure}[!htb]
	\graphicspath{{fig//}}
	\begin{center}
		\subfloat[Convergence rate of normalized ergodic Bregman distance concerning epochs.]
		{\includegraphics[width=.22\paperwidth]{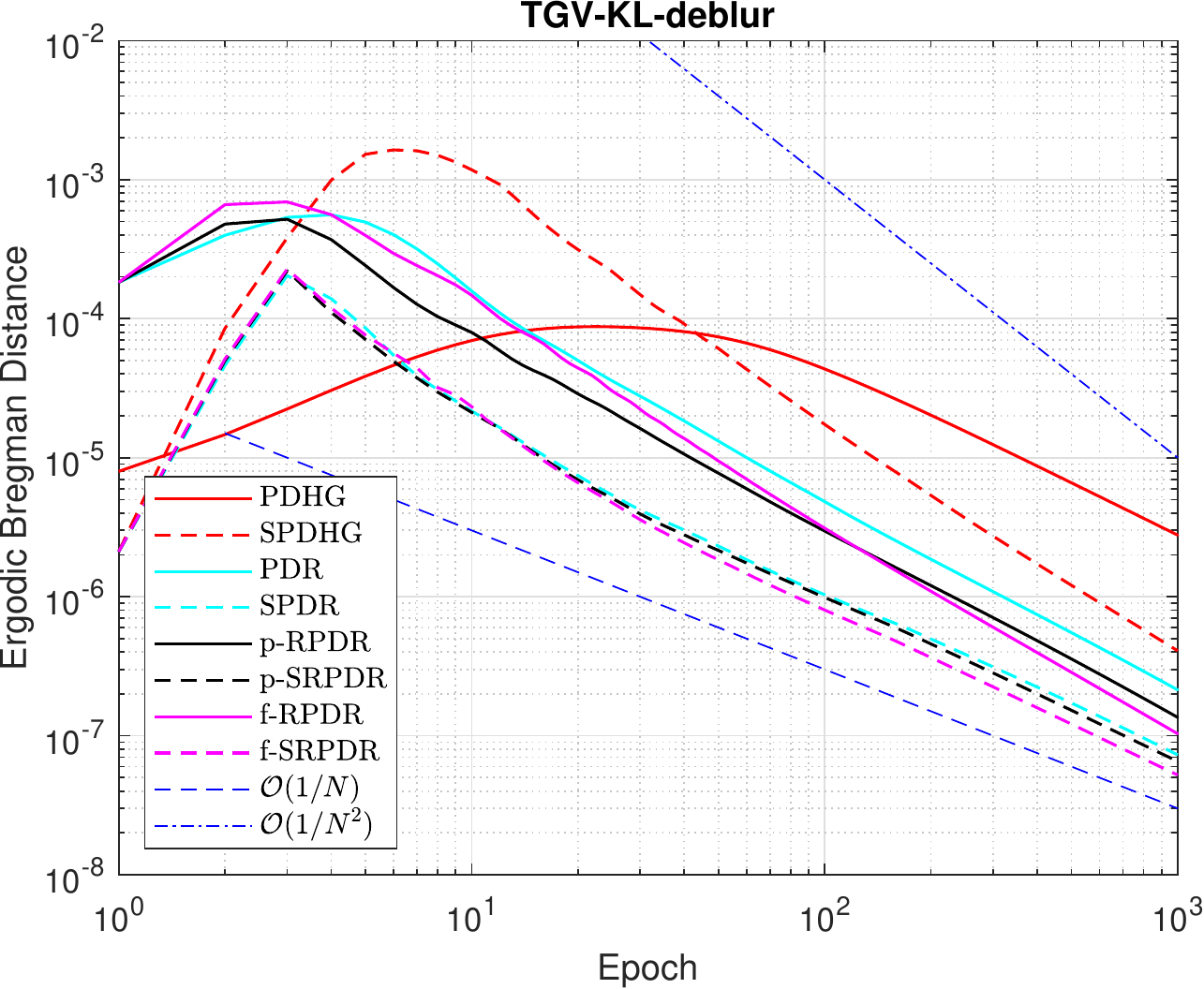}}\quad
		\subfloat[Convergence rate of normalized primal error concerning epochs.]
		{\includegraphics[width=.22\paperwidth]{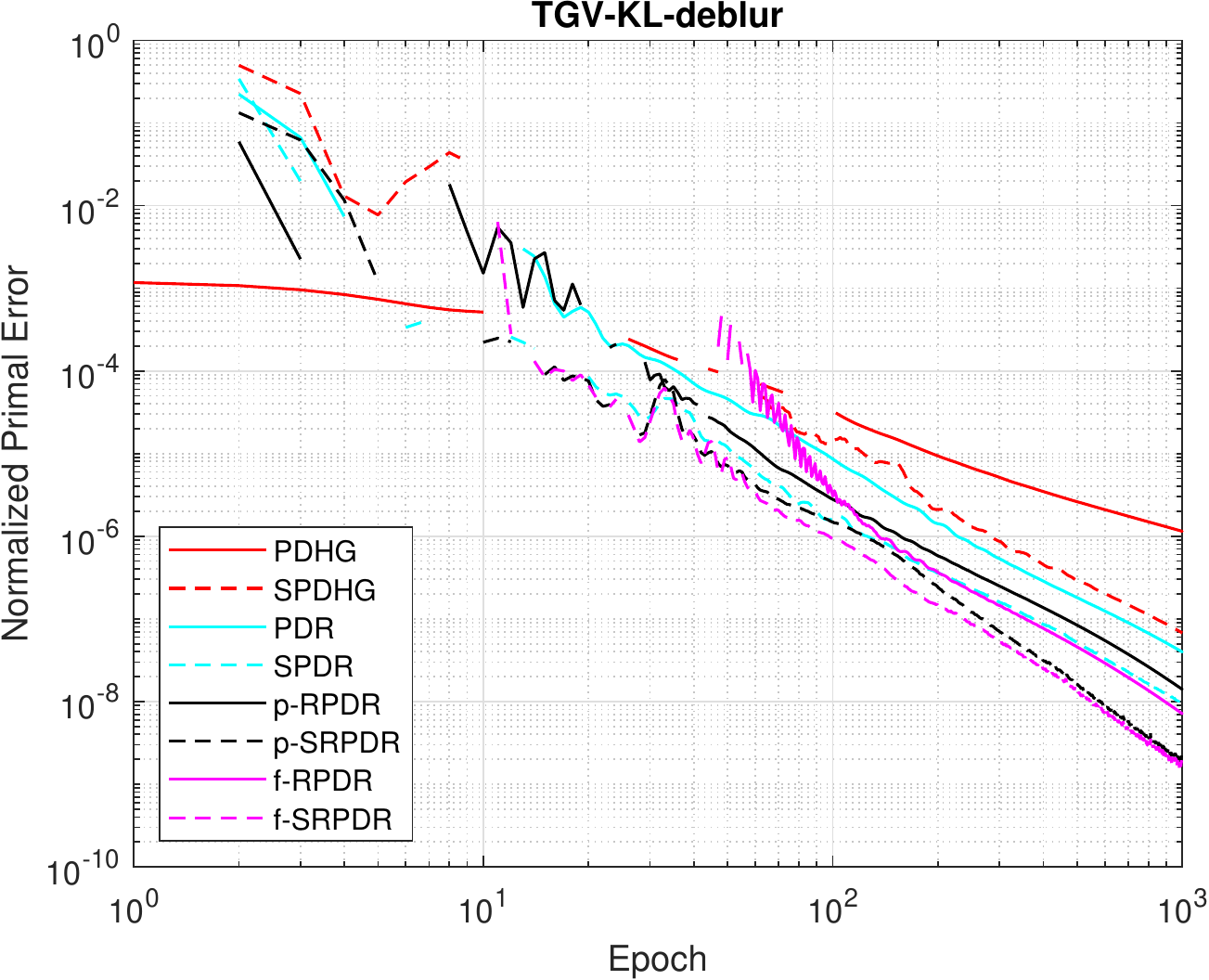}} \label{primal}\quad
		\subfloat[PSNR with respect to epochs.]
		{\includegraphics[width=.22\paperwidth]{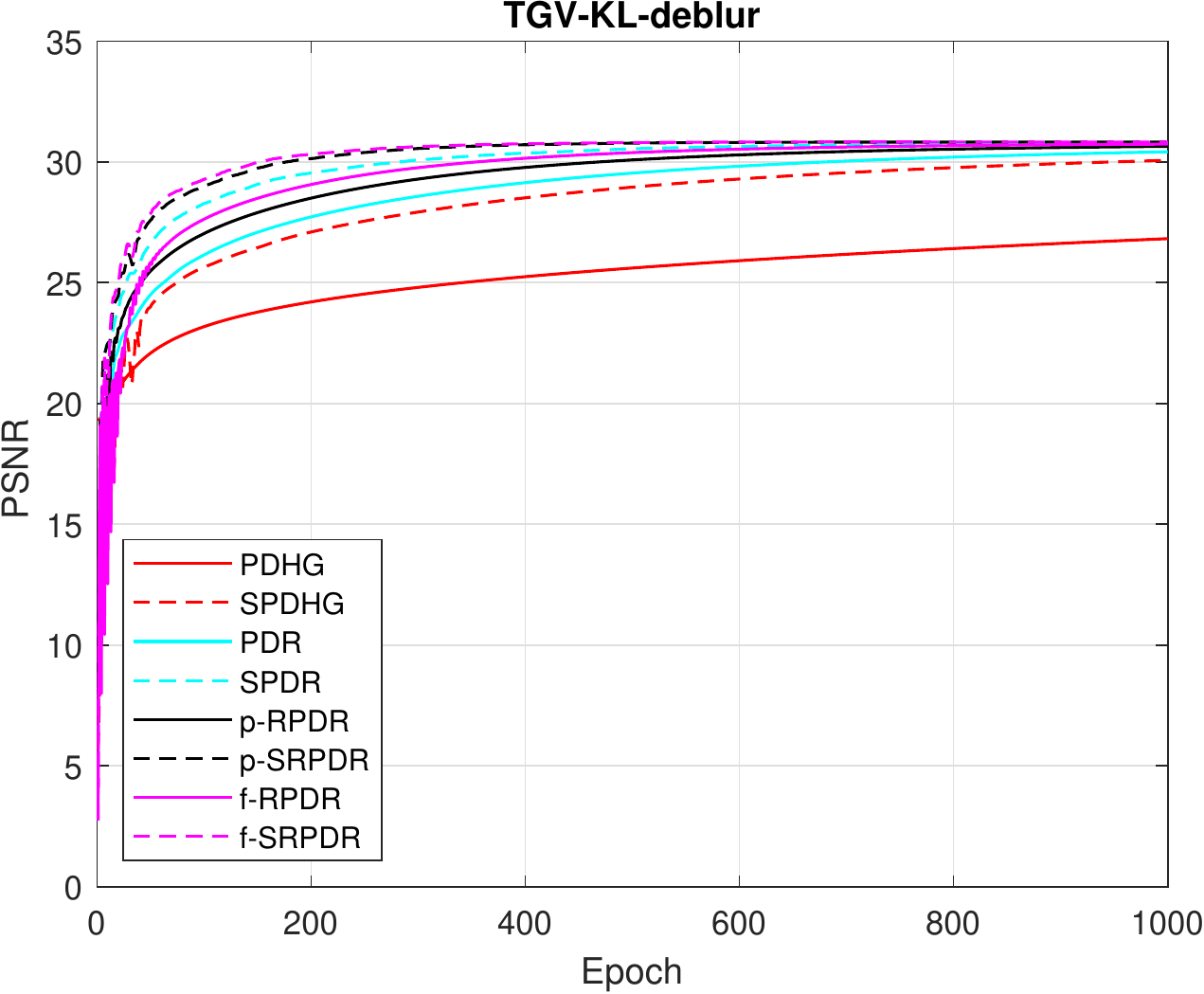}} \quad
	\end{center}
	\caption{\small \it TGV-KL deblurring problem. $\emph{(a)}$, $\emph(b)$ are normalized by $d_1\times d_2$ and in a double-logarithmic scale. Image $\emph(b)$ shows the speed comparison measured in terms of primal error, which is given by $ ( \mathcal{P}(x^{K})-\mathcal{P}(x^*) )/ \mathcal{P}(x^*) $.} \label{fig:tgvkl-deblur}
\end{figure} 

\begin{figure}[!htb] 
	\graphicspath{{fig//}}
	\begin{center}
		
		\subfloat[\tiny{Degraded image}]
		{\includegraphics[width=.126\paperwidth]{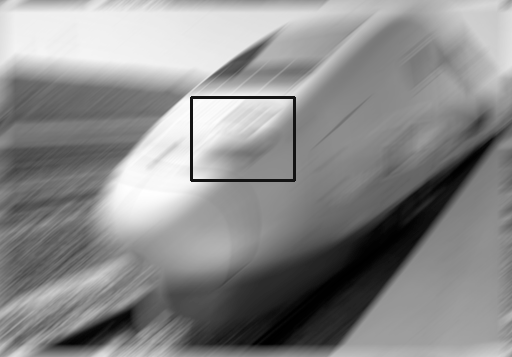}} \quad
		 		\subfloat[\tiny{PDHG, 300 epochs}]
				{\includegraphics[width=.126\paperwidth]{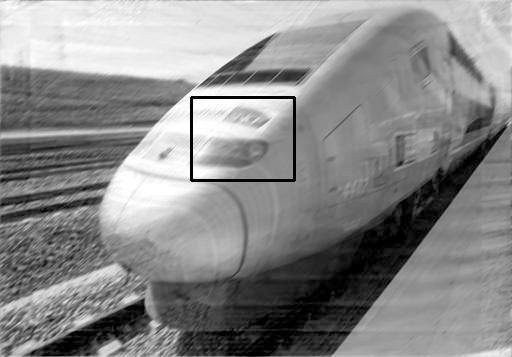}} \quad
		 		\subfloat[\tiny{PDR, 300 epochs}]
		 		{\includegraphics[width=.126\paperwidth]{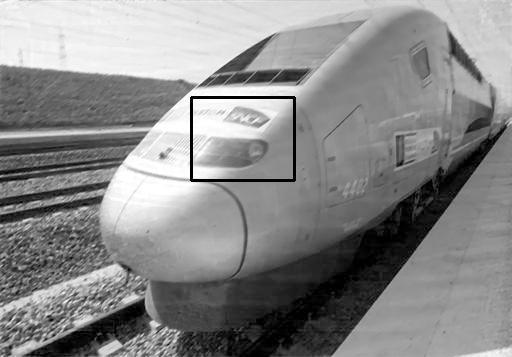}} \quad
\subfloat[\tiny{p-RPDR, 300 epochs}]
				{\includegraphics[width=.126\paperwidth]{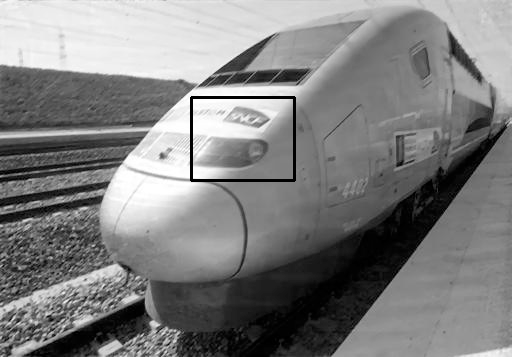}} \quad
		 		\subfloat[\tiny{f-RPDR, 300 epochs}]
		 		{\includegraphics[width=.126\paperwidth]{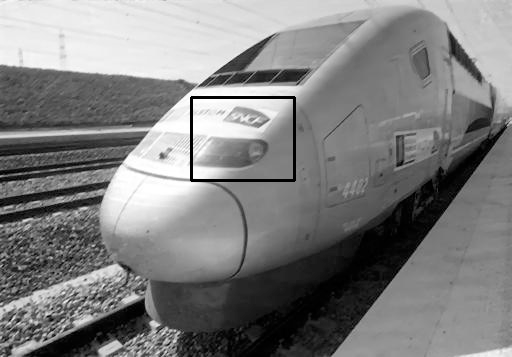}} \quad

		\subfloat[\tiny{Approximate primal part of saddle point}]
		{\includegraphics[width=.126\paperwidth]{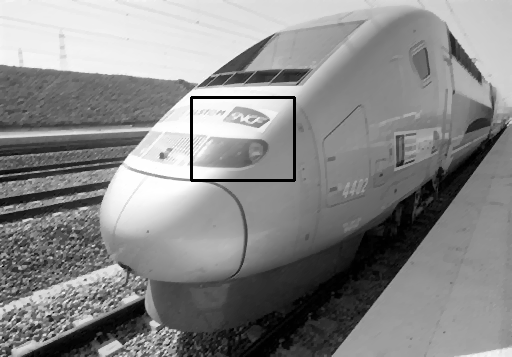}} \quad
		 		\subfloat[\tiny{SPDHG, 300 epochs}]
		 		{\includegraphics[width=.126\paperwidth]{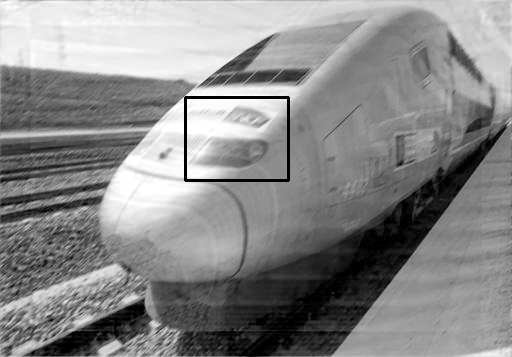}} \quad
		\subfloat[\tiny{SPDR, 300 epochs}]
		{\includegraphics[width=.126\paperwidth]{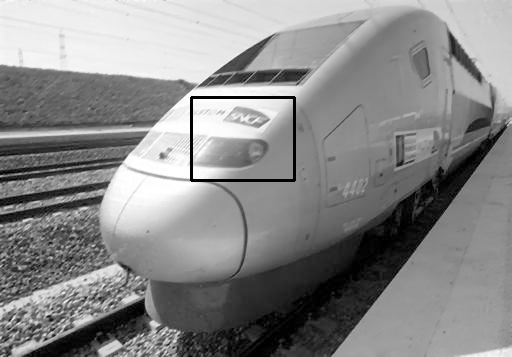}} \quad
	\subfloat[\tiny{p-SRPDR,300 epochs}]
		 		{\includegraphics[width=.126\paperwidth]{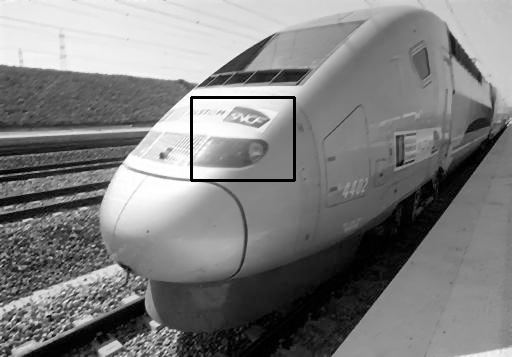}} \quad
		\subfloat[\tiny{f-SRPDR, 300 epochs}]
		{\includegraphics[width=.126\paperwidth]{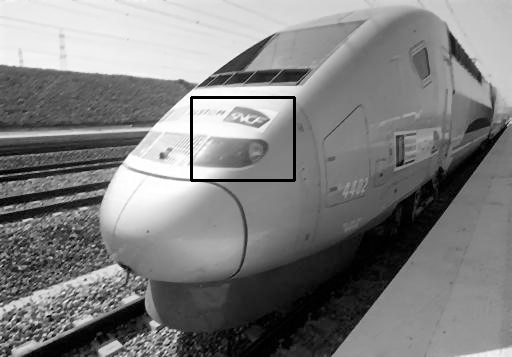}} \quad
		\subfloat[\tiny{Details of \textbf{a}}]
		{\includegraphics[width=.126\paperwidth]{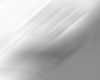}} \quad
		\subfloat[\tiny{Details of PDHG}]
		{\includegraphics[width=.126\paperwidth]{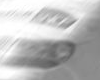}} \quad
		\subfloat[\tiny{Details of PDR}]
		{\includegraphics[width=.126\paperwidth]{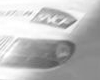}} \quad
		\subfloat[\tiny{Details of P-PDR}]
		{\includegraphics[width=.126\paperwidth]{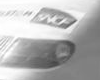}} \quad
  		\subfloat[\tiny{Details of f-RPDR}]
		{\includegraphics[width=.126\paperwidth]{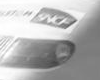}} \quad
		\subfloat[\tiny{Details of \textbf{f}}]
		{\includegraphics[width=.126\paperwidth]{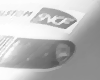}} \quad
		\subfloat[\tiny{Details of SPDHG}]
		{\includegraphics[width=.126\paperwidth]{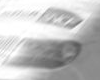}} \quad
		\subfloat[\tiny{Details of SPDR}]
		{\includegraphics[width=.126\paperwidth]{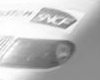}} \quad
  		\subfloat[\tiny{Details of p-SRPDR}]
		{\includegraphics[width=.126\paperwidth]{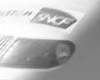}} \quad
		\subfloat[\tiny{Details of f-SRPDR}]
		{\includegraphics[width=.126\paperwidth]{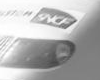}} \quad
		
	\end{center}	
	\caption{\small \it Examples of variational TGV-KL deblurring after 300 epochs with uniform sampling. }\label{fig:tgvkl-deblur-exam}	
\end{figure}
Figure \ref{fig:tgvkl-deblur} (a) shows that the Bregman distance converges at a rate of $\mathcal{O}(1/K)$ in practice. Figure \ref{fig:tgvkl-deblur} (b) shows that the primal errors are not continuous at the beginning epochs. This happens because $\mG_i$ in \eqref{eq:bredies_kl} can be infinite. Additionally, the proposed \ref{sto-pdr} and \ref{sto-rpdr} methods demonstrate high efficiency and strong performance compared to SPDHG, both in terms of Bregman distance and normalized primal error. Notably, f-SRPDR achieves a higher and more stable PSNR, as shown in Figure \ref{fig:tgvkl-deblur}.
Figure \ref{fig:tgvkl-deblur-exam} shows that the proposed SPDR, p-SRPDR, and f-SRPDR can give clearer reconstructions compared with PDR, f-RPDR, f-RPDR, SPDHG, and PDHG with the same epochs.

\subsection{Binary Classification with real data}\label{subsec:libsvm}In this subsection, we present the results of solving a binary classification problem using the real datasets \textit{Gisette} \cite{nips2003} and \textit{Madelon} \cite{nips2003}, both part of the \textit{LIBSVM} data collection \cite{libsvm}. These two datasets are selected to illustrate different relations between sample size and feature dimensionality.  Each data point is represented as a pair $(a_i, b_i)$, where $a_i \in \mathbb{R}^d$ is the feature vector and $b_i$ is the binary class label, taking values in $\{-1, 1\}$. Our objective is to solve the following regularized empirical risk minimization problem:
\begin{equation}\label{eq:min-linxiao}
    \min_{x\in \mathbb{R}^d} \frac{1}{n} \sum^n_{i=1} \phi_i(a_i^Tx) + \frac{\lambda}{2}\|x\|_2^2, \ \ \text{where} \ \phi_i(z)  =
\begin{cases}
0 & \text{if } b_i z \geq 1 \\
\frac{1}{2} - b_i z & \text{if } b_i z \leq 0 \\
\frac{1}{2} (1 - b_i z)^2 & \text{otherwise},
\end{cases}
\end{equation}
is the smoothed hinge loss, originally introduced by \cite{ShZ} and addressed with data term by \cite[Section 7.2]{Zhangxiao}. Define $A = [a_1,\dots,a_n]^T$. The corresponding $\mathcal{K}$, $\mathcal{G}_i^{*}(\mK_i x) $ and $\mathcal{F}(x) $ in the form of \eqref{eq:primal-dual-problem} read as 
$$
\mathcal{K} = \frac{1}{n}A, \quad \mG^*_i(\mK_i x) = \frac{1}{n}\phi_i(n\mK_i x) = \frac{1}{n}\phi_i(a_i^T x), \quad \mF(x) = \frac{\lambda}{2}\|x\|^2_2.
$$
The conjugate function of hinge lose $\phi_i$ is given as $\phi^*_i(z) = b_i z+ \frac{1}{2}z^2 + \mI_{\{b_iz\in[-1,0]\}}(z)$, and $\mG_i(z) = (1/n) \phi^*_i(z)$. Then the corresponding saddle point problem of \eqref{eq:min-linxiao} is written as,
\begin{equation}\label{eq:saddle-xiaolin}
\min_{x\in \mathbb{R}^d}\max_{y\in \dom{\mG}} \frac{\lambda
}{2}\|x\|^2_2+ \sum^n_{i=1}(\langle\mK_i x,y_i \rangle - \frac{1}{n}\phi_i^*(y_i)),\ \text{with}\quad \mG(y):= \sum^n_{i =1} \mG_i(y_i),\ \text{and}\ y = (y_1,\ldots,y_n)^T.
\end{equation}
Further, the resolvent of $ \mF$ is derived as $(I+\tau \partial \mF)^{-1}(z)={z}/(1+\lambda\tau)$. Set $l_i = \min(-b_i,0)$ and $r_i = \max(-b_i,0)$ thus the resolvent of $ \mG_i$ is derived as $(I+\sigma \partial \mG_i)^{-1}(y)=\text{Project}_{[l_i,r_i]}((ny-\sigma b_i)/(n+\sigma))$.  

\subsubsection{Characteristics of datasets and parameters for the algorithms} The \textit{Gisette} dataset represents the case where $n < d$, with $n = 1000$ and $d = 5000$, while the \textit{Madelon} dataset represents the case where $n > d$, with $n = 2000$ and $d = 500$. We will evaluate our algorithms on these two datasets using two regularization parameters: $\lambda = 10^{-6}$ and $\lambda = 10^{-4}$. The initial values for all algorithms are set to be zero matrices. The step sizes for each algorithm are provided below.\par

\setlength{\hangindent}{2em}
\ \ $\bullet$ 
PDHG, SPDHG: 
We use a step size of $\sigma = {2}/\|\mathcal{K}\|$ and $\tau = {1}/({2\|\mathcal{K}\|})$ for PDHG and SPDHG on the dataset \textit{Gisette} to achieve better performance. For the dataset \textit{Madelon}, we employ $\sigma = {2}/\|\mathcal{K}\|$ and $\tau = {1}/({2\|\mathcal{K}\|})$ for PDHG. In contrast, for SPDHG, we set $\sigma = \tau = {1}/(15\|\mathcal{K}\|)$, due to the significant fluctuations observed with larger step sizes on this dataset.

\setlength{\hangindent}{2em}
\ \ $\bullet$ \ref{PDR}, \ref{sto-pdr}: The step sizes are chosen as $\sigma = \tau=0.3$ for both algorithms and both datasets \textit{Madelon} and \textit{Gisette}. The precondition iteration is set to be the two-step Richardson \cite[Example 2.11, Proposition 2.15]{BSCC}. Thus the corresponding preconditioner $M$ in preconditioned iteration \eqref{preconditioner:from} is given by $(1+\sigma\tau\|\mK\|^2)I$.

\setlength{\hangindent}{2em}
\ \ $\bullet$ p-RPDR (f-RPDR), p-SRPDR (f-SRPDR): The step sizes and the preconditioners are set to be the same as those used in PDR and SPDR. The relaxation parameters are all set to 1.9. \\

\begin{figure}[!htb]
	\graphicspath{{fig//}}
	\begin{center}
		{\includegraphics[width=.33\paperwidth]{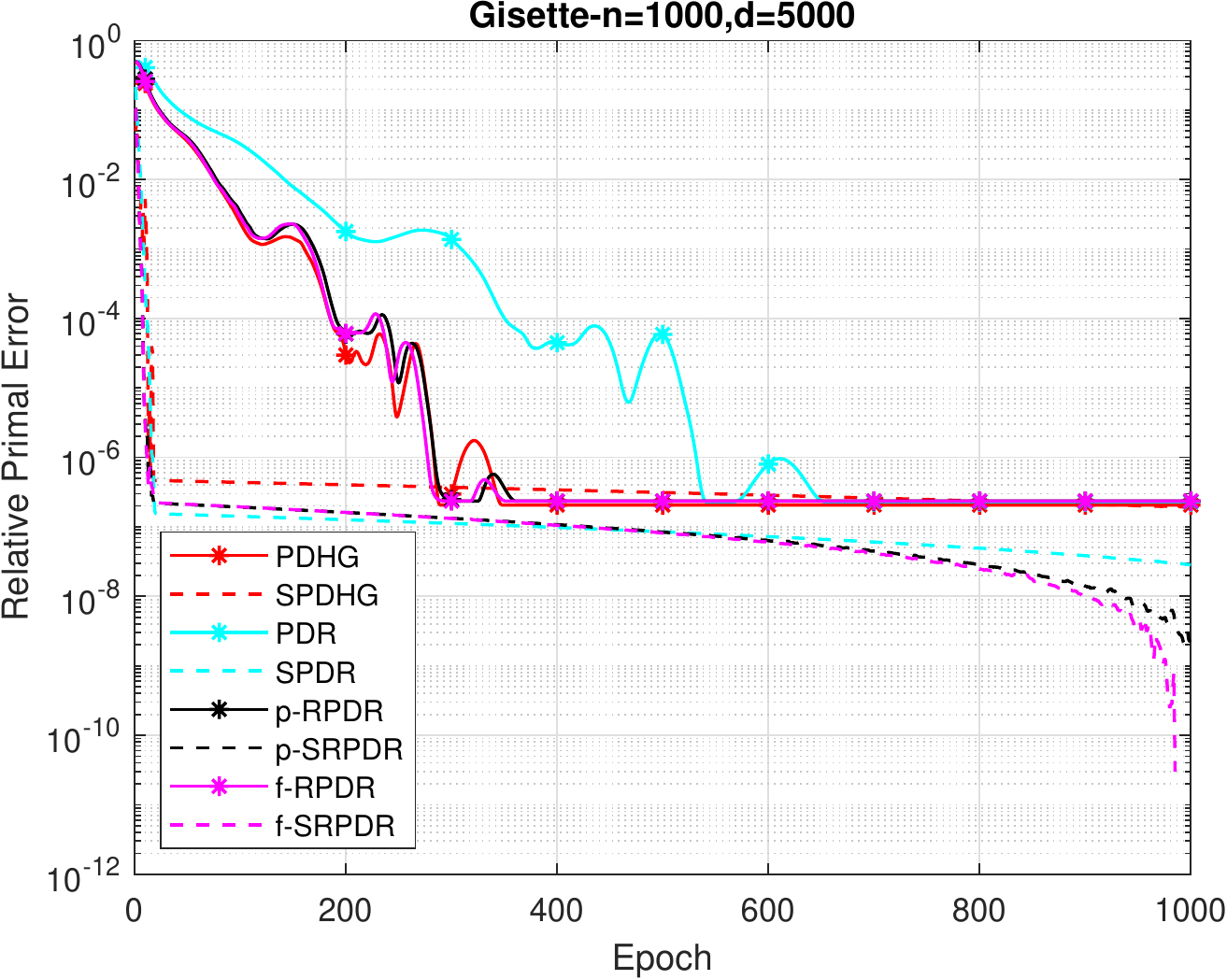}}\quad
{\includegraphics[width=.33\paperwidth]{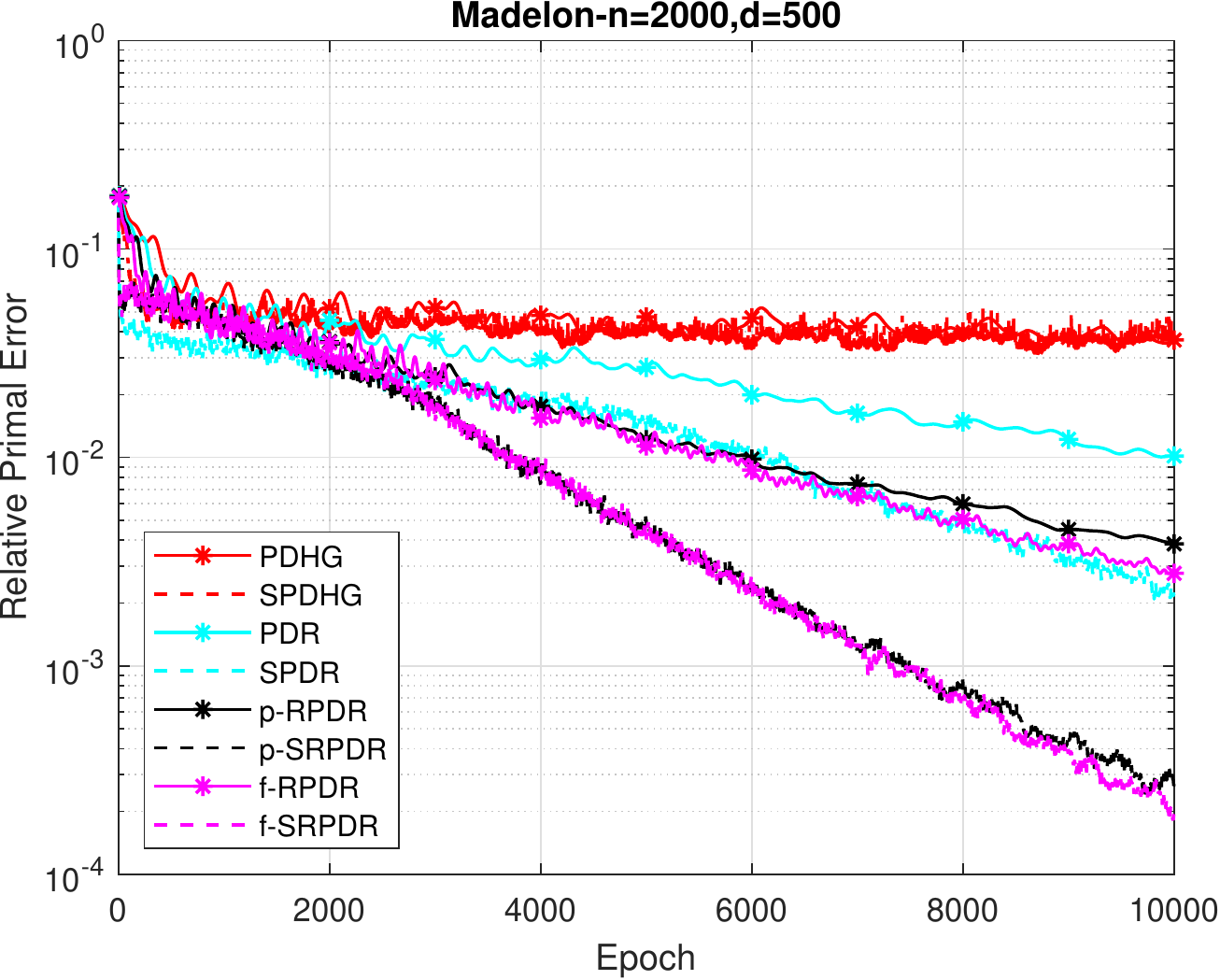}}\quad

    \end{center}
   \caption{Convergence rate of relative primal error concerning epochs for dataset \textit{Gisette} and \textit{Madelon}. The regularization parameter is set to be $1 \times 10^{-6}$.}
    \label{fig:libsvm-1e6}
\end{figure}

\begin{figure}[!htb]
	\graphicspath{{fig//}}
	\begin{center}
		{\includegraphics[width=.33\paperwidth]{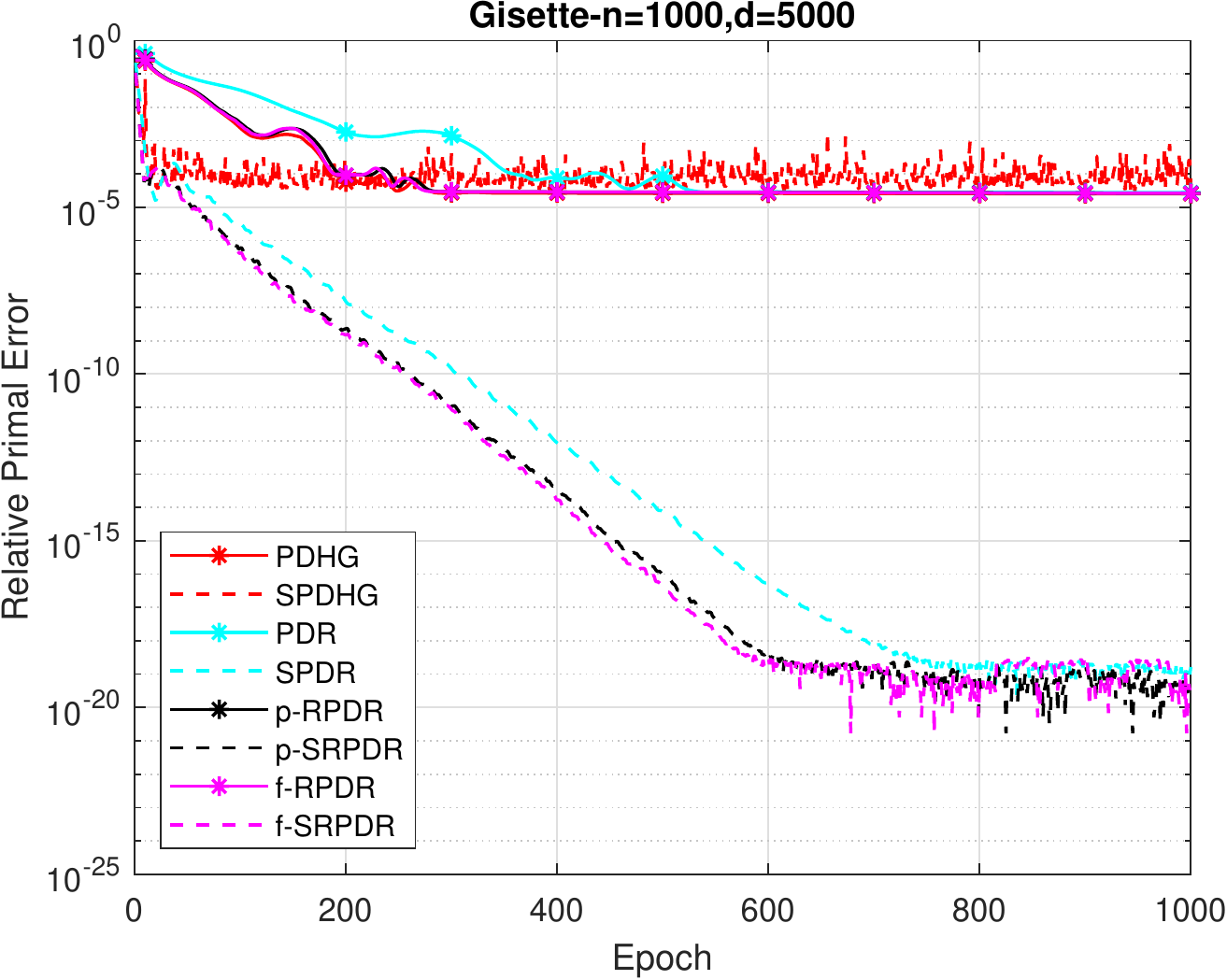}}\quad
{\includegraphics[width=.33\paperwidth]{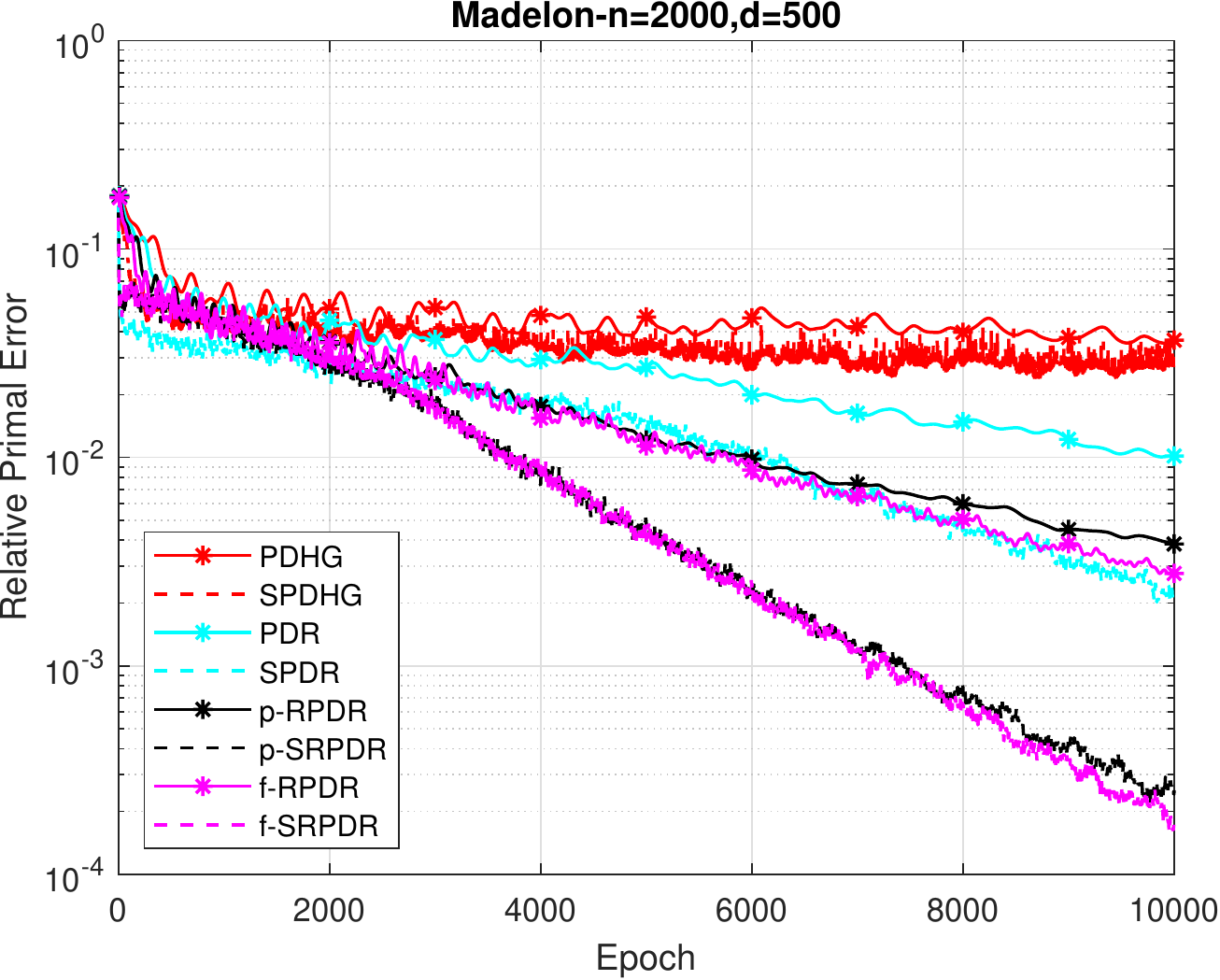}}\quad

    \end{center}
   \caption{Convergence rate of relative primal error concerning epochs for dataset \textit{Gisette} and \textit{Madelon}. The regularization parameter is set to be $1 \times 10^{-4}$. }
    \label{fig:libsvm-1e4}
\end{figure}

With these settings, we tested the primal error of eight algorithms in these two datasets. In Figure \ref{fig:libsvm-1e6} and Figure \ref{fig:libsvm-1e4}, the vertical axis is the relative primal error calculated by $\mathcal{P}(x^k) - \mathcal{P}(x^*) $ in logarithmic scale. The minimum loss $\mathcal{P}(x^*)$ for both datasets is calculated by f-RPDR after 750,000 epochs for \textit{Gisette} and after 500,000 epochs for \textit{Madelon}.


From Figure \ref{fig:libsvm-1e6} and Figure \ref{fig:libsvm-1e4}, We observe that our proposed \ref{sto-rpdr} and \ref{sto-pdr} methods are highly efficient and stable, with fewer fluctuations compared to SPDHG. We also find that the performance of these algorithms depends a lot on step sizes. For some algorithms, especially deterministic ones, the primal errors decrease quickly at first but then slow down considerably. As a result, the loss for our proposed stochastic algorithms after 1000 epochs can reach a level nearly the same as f-RPDR after 750,000 epochs on dataset \textit{Gisette}. This is why the relative primal errors of our stochastic algorithms can drop below  $10^{-16}$ on dataset \textit{Giesette}. It is worth noting that \ref{sto-pdrq} can also be used here. SRPDRQ is more compact because it does not require the updates of $\bar x ^{k+1}$ compared to SRPDR. However, since SRPDR is commonly used and the difference in primal error between SRPDR and SRPDRQ is not significant, we only present SPDR, f-SRPDR and p-SRPDR in this experiment. A comparison between SRPDR and SRPDRQ is included in the Appendix \ref{appen:com-spdrq-spdr}.

\section{Conclusion}\label{sec:conclusion}
We proposed a stochastic and relaxed PDR framework to solve general saddle-point problems for separable dual variables cases. We showed the almost sure convergence of  SRPDR along with SRPDRQ. We gave the $\mathcal{O}(1/K)$ convergence rate of the ergodic iteration sequences regarding the restricted primal-dual gap functions and the primal errors. For SRPDR, we will consider local linear convergence involving metric subregularity as in \cite{AFC} and design highly efficient preconditioners for more complicated operators in applications including the Radon transforms \cite{RGMGK} for tomography.

\section{Appendix}\label{sec:append} \begin{proof}[Proof of Proposition \ref{pro:prop-rpdr}] \label{txt:fur0-1}
    By  the equation \eqref{eq:pdr:relax:deter} and the iteration of \ref{rPDR}, we see 
    \begin{equation}\label{eq:uk:recover}
        \mT  u^k= I_{\varrho_k}^{-1}(\mT_{R} -I)u^k + u^k = I_{\varrho_k}^{-1}(u^{k+1}-u^k) + u^k . 
    \end{equation}
   Substituting $u_{t}^{k+1} = \mT  u^k$ to the equation \eqref{eq:polari} in Proposition \ref{pro:prop-pdr}, we arrive at
   \begin{align}
   \mathcal{L}(x_{\test}^{k+1}, y)-\mathcal{L}(x,y_{\test}^{k+1})&\leq\langle u^k-(I_{\varrho_k}^{-1}(u^{k+1}-u^k) + u^k),(I_{\varrho_k}^{-1}(u^{k+1}-u^k) + u^k)-u\rangle_\mathbf{M} \notag\\
   & = \langle I_{\varrho_k}^{-1}(u^{k}-u^{k+1}), I_{\varrho_k}^{-1}(u^{k+1}-u^k) + u^k-u  \rangle_\mathbf{M} \notag \\
   &= \langle I_{\varrho_k}^{-1}(u^{k}-u^{k+1}), (I_{\varrho_k}^{-1}-I)(u^{k+1}-u^k) + u^{k+1}-u  \rangle_\mathbf{M} \notag \\
   & = \langle I_{\varrho_k}^{-1}(u^{k}-u^{k+1}),  u^{k+1}-u  \rangle_\mathbf{M} 
   - \langle I_{\varrho_k}^{-1}(u^{k+1}-u^{k}), (I_{\varrho_k}^{-1}-I)(u^{k+1}-u^k)\rangle_\mathbf{M} \notag \\
  & = \langle u^{k}-u^{k+1},  u^{k+1}-u  \rangle_{\mathbf{M} I_{\varrho_k}^{-1}} - \|u^{k+1}-u^k\|_{\mathbf{M} I_{\varrho_k}^{-1} (I_{\varrho_k}^{-1}-I)}^2.  \label{eq:polar:1} 
   \end{align}
   By the polarization identity and the positive semidefiniteness of $\mathbf{M} I_{\varrho_k}^{-1}$, we obtain
   \begin{equation}\label{eq:polari:relax:separate}    
   \langle u^{k}-u^{k+1},  u^{k+1}-u  \rangle_{\mathbf{M} I_{\varrho_k}^{-1}} = \frac{1}{2}\big( \|u^k-u\|_{\mathbf{M} I_{\varrho_k}^{-1}}^2 - \|u^{k+1}-u\|_{\mathbf{M} I_{\varrho_k}^{-1}}^2  - \|u^{k+1}-u^k\|_{\mathbf{M} I_{\varrho_k}^{-1}}^2 \big).
   \end{equation}

   Substituting the above equality to \eqref{eq:polar:1} and noting $\mathbf{M} I_{\varrho_k}^{-1}+ 2\mathbf{M} I_{\varrho_k}^{-1} (I_{\varrho_k}^{-1}-I)=\mathbf{M}  (2I-I_{\varrho_k})I^{-2}_{\rho_k}$, we get \eqref{eq:polari-rpdr}.
\end{proof}

We need the following lemma first to prove Lemma \ref{lemma:fix-rpdr}.
\begin{lemma}\label{lem:weak:convergence:sequ}\label{txt:fur0-2}
Suppose $\rho_{k} \in (0,2)$ be a nondecreasing sequence, and  $\rho_{k} \rightarrow \bar \rho$ with $\bar \rho \in (0,2)$. If $\mathbf{c}^{k} \rightharpoonup \mathbf{c}$, $\bar \rho  \neq 1$, and the sequence $\{\mathbf{a}^k\}$ satisfy
\begin{equation}\label{eq:sesaro}
\mathbf{a}^{k+1} - (1-\rho_{k})\mathbf{a}^{k} = \rho_{k} \mathbf{c}^{k},
\end{equation}
 then we have $\mathbf{a}^{k} \rightharpoonup \mathbf{c}$.
For the special case $\rho_k \equiv 1$ with $\bar \rho =1$, we also have $\mathbf{a}^{k}\rightharpoonup \mathbf{c} $.
\end{lemma}
\begin{proof}
It can be seen that
\[
\rho_k \mathbf{c}^k = (\rho_{k} - \bar \rho)\mathbf{c}^{k} + \bar \rho \mathbf{c}^{k}.
\]
Since $\mathbf{c}^{k}$ is bounded and $\rho_{k}\rightarrow \bar \rho$, we have
$\rho_k \mathbf{c}^k \rightharpoonup \bar \rho \mathbf{c}$ with the assumption $\mathbf{c}^{k} \rightharpoonup \mathbf{c}$.
While $\rho_k \equiv 1$ with $\bar \rho =1$, we conclude that $\mathbf{a}^{k+1} =\mathbf{c}^{k}\rightharpoonup   \mathbf{c}$.

If $\bar \rho \neq 1$, then there exists $k_0$, such that either $\rho_k \in (1,2)$ while $\bar \rho \in (1,2)$ or $\rho_k \in (0, \bar \rho )$ while $\bar \rho \in (0,1)$ for all $k \geq k_0$.  Henceforth, denote $\alpha_{k+1}:=(1-\rho_{k})$, $\mathbf{b}^{k+1} := \mathbf{a}^{k+1} - \alpha_{k+1} \mathbf{a}^{k}$ for $k\geq k_0$ and define $\mathbf{b}^{k_0} := \mathbf{a}^{k_0}$ along with $\alpha_{k_0}=1$ for convenience.
By induction with direct calculation from $k_0$, we obtain
\[
\mathbf{a}^{k_0+k} = \bigl(\sum_{k'=0}^k \mathbf{b}^{k_0+k'} \prod_{k''=0}^{k'} \alpha_{k_0+k''}^{-1}\bigr)/\prod_{k'=0}^{k}
\alpha_{k_0+k'}^{-1}, \ \forall k \geq 0.
\]
 It can be seen that $\mathbf{b}^{k+1} \rightharpoonup \mathbf{b}:=\bar \rho \mathbf{c}$ due to $\mathbf{b}^{k+1}=\rho_{k} \mathbf{c}^{k} $ for $k> k_0$ with \eqref{eq:sesaro} and $\alpha_k \rightarrow \alpha:=(1-\bar \rho)$.
For the case $\bar \rho <1$, i.e.,  $\alpha_k <1$ for $k > k_0$,  for any $f$, applying Stolz-Ces\`{a}ro theorem, we have
\[
\lim_{k \rightarrow \infty} \langle \mathbf{a}^{k_0+k}, f \rangle = \lim_{k \to \infty}
  \frac{\langle \mathbf{b}^{k_0+k}, f\rangle  \prod_{k''=0}^{k} \alpha_{k_0+k''}^{-1}}{\prod_{k''=0}^{k} \alpha_{k_0+k''}^{-1} - \prod_{k''=0}^{k-1} \alpha_{k_0+k''}^{-1}} = \lim_{k \to \infty} \frac{\langle \mathbf{b}^{k_0+k}, f\rangle }{1 - \alpha_{k_0+k}}
  = \frac{\langle \mathbf{b}, f\rangle }{1 - \alpha} = \langle \mathbf{c}, f\rangle.
\]
We thus conclude that $\mathbf{a}^k \rightharpoonup \mathbf{c}$. For the case $\bar \rho \in (1,2)$, since  $\alpha_k < 0$ for $k \geq k_0$, the proof is similar and is omitted here. 
\end{proof}

\begin{proof}[Proof of Lemma \ref{lemma:fix-rpdr}] \label{txt:fur0-3}
Regarding the first point, notice the fixed points of \eqref{eq:pdr:relax:deter} (i.e., \ref{rPDR}) satisfy
 \[
 u^* = (I - I_{\varrho_k}) u^* + I_{\varrho_k} \mathcal{T}u^*,   
 \]
which is equivalent to $\mathcal{T}u^* = u^*$ by the invertibility of $I_{\varrho_k}$. Thus the first point can be obtained with \cite[Lemma 2.7]{BSCC} (or \cite[Lemma 3]{BS1}).

For the second point of this lemma, writing \eqref{eq:pdr:relax:deter} component-wisely and denoting $\bar{I}_{k_i} :=  \text{Diag}[\rho_{k_i,x}, \rho_{k_i,y}]$, we have
\begin{align}
\begin{pmatrix}
x^{k_i+1} \\ y^{k_i+1}
\end{pmatrix}
 &= (1-\rho_{k_i})\begin{pmatrix} x^{k_i} \\ y^{k_i}\end{pmatrix} + \rho_{k_i}
 \begin{pmatrix}
N_1 & \sigma K^* \\
-\tau K & N_2
 \end{pmatrix}^{-1}
 \begin{pmatrix}
 (N_1-I)x^{k_i} +\bar x^{k_i} \\
 (N_2-I)y^{k_i} +\bar y^{k_i}
 \end{pmatrix}, \label{eq:compo:xy:relax1} \\
\begin{pmatrix}
\bar x^{k_i+1} \\ \bar y^{k_i+1}
\end{pmatrix}
&= \begin{pmatrix}
\bar x^{k_i} \\ \bar y^{k_i}
\end{pmatrix}
+\bar{I}_{k_i}
\begin{pmatrix}
(I + \sigma \partial \mF)^{-1}(2x_{t}^{k_i+1} -\bar x^{k_i}) - x_{t}^{k_i+1} \\
(I + \tau \partial \mG)^{-1}(2y_{t}^{k_i+1} - \bar y^{k_i}) - y_{t}^{k_i+1}
\end{pmatrix}. \label{eq:compo:bar:xy:relax1}
\end{align}
Denote 
 $ A_{N}:= \begin{pmatrix}
N_1 & \sigma K^* \\
-\tau K & N_2
 \end{pmatrix}$ and $ A: =  \begin{pmatrix}
I & \sigma K^* \\
-\tau K & I
 \end{pmatrix}$. Multiplying $A_N$ to both sides of \eqref{eq:compo:xy:relax1}, by rearranging and cancelling some terms, we arrive at
 \begin{equation}
 \begin{pmatrix}
 N_1(x^{k_i+1} - x^{k_i}) + \sigma K^*(y^{k_i+1}-y^{k_i}) \\
 -\tau K(x^{k_i+1}-x^{k_i}) + N_2(y^{k_i+1}-y^{k_i})
 \end{pmatrix}
 =  \rho_{k_i} \begin{pmatrix}
 \bar x^{k_i} - x^{k_i} -\sigma K^*y^{k_i} \\
 \bar y^{k_i} -y^{k_i} + \tau K x^{k_i}
 \end{pmatrix}.
 \end{equation}
 Adding $(x^{k_i}-x^{k_i+1}, y^{k_i}-y^{k_i+1})^{T}$ to both sides of above equation, we have
  \begin{equation}\label{eq:relax:rela:deter}
 \begin{pmatrix}
I & \sigma K^* \\
-\tau K & I
 \end{pmatrix}
 \begin{pmatrix}
 x^{k_i+1} \\ y^{k_i+1}
 \end{pmatrix} + h^{k_i}
 =  \rho_{k_i} \begin{pmatrix}
 \bar x^{k_i}  \\
 \bar y^{k_i}
 \end{pmatrix}+ (1-\rho_{k_i})
  \begin{pmatrix}
I & \sigma K^* \\
-\tau K & I
 \end{pmatrix}
 \begin{pmatrix}
x^{k_i} \\
y^{k_i}
 \end{pmatrix},
 \end{equation}
where $h^{k_i}  = ((N_1-I)(x^{k_i+1}-x^{k_i}), (N_2-I)(y^{k_i+1}-y^{k_i}))^{T} \rightarrow 0$ by the assumption. We conclude that $A$ has a bounded inverse
  since $I + \sigma \tau K^*K$ is positive definite. We thus can rewrite \eqref{eq:relax:rela:deter} as follows
\begin{equation}\label{eq:weak:convergent}
 \begin{pmatrix}
 x^{k_i+1} \\ y^{k_i+1}
 \end{pmatrix} -
 (1-\rho_{k_i})
 \begin{pmatrix}
x^{k_i} \\
y^{k_i}
 \end{pmatrix}
 = \rho_{k}A^{-1}
 \begin{pmatrix}
 \bar x^{k_i}  \\
 \bar y^{k_i}
 \end{pmatrix} - A^{-1}h^{k_i}.
\end{equation}
By the assumption of this lemma, $\{\bar y^{k}\}_k$
and $\{\bar x^{k}\}_k$ are bounded. There thus exist weakly convergent subsequences and denote them as $\{ \bar x^{k}_i\}_k$ and $\{ \bar y^{k}_i\}_k$ with corresponding weak limit $\bar x$ and $\bar y$. Since $A^{-1}$ is linear and bounded, then $\rho_kA^{-1}(\bar x^{k}_i, \bar y^{k}_i)^{T} - A^{-1}h^{k}$ is weakly convergent with weak limit $\rho^* A^{-1}(\bar x, \bar y)^{T}$. Furthermore, applying Lemma \ref{lem:weak:convergence:sequ} to each component of \eqref{eq:weak:convergent}, we arrive at
\begin{equation}\label{eq:x:bar:limit:relation}
\begin{pmatrix}
x^{k}_i \\
y^{k}_i
\end{pmatrix}
\rightharpoonup
\begin{pmatrix}
 x^* \\
 y^*
\end{pmatrix}
:=A^{-1}
\begin{pmatrix}
\bar x^* \\
\bar y^*
\end{pmatrix}.
\end{equation}

We now turn to prove  $(x^*, y^*)$ is a saddle point for~\eqref{eq:saddle-point}. Actually we just need to prove that $(x_{\test}^{k_i+1}, y_{\test}^{k_i+1})$ weakly converge to $(x^*, y^*)$. With \eqref{eq:compo:bar:xy:relax1}, we can obtain
\begin{align}
&(\bar x^{k_i+1} -\bar x^{k_i})/\rho_{k_i,x}  + x_t^{k_i+1} = x_{\test}^{k_i+1}, \label{eq:x:test:relation}\\
&(\bar y^{k_i+1} -\bar y^{k_i})/\rho_{k_i,y}  + y_t^{k_i+1} = y_{\test}^{k_i+1},
\end{align}
and $(x_t^{k_i+1}, y_t^{k_i+1})$ have the same weak limits as $(x_{\test}^{k_i+1}, y_{\test}^{k_i+1})$. Together with \eqref{eq:compo:xy:relax1}, we see $(x^{k_i}, y^{k_i})$ have the same weak limits with $(x_t^{k_i+1}, y_t^{k_i+1})$. We thus conclude that
\begin{equation}\label{eq:test:weak:limit:dr1}
x_{\test}^{k_i} \rightharpoonup x^*, \quad y_{\test}^{k_i} \rightharpoonup y^*.
\end{equation}
Still with \eqref{eq:compo:xy:relax1} and the definitions of $x_{\test}^{k_i+1}$ and $ y_{\test}^{k_i+1} $, we get
\begin{equation} \label{eq:inlcusion:dr1:test}
\begin{aligned}
2x_{t}^{k_i+1}-\bar x^{k_i} - x_{\test}^{k_i+1} &\in \sigma \partial
\mF(x_{\test}^{k_i+1}) \\
2y_{t}^{k_i+1}-\bar y^{k_i} - y_{\test}^{k_i+1} &\in \tau \partial
\mG(y_{\test}^{k_i+1}).
\end{aligned}
\end{equation}
Combining \eqref{eq:x:test:relation} and \eqref{eq:inlcusion:dr1:test}, we have
\begin{align*}
&\sigma \partial \mF(x_{\test}^{k_i+1}) + \sigma K^* y_{\test}^{k_i+1} \ni 2x_{t}^{k_i+1}-\bar x^{k_i} - x_{\test}^{k_i+1} + \sigma K^* y_{\test}^{k_i+1} \\
&=x_{t}^{k_i+1} -\bar x^{k_i} + \sigma K^*y_t^{k_i+1} - (\bar x^{k_i+1} -\bar x^{k_i})/\rho_{k_i,x} + \sigma K^*(\bar y^{k_i+1} -\bar y^{k_i})/\rho_{k_i,y}.
\end{align*}
Still with \eqref{eq:compo:xy:relax1}, we obtain
\begin{equation}\label{eq:x:t:k:relation}
x^{k_i+1}-x^{k_i} = \rho_k(x_t^{k_i+1} - x^{k_i}), \quad y^{k_i+1}-y^{k_i} = \rho_k(y_t^{k_i+1} - y^{k_i}).
\end{equation}
Remembering the updates of $x_t^{k_i+1}$ and $y_t^{k_i+1}$ with $u^{t+1} = \mT u^k$, we have
\begin{equation}
x_{t}^{k_i+1} -\bar x^{k_i} + \sigma K^*y_t^{k_i+1} = (N_1-I)(x^{k_i} - x_t^{k_i+1}) = (N_1-I)(x^{k_i} - x^{k_i+1})/\rho_{k_i}.
\end{equation}
Together with \eqref{eq:x:test:relation} and \eqref{eq:x:t:k:relation}, we see
\[
\sigma \partial \mF(x_{\test}^{k_i+1}) + \sigma K^* y_{\test}^{k_i+1} \rightarrow 0.
\]
Similarly, we can get 
$\tau \partial \mG(y_{\test}^{k_i+1}) -\tau Kx_{\test}^{k_i+1} \rightarrow 0$.
 Since $(x,y) \mapsto \bigl(\mK^* y + \partial \mF(x), -\mK x
  + \partial \mG(y) \bigr)$ is maximally
  monotone, it is  particular weak-strong closed. Together with \eqref{eq:test:weak:limit:dr1}, we finally arrive at
  \[
  0 \in \mK^* y^* + \partial \mF(x^*), \quad 0 \in -\mK x^*
  + \partial \mG(y^*).
  \]
\end{proof}
    \begin{proof}[Proof of Theorem \ref{thm:convergence-rpdr}]\label{txt:fur0-4}
 With the convexity of  $\mF$ and $\mG$, we conclude that $\mL(x,y)$ is convex in $x$ and concave in $y$. Together with the definitions of $x_{\test}^K$ and $y_{\test}^K$, noting that $\mathbf{M} I_{\varrho_k}^{-1}  \succeq \mathbf{M} I_{\varrho_{k+1}}^{-1}$ for $k=0,\ldots, K$ and $K_1: = K+1$,
    it can be seen that 
\begin{align*}
&\mH(z_{\test,K},z)=\mL(x_{\test}^{K},y) - \mL(x, y_{\test}^{K}) \leq \frac{1}{K_1}\sum_{k=0}^{K} (\mL(x_{\test}^{k}, y) - \mL(x, y_{\test}^{k})) \notag  \\
& \leq \frac{1}{2K_1}\sum_{k=0}^{K} \left( {\|u^k-u\|_{\mathbf{M} I_{\varrho_k}^{-1}}^2 - \|u^{k+1}-u\|_{\mathbf{M} I_{\varrho_k}^{-1}}^2  - \|u^{k+1}-u^k\|_{\mathbf{M} I_{\varrho_k}^{-1}}^2 } \right) \label{eq:ergodic:convex}\\
& \leq \frac{1}{2K_1} \left( {\|u^0-u\|_{\mathbf{M} I_{\varrho_0}^{-1}}^2 - \|u^{K+1}-u\|_{\mathbf{M} I_{\varrho_{K}}^{-1}}^2  -\sum_{k'=0}^{K} \|u^{k+1}-u^k\|_{\mathbf{M} I_{\varrho_k}^{-1}}^2 } \right) \leq  \frac{1}{2K_1 \rho_l} \|u^0-u\|_{\mathbf{M}}.  \notag 
\end{align*}
Taking supremum regarding $z \in \mathbb{B}$ in the above inequality leads to this theorem.
    \end{proof}

\begin{remark} \label{appen:com-spdrq-spdr}
Our experiments shows that \ref{sto-rpdr} and \ref{sto-pdrq} have no significant differences in primal errors with respect to the number of epochs. However SRPDRQ has lower complexity compared to SRPDR; further details can be found in Section \ref{sec:sto:pdrq}.

To illustrate, we consider a single experiment on the binary classification problem \eqref{eq:min-linxiao} using the dataset \textit{Madelon}. The regularization parameter is set to be $1 \times 10^{-6}$. We compare the relative primal error between SPDRQ (f-SRPDRQ) and SPDR (f-SRPDR), furthermore, we compute the absolute difference in primal error between SPDRQ (f-SRPDRQ) and SPDR (f-SRPDR). The relaxation parameters for stochastic relaxed algorithms are all set to 1.9. The results are shown in Figure \ref{fig:spdrq}. We can see that the absolute difference between the given methods is of level $10^{-10}$ to $10^{-5}$, which does not have a significant difference.

\begin{figure}[!htb] 
       \centering
\subfloat[\scriptsize{Relative primal error of SPDR and SPDRQ}]{\includegraphics[width=.23\paperwidth]{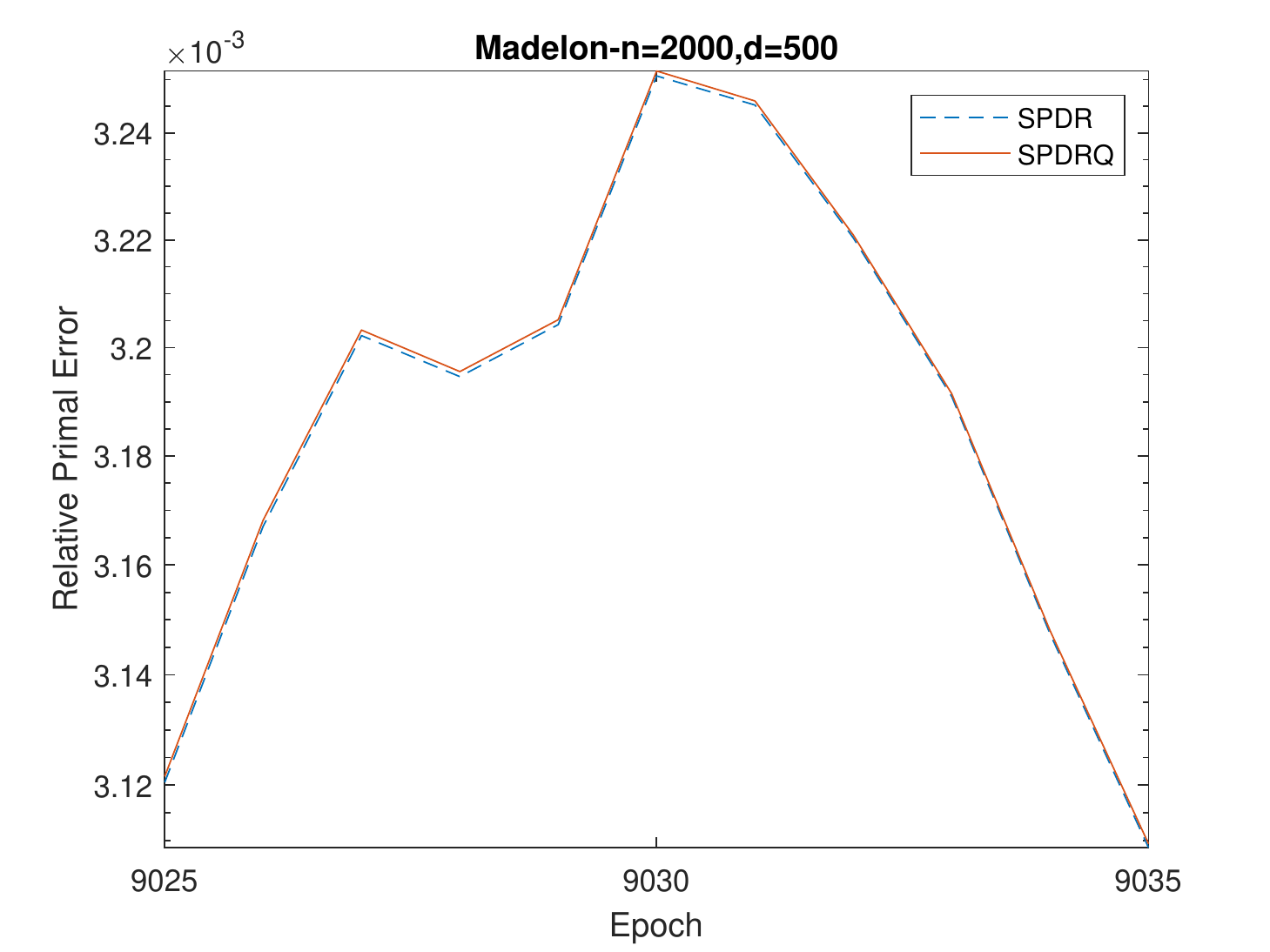}}\
\subfloat[\scriptsize{Relative primal error of f-SRPDR and f-SRPDRQ}]{\includegraphics[width=.23\paperwidth]{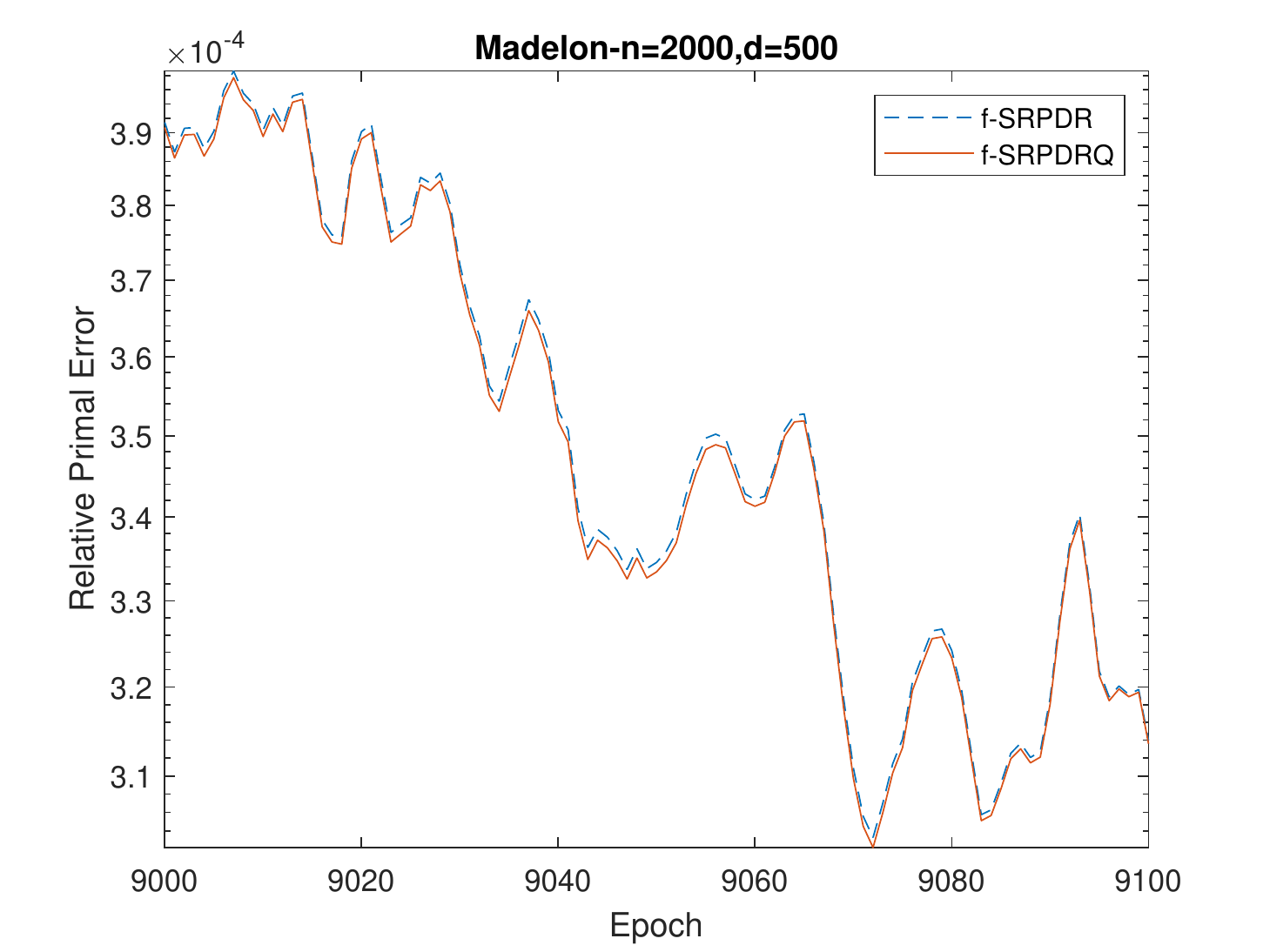}}\ \subfloat[\scriptsize{Absolut difference on primal error between SPDRQ (f-SRPDRQ) and SPDR (f-SRPDR)}]
{\includegraphics[width=.23\paperwidth]{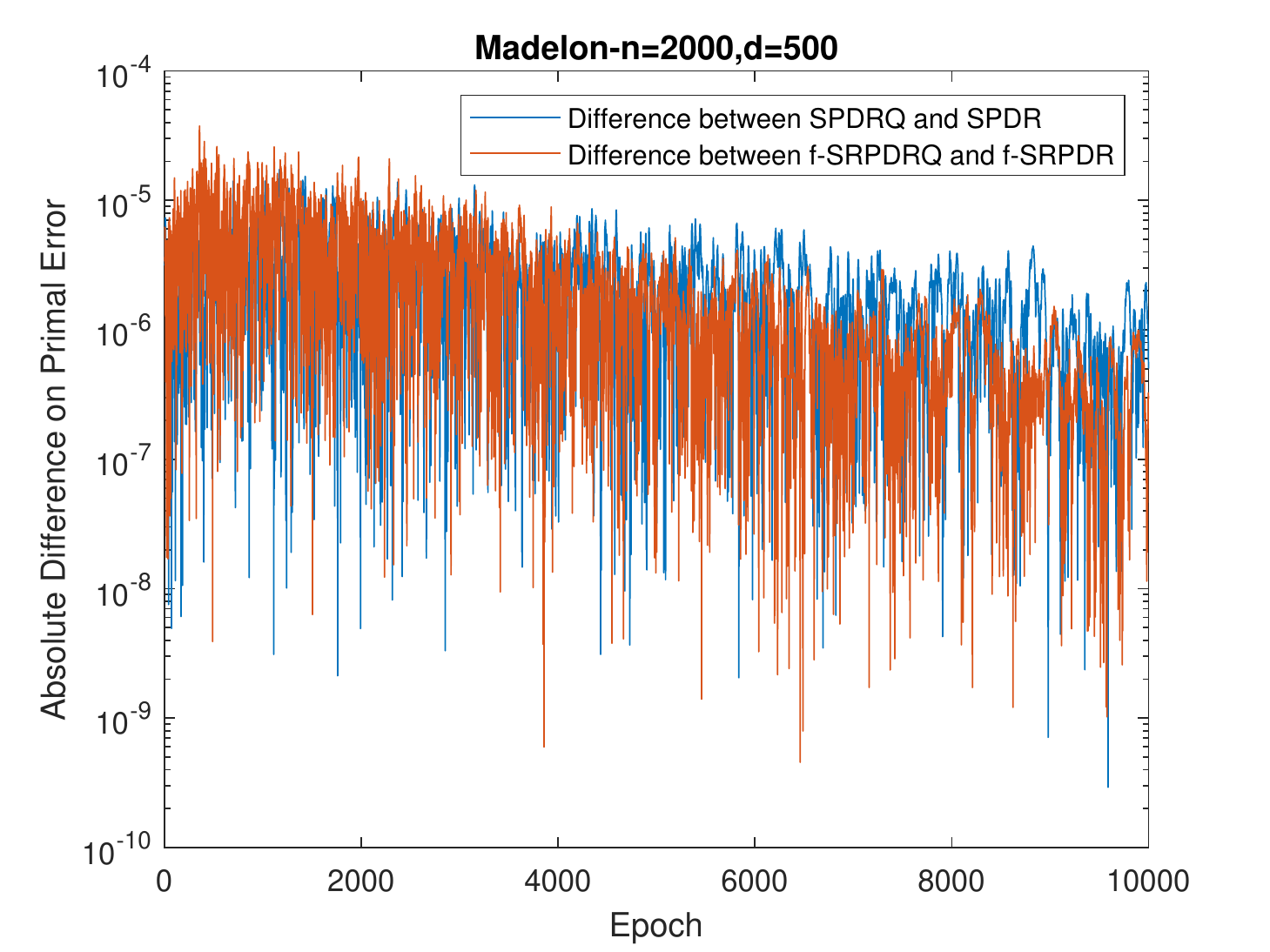}}\caption{Comparison between SRPDR and SPRDRQ on primal errors with respect to the epochs.}
\label{fig:spdrq}
\end{figure}
\end{remark}

\bibliographystyle{plain}
\bibliography{sto_pdr}

\begin{thebibliography}{10}

\bibitem{AFC}
Ahmet Alacaoglu, Olivier Fercoq, and Volkan Cevher.
\newblock On the convergence of stochastic primal-dual hybrid gradient.
\newblock {\em SIAM Journal on Optimization}, 32(2):1288--1318, 2022.

\bibitem{Artacho}
Francisco~J. Arag{\'o}n~Artacho, Jonathan~M. Borwein, and Matthew~K. Tam.
\newblock Global behavior of the douglas--rachford method for a nonconvex feasibility problem.
\newblock {\em Journal of Global Optimization}, 65(2):309--327, Jun 2016.

\bibitem{AC}
Hedy Attouch and Alexandre Cabot.
\newblock Convergence rates of inertial forward-backward algorithms.
\newblock {\em SIAM Journal on Optimization}, 28(1):849--874, 2018.

\bibitem{HBPL2}
H.~H. Bauschke and P.~L. Combettes.
\newblock {\em Convex Analysis and Monotone Operator Theory in Hilbert Spaces}.
\newblock Springer, Cham, 2017.

\bibitem{DPB}
Dimitri~P. Bertsekas.
\newblock Incremental proximal methods for large scale convex optimization.
\newblock {\em Mathematical Programming}, 129(2):163--195, Oct 2011.

\bibitem{BPK}
K.~Bredies, K.~Kunisch, and T.~Pock.
\newblock Total generalized variation.
\newblock {\em SIAM Journal on Imaging Sciences}, 3(3):492--526, 2010.

\bibitem{BS1}
K.~Bredies and H.~Sun.
\newblock Accelerated douglas-rachford methods for the solution of convex-concave saddle-point problems.
\newblock {\em arXiv:1604.06282v1}, 2016.

\bibitem{BrediesELN}
Kristian Bredies, Enis Chenchene, Dirk~A. Lorenz, and Emanuele Naldi.
\newblock Degenerate preconditioned proximal point algorithms.
\newblock {\em SIAM Journal on Optimization}, 32(3):2376--2401, 2022.

\bibitem{BS}
Kristian Bredies and Hongpeng Sun.
\newblock Preconditioned douglas--rachford algorithms for tv- and tgv-regularized variational imaging problems.
\newblock {\em Journal of Mathematical Imaging and Vision}, 52(3):317--344, Jul 2015.

\bibitem{BSCC}
Kristian Bredies and Hongpeng Sun.
\newblock Preconditioned douglas--rachford splitting methods for convex-concave saddle-point problems.
\newblock {\em SIAM Journal on Numerical Analysis}, 53(1):421--444, 2015.

\bibitem{arias}
Luis~M. Brice\~{n}o Arias and Fernando Rold\'{a}n.
\newblock Split-douglas--rachford algorithm for composite monotone inclusions and split-admm.
\newblock {\em SIAM Journal on Optimization}, 31(4):2987--3013, 2021.

\bibitem{CP}
A.~Chambolle and T.~Pock.
\newblock A first-order primal-dual algorithm for convex problems with applications to imaging.
\newblock {\em Journal of Mathematical Imaging and Vision}, 40(1):120--145, 2011.

\bibitem{CERS}
Antonin Chambolle, Matthias~J. Ehrhardt, Peter Richt\'{a}rik, and Carola-Bibiane Sch\"{o}nlieb.
\newblock Stochastic primal-dual hybrid gradient algorithm with arbitrary sampling and imaging applications.
\newblock {\em SIAM Journal on Optimization}, 28(4):2783--2808, 2018.

\bibitem{libsvm}
Chih-Chung Chang and Chih-Jen Lin.
\newblock Libsvm: A library for support vector machines.
\newblock {\em ACM Trans. Intell. Syst. Technol.}, 2(3), may 2011.

\bibitem{CPE}
Patrick~L. Combettes and Jean-Christophe Pesquet.
\newblock Stochastic quasi-fejér block-coordinate fixed point iterations with random sweeping.
\newblock {\em SIAM Journal on Optimization}, 25(2):1221--1248, 2015.

\bibitem{CPE1}
Patrick~L. Combettes and Jean-Christophe Pesquet.
\newblock Stochastic quasi-fej{\'e}r block-coordinate fixed point iterations with random sweeping ii: mean-square and linear convergence.
\newblock {\em Mathematical Programming}, 174(1):433--451, Mar 2019.

\bibitem{Condat1}
Laurent Condat, Daichi Kitahara, Andr\'{e}s Contreras, and Akira Hirabayashi.
\newblock Proximal splitting algorithms for convex optimization: A tour of recent advances, with new twists.
\newblock {\em SIAM Review}, 65(2):375--435, 2023.

\bibitem{DR}
Jim Douglas and H.~H. Rachford.
\newblock On the numerical solution of heat conduction problems in two and three space variables.
\newblock {\em Transactions of the American Mathematical Society}, 82(2):421--439, 1956.

\bibitem{EP}
J.~Eckstein and D.~P. Bertsekas.
\newblock On the {Douglas--Rachford} splitting method and the proximal point algorithm for maximal monotone operators.
\newblock {\em Mathematical Programming}, 55(1-3):293--318, 1992.

\bibitem{YMV}
Yu.~M. Ermol’ev.
\newblock On the method of generalized stochastic gradients and quasi-fejér sequences.
\newblock {\em Cybernetics}, 5:208--220, 1969.

\bibitem{FercoqBianchi}
Olivier Fercoq and Pascal Bianchi.
\newblock A coordinate-descent primal-dual algorithm with large step size and possibly nonseparable functions.
\newblock {\em SIAM Journal on Optimization}, 29(1):100--134, 2019.

\bibitem{guohan}
Ke~Guo, Deren Han, and Xiaoming Yuan.
\newblock Convergence analysis of douglas--rachford splitting method for “strongly + weakly” convex programming.
\newblock {\em SIAM Journal on Numerical Analysis}, 55(4):1549--1577, 2017.

\bibitem{nips2003}
Isabelle Guyon, Steve Gunn, Asa Ben-Hur, and Gideon Dror.
\newblock Result analysis of the nips 2003 feature selection challenge.
\newblock In L.~Saul, Y.~Weiss, and L.~Bottou, editors, {\em Advances in Neural Information Processing Systems}, volume~17. MIT Press, 2004.

\bibitem{hanyan}
Deren Han, Yansheng Su, and Jiaxin Xie.
\newblock Randomized douglas–rachford methods for linear systems: Improved accuracy and efficiency.
\newblock {\em SIAM Journal on Optimization}, 34(1):1045--1070, 2024.

\bibitem{he1}
Bingsheng He and Xiaoming Yuan.
\newblock On the $\mathcal{O}(1/n)$ convergence rate of the douglas–rachford alternating direction method.
\newblock {\em SIAM Journal on Numerical Analysis}, 50(2):700--709, 2012.

\bibitem{he2}
Bingsheng He and Xiaoming Yuan.
\newblock On non-ergodic convergence rate of douglas--rachford alternating direction method of multipliers.
\newblock {\em Numerische Mathematik}, 130(3):567--577, Jul 2015.

\bibitem{he3}
Bingsheng He and Xiaoming Yuan.
\newblock On the convergence rate of douglas--rachford operator splitting method.
\newblock {\em Mathematical Programming}, 153(2):715--722, Nov 2015.

\bibitem{kl2}
Thorsten Hohage and Frank Werner.
\newblock Inverse problems with poisson data: statistical regularization theory, applications and algorithms.
\newblock {\em Inverse Problems}, 32(9):093001, jul 2016.

\bibitem{RGMGK}
Richard Huber, Georg Haberfehlner, Martin Holler, Gerald Kothleitner, and Kristian Bredies.
\newblock Total generalized variation regularization for multi-modal electron tomography.
\newblock {\em Nanoscale}, 11(12):5617--5632, February 2019.

\bibitem{LT}
M.~Ledoux and M.~Talagrand.
\newblock {\em Probability in Banach Spaces: Isoperimetry and Processes}.
\newblock Number~1 in Classics in Applied Mathematics. Springer-Verlag Berlin Heidelberg, 1991.

\bibitem{tkpong}
Guoyin Li and Ting~Kei Pong.
\newblock Douglas--rachford splitting for nonconvex optimization with application to nonconvex feasibility problems.
\newblock {\em Mathematical Programming}, 159(1):371--401, Sep 2016.

\bibitem{LM}
P.~L. Lions and B.~Mercier.
\newblock Splitting algorithms for the sum of two nonlinear operators.
\newblock {\em SIAM Journal on Numerical Analysis}, 16(6):964--979, 1979.

\bibitem{luke}
D.~Russell Luke and Anna-Lena Martins.
\newblock Convergence analysis of the relaxed douglas--rachford algorithm.
\newblock {\em SIAM Journal on Optimization}, 30(1):542--584, 2020.

\bibitem{NJL}
A.~Nemirovski, A.~Juditsky, G.~Lan, and A.~Shapiro.
\newblock Robust stochastic approximation approach to stochastic programming.
\newblock {\em SIAM Journal on Optimization}, 19(4):1574--1609, 2009.

\bibitem{Op1}
Zdzisław Opial.
\newblock {Weak convergence of the sequence of successive approximations for nonexpansive mappings}.
\newblock {\em Bulletin of the American Mathematical Society}, 73(4):591--597, 1967.

\bibitem{sample1}
Zheng Qu, Peter Richtarik, and Tong Zhang.
\newblock Quartz: Randomized dual coordinate ascent with arbitrary sampling.
\newblock In C.~Cortes, N.~Lawrence, D.~Lee, M.~Sugiyama, and R.~Garnett, editors, {\em Advances in Neural Information Processing Systems}, volume~28. Curran Associates, Inc., 2015.

\bibitem{sample2}
Peter Richtárik and Martin Takáč.
\newblock Parallel coordinate descent methods for big data optimization.
\newblock {\em Mathematical Programming}, 156:433--484, 2016.

\bibitem{RM2}
H.~Robbins and D.~Siegmund.
\newblock A convergence theorem for non negative almost supermartingales and some applications.
\newblock In Jagdish~S. Rustagi, editor, {\em Optimizing Methods in Statistics}, pages 233--257. Academic Press, 1971.

\bibitem{ShZ}
Shai Shalev-Shwartz and Tong Zhang.
\newblock Stochastic dual coordinate ascent methods for regularized loss.
\newblock {\em J. Mach. Learn. Res.}, 14(1):567–599, feb 2013.

\bibitem{SAN}
A.~N. Shiryaev.
\newblock {\em Probability-2, 3 edn.}
\newblock GTM, vol. 95, Springer, New York, 2019.

\bibitem{BF}
Benar~F. Svaiter.
\newblock A simplified proof of weak convergence in douglas–rachford method.
\newblock {\em Operations Research Letters}, 47(4):291--293, 2019.

\bibitem{patrinos}
Andreas Themelis and Panagiotis Patrinos.
\newblock Douglas--rachford splitting and admm for nonconvex optimization: Tight convergence results.
\newblock {\em SIAM Journal on Optimization}, 30(1):149--181, 2020.

\bibitem{Zhangxiao}
Yuchen Zhang and Lin Xiao.
\newblock Stochastic primal-dual coordinate method for regularized empirical risk minimization.
\newblock {\em Journal of Machine Learning Research}, 18(84):1--42, 2017.

\end{thebibliography}

\end{document}